\documentclass[a4paper,french,twoside,12pt]{smfart}

\usepackage{smfenum}
\usepackage{smfthm}
\usepackage{amssymb}
\usepackage{mathrsfs}
\usepackage{hyperref}
\input{xypic}

\numberwithin{equation}{section} 
\theoremstyle{plain}

\def\CC{\mathbb{C}}

\def\QQ{\mathbb{Q}}
\def\RR{\mathbb{R}}
\def\ZZ{\mathbb{Z}} 

\def\A{{\rm A}}
\def\B{{\rm B}}
\def\C{{\rm C}}
\def\D{{\rm D}}
\def\E{{\rm E}}
\def\F{{\rm F}}
\def\G{{\rm G}}
\def\H{{\rm H}}
\def\I{{\rm I}}
\def\J{{\rm J}}
\def\K{{\rm K}}
\def\L{\Lambda}
\def\M{{\rm M}}
\def\N{{\rm N}}
\def\P{{\rm P}}

\def\R{{\rm R}}

\def\U{{\rm U}}
\def\V{{\rm V}}
\def\W{{\rm W}}
\def\X{{\rm X}}

\def\AA{\mathfrak{A}}
\def\BB{\mathfrak{B}}
\def\HH{\mathfrak{H}}
\def\JJ{\mathfrak{J}}
\def\KK{\mathfrak{K}}
\def\LL{{\rm L}}
\def\MM{\mathfrak{m}}
\def\PP{\mathfrak{a}}
\def\PB{\mathfrak{b}}
\def\SS{{\rm S}}

\def\Cc{\mathscr{C}}

\def\Ff{\mathscr{F}}
\def\Gg{\mathscr{G}}

\def\Ll{\mathscr{L}}

\def\Oo{\mathscr{O}}
\def\Pp{\mathscr{P}}

\def\Uu{\mathscr{U}}

\def\a{\alpha} 
\def\b{\beta}
\def\d{\delta}
\def\e{{\rm e}}
\def\f{\rightarrow}
\def\g{\gamma}

\def\ii{{\iota}}

\def\k{\kappa}
\def\l{\lambda}

\def\n{\eta}
\def\o{\mathfrak{o}}
\def\o{\mathscr{O}}
\def\p{\mathfrak{p}}
\def\r{\rho}
\def\s{\sigma}
\def\t{\theta}
\def\v{\upsilon}
\def\w{\varpi}

\def\ie{c'est-\`a-dire }

\def\>{\geqslant}
\def\<{\leqslant}

\def\Hom{{\rm Hom}}
\def\End{{\rm End}}
\def\Aut{{\rm Aut}}
\def\Mat{{\rm M}}
\def\GL{{\rm GL}}
\def\Gal{{\rm Gal}}

\def\tr{{\rm tr}}

\def\Ind{{\rm Ind}}
\def\mult#1{{#1}^{\times}}

\def\Seq{{\bf S}}
\def\Seq{\Ll}

\def\rr{r}
\def\ss{s}
\def\Ga{\Gamma}

\def\epsilo{\epsilon}
\def\epsilo{\rho}
\def\Om{\Omega}
\def\pt{\circ}

\def\audessus{{coh\'erente avec}}
\def\unif{\w}

\author{V. S\'echerre}
\address{Institut de Math\'ematiques de Luminy\\
CNRS UMR $6206$\\
Universit\'e de la M\'editerran\'ee, 163 avenue de Luminy\\
$13288$ Marseille Cedex $09$\\
France}
\email{secherre@iml.univ-mrs.fr}

\author{S. Stevens}
\address{School of Mathematics\\
University of East Anglia\\
Norwich NR4 7TJ\\
United Kingdom}
\email{Shaun.Stevens@uea.ac.uk}

\title[Repr\'esentations supercuspidales de $\GL_m(\D)$]
{Repr\'esentations lisses de $\GL_m(\D)$\\
IV~: repr\'esentations supercuspidales}

\begin{abstract}
Soit $\F$ un corps commutatif localement compact non archi\-m\'edien,
et soit $\D$ une alg\`ebre \`a division de centre $\F$.
Nous prouvons que toute repr\'esentation irr\'eductible super\-cuspidale
du groupe $\GL_m(\D)$, de niveau non nul, est l'induite compacte d'une
repr\'esentation d'un sous-groupe ouvert compact modulo le centre de 
$\GL_m(\D)$.
Plus pr\'ecis\'ement, nous prouvons que de telles repr\'esentations
contiennent un type simple maxi\-mal au sens de \cite{VS3}.
\end{abstract}

\begin{altabstract}
Let $\F$ be a non Archimedean locally compact field and let $\D$ be a 
central $\F$-division algebra.
We prove that any positive level supercuspidal irreducible
representation of the group $\GL_m(\D)$ is compactly induced from a 
representation of a compact mod center open subgroup of $\GL_m(\D)$.
More precisely, we prove that such representations contain a maximal
simple type in the sense of \cite{VS3}.
\end{altabstract}

\thanks
{Ce travail a b\'en\'efici\'e d'un financement de la part de 
l'EPSRC (grant GR/T21714/01).}

\begin{document}

\maketitle
\setcounter{tocdepth}{2}

\section*{Introduction}

Soit $\F$ un corps commutatif localement compact non archim\'edien, 
et soit $\G$ une forme in\-t\'e\-rieu\-re de $\GL_n(\F)$, $n\>1$.
C'est un groupe de la forme $\GL_m(\D)$, o\`u $\D$ est une
$\F$-alg\`ebre \`a division, de dimension $d^2$ sur son centre $\F$,
et o\`u $n=md$.
Cet article, qui fait suite au travail entrepris par le premier auteur
dans \cite{VS1,VS2,VS3}, met un terme \`a la classification des blocs
simples de la cat\'egorie des repr\'esentations lisses complexes de
$\G$ au moyen de la th\'eorie des types de Bushnell et Kutzko. 

Notre r\'esultat principal peut \^etre formul\'e ainsi : si $\rho$ est
une repr\'esentation irr\'eductible supercuspidale de $\G$, alors il
existe un type pour la classe inertielle de $\rho$.
En d'autres termes, nous prouvons qu'il existe un sous-groupe ouvert
compact $\J$ de $\G$ et une repr\'esentation irr\'eductible $\l$
de $\J$ telle que les repr\'esentation ir\-r\'e\-duc\-ti\-bles de $\G$
dont la res\-tric\-tion \`a $\J$ contient $\l$ sont exac\-tement celles
qui sont \'equivalentes \`a $\rho\otimes\chi$ pour un caract\`ere non
ramifi\'e $\chi$ de $\G$.
Plus pr\'ecis\'ement, nous prouvons qu'on peut choisir pour $(\J,\l)$
un type simple maximal au sens de \cite{VS3}.

Le probl\`eme de la classification des repr\'esentations lisses
complexes de $\G$ par la th\'eorie des types a d\'ej\`a \'et\'e
abord\'e par plusieurs auteurs.
Bien entendu, il faut mentionner en premier lieu les travaux
fondateurs de Bushnell et Kutzko \cite{BK,BK2} concernant le groupe
d\'eploy\'e $\GL_n(\F)$, qui ont 
donn\'e le ton \`a tous les trauvaux ult\'erieurs sur le sujet.
Ensuite, les premiers travaux concernant les formes int\'erieures non
d\'eploy\'ees de $\GL_n(\F)$ sont ceux de E.-W.\ Zink \cite{Zi5} et
de Broussous \cite{Br2}~: tous deux donnent une classification des
repr\'esentations de $\GL_1(\D)$, le premier
lorsque $\F$ est de caract\'eristique nulle, le se\-cond sans
restriction sur la caract\'eristique.
Dans \cite{GSZ}, Grabitz, Silberger et Zink traitent le cas tr\`es
particulier du niveau z\'ero, \ie des repr\'esentations
irr\'eductibles de $\GL_m(\D)$ poss\'edant un vecteur non nul
invariant par le sous-groupe $1+\Mat_m(\p_\D)$, o\`u $\p_\D$ d\'esigne
l'id\'eal maximal de l'anneau des entiers de $\D$.
Concernant le cas g\'en\'eral, \ie les repr\'esentations
irr\'eductibles de $\GL_m(\D)$ de niveau quelconque, on trouve un
certain nombre de r\'esultats dans les travaux de Broussous
\cite{Br3,Br1,Br4}, Broussous-Grabitz \cite{BG} et Grabitz
\cite{Grabitz,Grabitz2}.

Ce travail | les articles \cite{VS1,VS2,VS3} auxquels vient s'ajouter
le pr\'esent article | suit la m\'ethode g\'en\'erale de construction
de types \'elabor\'ee par Bushnell et Kutzko dans \cite{BK} et
am\'elior\'ee dans \cite{BK1} par la th\'eorie des paires couvrantes. 
D\'ecrivons-en bri\`evement l'organisation.
On fixe une fois pour toutes une strate simple $[\AA,n,0,\b]$ de
la $\F$-alg\`ebre $\Mat_m(\D)$ 
({\it cf.} D\'e\-fi\-nition \ref{stratepure}). 
Rappelons simplement ici que $\b$ est un \'el\'ement de $\Mat_m(\D)$
tel que la $\F$-alg\`ebre $\E=\F[\b]$ soit un corps et que $\AA$ est un
$\o_{\F}$-ordre h\'er\'editaire de $\Mat_m(\D)$ normalis\'e par
$\mult\E$. 
\begin{enumerate}
\item[(i)]
Dans une premi\`ere \'etape (\cite{VS1}), on associe \`a la strate
simple $[\AA,n,0,\b]$ un ensemble fini $\Cc(\b,\AA)$ de 
{\it caract\`eres simples}.
Il s'agit de caract\`eres d\'efinis sur un certain sous-groupe ouvert
compact $\H^1=\H^1(\b,\AA)$ de $\G$ jouissant de remarquables
propri\'et\'es de fonctorialit\'e connues sous le nom de
propri\'et\'es de {\it transfert}.  
Plus pr\'ecis\'ement, si $[\AA',n',0,\b]$ est n'importe quelle strate
simple d'une $\F$-alg\`ebre centrale simple dans laquelle est plong\'e
$\E$, il existe une bijection canonique de $\Cc(\b,\AA)$ sur
$\Cc(\b,\AA')$.
\item[(ii)]
Dans une seconde \'etape (\cite{VS2}), on construit, pour chaque 
caract\`ere simple $\t\in\Cc(\b,\AA)$, une famille finie de 
$\beta$-{\it extensions}. 
Il s'agit de repr\'esentations irr\'eductibles d'un sous-groupe ouvert
compact $\J=\J(\b,\AA)$ de $\G$ dont la res\-tric\-tion \`a
$\H^1(\b,\AA)$ contient $\t$ et | surtout | dont l'entrelacement est
le m\^eme que celui de $\t$.
\item[(iii)]
Dans une troisi\`eme \'etape (\cite{VS3}), lorsque $\AA$ est un ordre
prin\-ci\-pal, on cons\-truit pour chaque $\b$-extension $\k$
d'un caract\`ere simple $\t\in\Cc(\b,\AA)$ une fa\-mil\-le finie de 
{\it types simples}. 
Ce sont des repr\'esentations irr\'eductibles de $\J(\b,\AA)$ de la
forme $\k\otimes\s$, o\`u $\s$ est l'inflation \`a $\J(\b,\AA)$ d'une
repr\'esentation irr\'eductible super\-cus\-pi\-da\-le du groupe
r\'eductif fini~:
\begin{equation*}
\label{grfinired}
\J(\b,\AA)/\J^1(\b,\AA)\simeq\GL_s(k)^r,
\end{equation*}
o\`u $r,s$ sont des entiers $\>1$ tels que le produit $rs$ divise $m$, 
o\`u $k$ est une extension finie du corps r\'esiduel de $\E$ et o\`u
$\s$ est de la forme $\s_0^{\otimes r}$, avec $\s_0$ une
repr\'esentation irr\'eductible super\-cus\-pi\-da\-le de $\GL_{s}(k)$. 
Chaque type simple $(\J(\b,\AA),\l)$ ainsi construit est un type pour
une classe inertielle simple de $\G$ de la forme
$[\G_0^r,\rho_0^{\otimes r}]_\G$, o\`u $\rho_0$ est une
repr\'esentation irr\'eductible supercuspidale du groupe
$\G_0=\GL_{m/r}(\D)$.
L'alg\`ebre de Hecke de $\G$ relative \`a $\l$ est une alg\`ebre de
Hecke affine. 
Plus pr\'ecis\'ement, elle est isomorphe \`a l'alg\`ebre
de Hecke-Iwahori de $\GL_r(\K)$, o\`u $\K$ est une extension non
ramifi\'ee de $\E$ dont le corps r\'esiduel est une extension de
degr\'e $s$ de $k$. 
\item[(iv)]
La quatri\`eme \'etape est celle qui occupe le pr\'esent article.
L'objectif en est l'exhaustion des re\-pr\'e\-sen\-ta\-tions
irr\'eductibles supercuspidales par les types simples.
Plus pr\'ecis\'ement, nous prou\-vons le r\'esultat suivant 
({\it cf.} Th\'eor\`eme \ref{ExhaustionSupercuspidale})~:

\medskip

\noindent\textsc{Th\'eor\`eme :} 
Soit $\rho$ une repr\'esentation irr\'eductible supercuspidale 
de niveau non nul de $\G$.
Il existe un type simple maximal $(\J(\b,\AA),\l)$, au sens de
\cite{VS3}, tel que la restriction de $\rho$ \`a $\J(\b,\AA)$ 
contienne $\l$.

\medskip

On en d\'eduit imm\'ediatement, \`a partir de 
\cite[Th\'eor\`eme 5.6]{VS3}, que si 
$\mathfrak{s}=[\G_0^r,\rho_0^{\otimes r}]_\G$ est une classe
inertielle simple de niveau non nul de $\G$, \ie que $r$ est un
diviseur de $m$ et $\rho_0$ une repr\'esentation irr\'eductible
supercuspidale de niveau non nul du groupe $\G_0=\GL_{m/r}(\D)$, 
il existe un type simple $(\J(\b,\AA),\l)$ qui est un type pour 
$\mathfrak{s}$.
La structure de l'alg\`ebre de Hecke de $\G$ relative \`a
$(\J(\b,\AA),\l)$ est donn\'ee par \cite[Th\'eor\`eme 4.6]{VS3}.
Comme mentionn\'e \`a l'\'etape (iii), elle est isomorphe \`a une
alg\`ebre de Hecke-Iwahori.
\end{enumerate}

Nous passons maintenant \`a la description des m\'ethodes
utilis\'ees. 
Notre premi\`ere t\^ache | qui est aussi la plus difficile | est de
s\'eparer les re\-pr\'e\-sen\-ta\-tions irr\'eductibles de niveau non
nul de $\G$ en trois cat\'egories, comme suit
({\it cf.} Th\'eor\`eme \ref{??})~: une repr\'esentation
irr\'eductible de niveau non nul de $\G$ contient ou bien une strate
scind\'ee, ou bien un caract\`ere scind\'e, ou bien un
ca\-rac\-t\`e\-re simple $\t\in\Cc(\b,\AA)$ pour une strate simple 
$[\AA,n,0,\b]$ de $\A$.
(On renvoie au \S\ref{SuperSection} pour les d\'efinitions de strate
scind\'ee et de caract\`ere scin\-d\'e | appel\'e {\it type scind\'e}
dans \cite{BK}.)
Ce travail est d\'ej\`a amorc\'e par Broussous \cite{Br4}, qui prouve
qu'une repr\'esentation irr\'eductible de niveau non nul $\pi$ de $\G$ 
ne contenant pas de strate scind\'ee contient un ca\-rac\-t\`e\-re
simple de niveau $m\>0$, \ie la restriction d'un ca\-rac\-t\`e\-re
simple $\t\in\Cc(\b,\AA)$ \`a un sous-groupe $\H^{m+1}(\b,\AA)$
\'eventuellement plus petit que $\H^1(\b,\AA)$.
Si $m\>1$, on cherche \`a construire, en proc\'edant par raffinement
comme dans \cite{BK}, un caract\`ere simple $\t'\in\Cc(\b,\AA')$
relatif \`a une autre strate simple $[\AA',n',0,\b]$ de $\A$, contenu
dans $\pi$ \`a un niveau (normalis\'e) strictement moindre.
Mais, contrairement \`a ce qui se passe dans le cas d\'eploy\'e, 
il n'est pas possible de passer de $\t$ \`a $\t'$ en une seule \'etape
| techniquement, c'est \cite[Proposition 1.2.8]{BK} qui fait
d\'efaut~: voir \cite[\S1.5.3]{VS3}.
En d\'ecoupant le segment $[\AA,\AA']$ dans l'immeuble de Bruhat-Tits
de $\G$ en morceaux suffisamment petits 
({\it cf.} \S\ref{LeSoleilBrilleAlaPlage}, Hypoth\`ese (H)), 
on obtient le r\'esultat crucial \ref{Transfer}, clef du processus de
raffinement et analogue de \cite[Lemma 5.4]{St4}.
Tout ceci est assez tech\-ni\-que et n\'ecessite~:
\begin{enumerate}
\item[(i)]
d'employer le langage des suites de r\'eseaux, plus g\'en\'eral que
celui des ordres h\'er\'editaires et qui permet une description des
points rationnels de l'immeuble de Bruhat-Tits de $\G$ 
({\it cf.} \cite{Br3}) ;
\item[(ii)]
de d\'efinir les caract\`eres simples relatifs \`a une suite de
r\'eseaux et d'\'e\-ten\-dre le transfert des caract\`eres simples 
\`a ce cadre : c'est l'objet de la section \ref{CSsuites} de cet
article ; 
\item[(iii)]
de g\'en\'eraliser la notion de strate d\'eriv\'ee 
({\it cf.} D\'efinition \ref{StrateDerivee!}), conduisant 
elle-m\^eme \`a celle de caract\`ere scind\'e 
({\it cf.} D\'efinition \ref{CarScinde}).
\end{enumerate}

Mentionnons un point important~: l'introduction des caract\`eres
simples relatifs \`a une suite de r\'eseaux 
a ceci de g\^enant que les niveaux normalis\'es deviennent des nombres
rationnels dont le d\'enominateur n'est pas born\'e | ce qui
risquerait d'emp\^echer le processus de raffinement de terminer.
Il est donc indispensable de v\'erifier que l'on passe de $\t$ \`a un
caract\`ere simple $\t'$ qui, lui, est relatif \`a un ordre
h\'er\'editaire, m\^eme si les \'etapes interm\'ediaires peuvent
mettre en jeu des caract\`eres simples relatifs \`a des suites de
r\'eseaux quelconques.
Ceci justifie le th\'eor\`eme \ref{Edescent}, et tous les r\'esultats
pr\'eparatoires des \S\S\ref{Diotime}--\ref{Parmenide}.

Remarquons \'egalement que la distinction entre strate scind\'ee et
caract\`ere scind\'e est assez superficielle~: c'est \`a peu pr\`es 
la m\^eme que celle que l'on fait entre types simples de niveau $0$ 
et de niveau $>0$. 

Notre seconde t\^ache, qui occupe la section \ref{AymeJacquet}, est de
montrer qu'une re\-pr\'e\-sen\-ta\-tion irr\'eductible de niveau non
nul de $\G$ contenant une strate scind\'ee ou un caract\`ere scind\'e
a un module de Jacquet non trivial.
Ceci implique automatiquement que toute repr\'esentation
irr\'eductible supercuspidale de niveau non nul de $\G$ contient un 
ca\-rac\-t\`e\-re simple $\t\in\Cc(\b,\AA)$ pour une strate simple
$[\AA,n,0,\b]$ de $\A$.
Pour aboutir au th\'eor\`eme \ref{ExhaustionSupercuspidale}, il reste
alors \`a d\'ecrire le passage de $\H^1$ \`a $\J$, ce qui fait l'objet
de la section \ref{YoupiNiveauZero}.
Une repr\'esentation irr\'eductible supercuspidale de niveau non nul
de $\G$ contient {\it a priori} une repr\'esentation irr\'eductible 
de $\J$ de la forme $\vartheta=\k\otimes\s$, o\`u $\s$ est l'inflation
\`a $\J$ d'une repr\'esentation irr\'eductible de $\J/\J^1$.
En s'inspirant de \cite[\S1]{GSZ}, et en utilisant la notion de 
{\it coh\'erence} d\'ej\`a largement d\'evelopp\'ee dans
\cite{VS2,VS3}, on montre que $\vartheta$ est un type simple maximal, 
\ie que $\AA$ est principal, que $\s$ est cuspidale de la forme
$\s_0^{\otimes r}$ et que $\AA\cap\B$ est un ordre maximal de $\B$.

Nous terminons cette introduction en mentionnant deux probl\`emes qui
restent \`a traiter au sujet des repr\'esentations lisses de
$\G$~:  
(i) la cons\-truction d'endoclasses et (ii) la construction de types
pour n'importe quelle classe inertielle de $\G$.
Le premier revient essentiellement \`a prouver un r\'esultat de type
{\it ``entrelacement implique conjugaison''} pour les caract\`eres
simples~: si deux caract\`eres simples $\t_1\in\Cc(\AA,m_1,\b_1)$ et 
$\t_2\in\Cc(\AA,m_2,\b_2)$ s'entrelacent dans $\G$, alors ils sont
conjugu\'es sous $\G$.
Le second revient, de fa\c con analogue \`a \cite{BK2}, \`a construire
des {\it types semi-simples} pour $\G$.

\section*{Notations et conventions}

Soit $\F$ un corps commutatif localement compact non archim\'edien.
Toutes les $\F$-alg\`ebres sont suppos\'ees unitaires et de dimension
finie.
Par $\F$-\emph{alg\`ebre \`a division} on entend $\F$-alg\`ebre 
centrale dont l'anneau sous-jacent est un corps (pas n\'ecessairement
commutatif).

Si $\K$ est une extension finie de $\F$, ou plus g\'en\'eralement 
une alg\`ebre \`a division sur une extension finie de $\F$, 
on note $\o_\K$ son anneau d'entiers, $\p_\K$ son id\'eal maximal 
et $k_\K$ son corps r\'esiduel.

Si $\A$ est une alg\`ebre centrale simple sur une extension finie 
$\K$ de $\F$, on note ${\rm N}_{\A/\K}$ (resp. $\tr_{\A/\K}$) la 
norme (resp. la trace) r\'eduite de $\A$ sur $\K$.

Si $u$ est un nombre r\'eel, on note $\lceil u\rceil$ le plus petit
entier $\>u$ et $\lfloor u\rfloor$ le plus grand entier $\<u$, 
\ie la partie enti\`ere de $u$.

Un {\it caract\`ere} d'un groupe topologique $\G$ est un 
homomorphisme continu de $\G$ dans le groupe multiplicatif $\mult\CC$
du corps des nombres complexes.

Toutes les repr\'esentations sont suppos\'ees lisses et \`a coefficients
complexes. 


\section{Pr\'eliminaires}
\label{preliminaires}

Dans cette section, on rappelle le langage des strates dans une
$\F$-alg\`ebre centrale simple.
Pour plus de d\'etails, on renvoie le lecteur \`a 
\cite{Br4,BK,BK2,VS1}.

\subsection{}
\label{Lapinescence}

Soit $\A$ une $\F$-alg\`ebre centrale simple, et soit $\V$ un $\A$-module
\`a gauche simple.
L'alg\`ebre $\End_{\A}(\V)$ est une $\F$-alg\`ebre \`a division
dont l'alg\`ebre oppos\'ee est not\'ee $\D$. 
Aussi $\V$ est-il un $\D$-espace vectoriel \`a droite, 
et on a un isomorphisme canonique de $\F$-alg\`ebres entre
$\A$ et $\End_{\D}(\V)$.

\begin{defi}
Une {\it $\o_{\D}$-suite de r\'eseaux} de $\V$ est une suite 
d\'ecroissante $\L=(\L_k)_{k\in\ZZ}$ de $\o_\D$-r\'eseaux de 
$\V$ pour laquelle il existe un entier $e\>1$ tel que, pour 
tout $k\in\ZZ$, on ait $\L_{k+e}=\L_k\p_\D$.  
Cet entier $e$, unique, est appel\'e la {\it p\'eriode} de 
$\L$ sur $\o_\D$ et not\'e $e(\L|\o_\D)$.   
\end{defi}

Une $\o_\D$-suite de r\'eseaux est dite {\it stricte} si elle est
strictement d\'ecroissante.
On emploie aussi le terme de {\it cha\^\i ne} de r\'eseaux
pour suite de r\'eseaux stricte.
On note $\Seq(\V,\Oo_\D)$ l'ensemble des $\o_\D$-suites 
de r\'eseaux de $\V$.

Si $\W$ est un sous-$\D$-espace vectoriel de $\V$, l'application 
$k\mapsto\L_{k}\cap\W$ est une $\o_\D$-suite de r\'eseaux de $\W$
de m\^eme p\'eriode que $\L$.
On la note $\L\cap\W$.

Si $a$ est un entier $\>1$ et si $\L\in\Seq(\V,\Oo_\D)$, on note $a\L$
la suite $k\mapsto\L_{\lceil k/a\rceil }$.

\begin{rema}
\label{convGeneve}
Il est commode de prolonger les suites de r\'eseaux \`a $\RR$
tout entier, en posant $\L_x=\L_{\lceil x\rceil}$ pour $x\in\RR$.
De cette fa\c con, $\L$ devient une fonction d\'ecroissante \`a
valeurs dans l'ensemble des $\o_\D$-r\'eseaux de $\V$, continue \`a
gauche pour la topologie discr\`ete sur l'espace d'arriv\'ee.
Ses points de discontinuit\'e sont entiers.
L'application $x\mapsto\L(xe(\L|\o_\F))$, d\'efinie sur $\RR$, 
est une {\it $\o_\D$-fonc\-tion de r\'eseaux} au sens de \cite{Br3}.
\end{rema}

\subsection{}

\`A toute suite de r\'eseaux $\L\in\Seq(\V,\Oo_\D)$ on associe une
suite de r\'eseaux $\PP(\L)\in\Seq(\A,\Oo_\F)$ d\'efinie par~:
\begin{equation*}
\label{filsuite}
\PP_k(\L)=\{a\in\A\ |\ a\L_{l}\subset\L_{l+k},\ l\in\ZZ\}, 
\quad k\in\ZZ.
\end{equation*}
Elle caract\'erise la classe de translation de $\L$.
Deux $\o_\F$-r\'eseaux de cette suite sont d'une importance
particuli\`ere : $\AA(\L)=\PP_0(\L)$ est un ordre
h\'e\-r\'e\-di\-tai\-re de $\A$, et son radical de Jacobson
est $\mathfrak{P}(\L)=\PP_1(\L)$.
Ils ne d\'ependent que de l'en\-sem\-ble $\{\L_k\ |\ k\in\ZZ\}$.

On note $\KK(\L)$ le normalisateur de $\L$ dans $\mult\A$.
Si $g\in\KK(\L)$, on note $\v_{\L}(g)$ l'entier $n\in\ZZ$
d\'efini par $g(\L_{k})=\L_{k+n}$.
L'application $\v_{\L}$ est un morphisme de groupes de $\KK(\L)$
dans $\ZZ$, dont le noyau, not\'e $\U(\L)$, 
est le groupe multiplicatif de l'ordre h\'er\'editaire $\AA(\L)$.
On pose $\U_0(\L)=\U(\L)$ et, pour $k\>1$, on pose 
$\U_k(\L)=1+\PP_k(\L)$.

\subsection{}
\label{PureteDeMoeurs}

Soit $\E$ une extension finie de $\F$ 
telle qu'il existe un homomorphisme de $\F$-alg\`ebres $\ii:\E\f\A$,
qu'on suppose fix\'e dans la suite.
On identifie ainsi $\E$ \`a la sous-$\F$-alg\`ebre $\ii\E$ de $\A$.

\begin{defi}
Une suite $\L\in\Seq(\V,\Oo_\D)$ est dite 
{\it $\E$-pure} si elle est nor\-ma\-lis\'ee par $\mult\E$, 
\ie si c'est une $\o_\E$-suite de r\'eseaux de $\V$.
\end{defi}

Le commutant de $\E$ dans $\A$, qu'on note $\B$, est une 
$\E$-al\-g\`ebre centrale simple.
On fixe un $\B$-module \`a gauche simple $\V_\E$ et on note 
$\D_{\E}$ l'alg\`ebre oppos\'ee \`a $\End_\B(\V_{\E})$.
Si $\Ga$ est une $\Oo_{\D_\E}$-suite de r\'eseaux de $\V_\E$, 
on note $\PB(\Ga)$ la $\Oo_{\E}$-suite de r\'eseaux de $\B$ qu'elle
d\'efinit. 

Le th\'eor\`eme suivant traduit en termes de suites de r\'eseaux 
la cor\-respondance \cite{Br3} entre $\Oo_\D$-fonctions de r\'eseaux
de $\V$ invariantes par $\mult\E$ et $\Oo_{\D_\E}$-fonctions de
r\'eseaux de $\V_\E$. 
Il \'etend aux suites $\E$-pures quelconques la correspondance
d\'efinie par \cite[Th\'eor\`eme 1.3]{Br1}. 

\begin{theo}
\label{CNNBr1}
Soit $\L\in\Seq(\V,\Oo_\D)$ une suite $\E$-pure.
Il exis\-te une $\Oo_{\D_\E}$-suite $\Ga$ de r\'eseaux de $\V_\E$,
unique \`a translation pr\`es, telle que~:
\begin{equation*}
\label{Edescen}
\PP_{k}(\L)\cap\B=\PB_{k}(\Ga), \quad k\in\ZZ.
\end{equation*}
Le normalisateur de $\Ga$ dans $\mult\B$ est \'egal \`a
$\KK(\L)\cap\mult\B$.
\end{theo}

\begin{proof}
On se r\'ef\`ere \`a \cite{Br3}, qui utilise le langage des fonctions 
de r\'eseaux.
On d\'esigne 
par $d$ le degr\'e r\'eduit de $\D$ sur $\F$,
par $e$ l'indice de ra\-mi\-fication de $\E/\F$ et 
par $r$ la p\'eriode de $\L$ sur $\Oo_\D$.
\`A la suite $\L$ correspond la $\Oo_{\D}$-fonction 
$\Ff:x\mapsto\L(rdx)$ 
({\it cf.} Remarque \ref{convGeneve}).
Il lui correspond une $\Oo_{\D_\E}$-fonction $\Gg$ de r\'eseaux de
$\V_\E$, unique \`a translation pr\`es, telle que~:
\begin{equation}
\label{LeSoleilBrille}
\PP_x(\Ff)\cap\B=\PB_{ex}(\Gg), \quad x\in\RR
\end{equation}
({\it cf.} \cite[Theorem II.1.1]{Br3}).
On note $d_\E$ le degr\'e r\'eduit de $\D_\E$ sur $\E$.
C'est un diviseur de $d$, et le quotient de $d$ par $d_\E$ est \'egal
au pgcd de $d$ et $[\E:\F]$ (voir par exemple \cite{Zi1}). 
Puisque $\L$ est $\E$-pure, l'entier $e$ divise $rd$, donc $ed_\E$
divise $rd$. 
On en d\'eduit que $\Ga:k\mapsto\Gg(ke/rd)$ est une
$\Oo_{\D_\E}$-suite de r\'eseaux de $\V_\E$ qui satisfait \`a la
condition voulue.
L'unicit\'e est imm\'ediate.
\end{proof}

\begin{rema}
En particulier, l'intersection $\BB(\L)=\AA(\L)\cap\B$ est un or\-dre
h\'er\'editaire de $\B$.
\end{rema}

\subsection{}

L'inconv\'enient de la correspondance $\L\mapsto\Ga$ d\'efinie par 
le th\'eor\`eme \ref{CNNBr1} est qu'elle n'est pas normalis\'ee, ce
qui fait que son image ne contient pas, en g\'en\'eral, toutes les
$\Oo_{\D_\E}$-cha\^\i nes.

\begin{exem}
\label{ExCanotSauvetage}
Soit $\A=\Mat_2(\D)$ et soit $\E/\F$ une extension quadratique non
ramifi\'ee incluse dans $\D$.
Une cha\^\i ne $\E$-pure de p\'eriode $2$ a pour image une suite non
stricte (de p\'eriode $4$), une cha\^\i ne $\E$-pure de p\'eriode $1$
a pour image une suite de p\'eriode $2$, et une suite $\E$-pure non
stricte a pour image une suite non stricte.
Donc, si $\Ga$ est une $\Oo_{\D_\E}$-cha\^\i ne de p\'eriode $1$, il ne
lui correspond aucune $\Oo_{\D}$-suite $\E$-pure $\L$.
\end{exem}

Si l'on veut r\'ecup\'erer toutes les $\Oo_{\D_\E}$-cha\^\i nes, 
il est n\'ecessaire de rajouter un facteur de normalisation.
Le r\'esultat suivant compl\`ete \cite[Proposition II.5.4]{Br3}.

\begin{theo}
\label{Edescent}
Soit $\Ga$ une $\Oo_{\D_\E}$-suite stricte de $\V_\E$.
Il exis\-te un uni\-que entier $\epsilo\>1$ et 
une suite $\L\in\Seq(\V,\Oo_\D)$ stricte et $\E$-pure, 
unique \`a trans\-la\-tion pr\`es, tels que~: 
\begin{equation*}
\PP_{k}(\L)\cap\B=\PB_{{k}/{\epsilo}}(\Ga), \quad k\in\ZZ.
\end{equation*}
Le normalisateur de $\Ga$ dans $\mult\B$ est \'egal \`a
$\KK(\L)\cap\mult\B$. 
\end{theo}

\begin{rema}
Autrement dit, il existe un unique entier $\epsilo$ tel que
$\epsilo\Ga$ soit dans l'image de la correspondance d\'efinie 
par le th\'eor\`eme \ref{CNNBr1} et ait une cha\^\i ne pour 
ant\'ec\'edent.
\end{rema}

\begin{proof}
On reprend en partie les notations de la preuve pr\'ec\'edente.
On note $s$ le pgcd de $d$ et du degr\'e r\'esiduel de $\E/\F$
et $r$ la p\'eriode de $\Ga$ sur $\Oo_{\E}$.
\`A la suite $\Ga$ correspond la $\Oo_{\D_\E}$-fonction 
$\Gg:x\mapsto\Ga(rx)$.
L'ensemble de ses points de discontinuit\'e est $r^{-1}\ZZ$.
Il lui correspond une $\Oo_\D$-fonction $\E$-pure $\Ff$ de r\'eseaux
de $\V$, unique \`a translation pr\`es, v\'erifiant
(\ref{LeSoleilBrille}).
Plus pr\'ecis\'ement, la fonction $\Ff$ est d\'efinie en 
\cite[Lemma II.3.1]{Br3}, o\`u l'on voit que ses points de
discontinuit\'e sont les nombres r\'eels de la forme~: 
\begin{equation*}
\frac{a}{d}+\frac{b}{er}, \quad 0\<a\<s-1\text{ et } b\in\ZZ.
\end{equation*}
On va voir qu'ils forment un id\'eal fractionnaire de $\ZZ$.
Le quotient de $d$ par $d_\E$ divise $se$. 
On en d\'eduit que l'ensemble des points de discontinuit\'e de $\Ff$
est~:
\begin{equation*}
\frac{1}{d}\ZZ+\frac{1}{er}\ZZ=
\frac{(d,er)}{der}\ZZ,
\end{equation*}
o\`u $(a,b)$ d\'esigne le pgcd de deux entiers $\>1$.
Si on pose $\epsilo=d/(d,er)$,
l'application $\L:k\mapsto\Ff(k/\epsilo er)$ est une
$\Oo_\D$-suite de r\'eseaux stricte qui r\'epond \`a la question.

Il reste \`a prouver l'unicit\'e.
Si le couple $(\epsilo',\L')$ est une autre solution, alors
$\PP_k(\epsilo\L')\cap\B=\PP_k(\epsilo'\L)\cap\B$ pour $k\in\ZZ$.
Les suites $\epsilo'\L$ et $\epsilo\L'$ d\'efinissent la m\^eme
fonction de r\'eseaux $\Ff$ et ont la m\^eme p\'eriode : elles sont
donc dans la m\^eme classe de translation. 
Puisque $\L$ et $\L'$ sont toutes les deux strictes, on en d\'eduit
que $\epsilo=\epsilo'$, puis que $\L$ et $\L'$ sont dans la m\^eme
classe de translation.
\end{proof}

\begin{exem}
\begin{itemize}
\item[(i)]
On reprend l'exemple \ref{ExCanotSauvetage}.
Si $\Ga$ est une cha\^\i ne de p\'eriode $2$, alors $\L$ est de
p\'eriode $1$ et $\epsilo=1$.
Si $\Ga$ est une cha\^\i ne de p\'eriode $1$, alors $\L$ est de
p\'eriode $1$ et $\epsilo=2$.
\item[(ii)]
Dans le cas d\'eploy\'e, \ie lorsque $\D=\F$, la cha\^\i ne $\L$ 
est simplement $\Ga$ vue comme une $\o_\F$-cha\^\i ne et on a 
$\epsilo=1$.
\end{itemize}
\end{exem}

\begin{rema}
\label{Tetragrammaton}
L'entier $\epsilo$ est \'egal au rapport de $e(\L|\Oo_\F)$ sur 
$e(\Ga|\Oo_\F)$.
\end{rema}

\subsection{}
\label{KarlKrauss}

Soit :
\begin{equation}
\label{PipDickens}
\V=\V^1\oplus\ldots\oplus\V^l
\end{equation}
une d\'ecomposition de $\V$ en une somme directe de $l\>1$ 
sous-$\D$-espaces vectoriels.
Pour chaque $1\<i\<l$, on note $\e^i$ le projecteur sur $\V^i$
parall\`element \`a $\bigoplus_{j\neq i}\V^j$.
On pose $\A^{ij}=\e^i\A\e^j=\Hom_{\D}(\V^j,\V^i)$
et $\A^i=\A^{ii}=\End_{\D}(\V^i)$.
On note :
\begin{eqnarray*}
\M&=&\Aut_{\D}(\V^1)\times\ldots\times\Aut_{\D}(\V^l)
\end{eqnarray*}
le stabilisateur de la d\'ecomposition (\ref{PipDickens}) dans
$\mult\A$.
C'est un sous-groupe de Levi de $\mult\A$.
Soit $\P$ un sous-groupe parabolique de $\mult\A$ de facteur de
Levi $\M$, et \'ecrivons $\P=\M\N$, o\`u $\N$ est le radical
unipotent de $\P$.
On note $\P^-$ le sous-groupe parabolique oppos\'e \`a $\P$,
et $\N^-$ son radical unipotent.

\begin{defi}[\cite{BK1}, 6.1]
Soit $\K$ un sous-groupe de $\mult\A$.
\begin{enumerate}
\item[(i)] 
On dit que $\K$ est {\it d\'ecompos\'e}, 
ou qu'il admet une {\it d\'e\-com\-po\-si\-tion d'Iwahori},
rela\-tivement \`a $(\M,\P)$, si~:
\begin{equation*}
\K=(\K\cap\N^-)\cdot(\K\cap\M)\cdot(\K\cap\N).
\end{equation*}
\item[(ii)] 
Soit $\tau$ une repr\'esentation de $\K$.
La paire $(\K,\tau)$ est dite {\it d\'ecompos\'ee} 
rela\-tivement \`a $(\M,\P)$ si $\K$ est d\'ecompos\'e 
rela\-tivement \`a $(\M,\P)$ et si $\K\cap\N$, 
$\K\cap\N^-$ sont dans le noyau de $\tau$.
\end{enumerate}
\end{defi}

On utilisera le lemme suivant pour prouver que certains
sous-groupes de $\mult\A$ admettent une d\'e\-com\-po\-sition 
d'Iwahori.

\begin{lemm}[\cite{BH1}, 10.4]
\label{ProtoIwahoriX}
Soit $\R$ un $\o_\F$-r\'eseau dans $\A$ tel que 
$1+\R$ soit un sous-groupe de $\mult\A$.
On suppose que $\e^i\R\e ^j\subset\R$ pour chaque paire 
$(i,j)$.
Alors $1+\R$ admet une d\'e\-com\-po\-sition d'Iwahori 
rela\-tivement \`a $(\M,\P)$, pour tout sous-groupe 
para\-bo\-lique $\P$ de $\mult\A$ de facteur de Levi $\M$.
\end{lemm}

Soit $\L$ une $\Oo_\D$-suite de r\'eseaux de $\V$.
Pour chaque $i$, on pose $\L^i=\L\cap\V^i$.

\begin{defi}
\label{defERG}
On dit que (\ref{PipDickens}) est {\it conforme \`a $\L$}, ou encore
que $\L$ est {\it d\'ecompos\'ee par} (\ref{PipDickens}), si~:
\begin{equation}
\label{PipDickensLambda}
\L=\L^1\oplus\ldots\oplus\L^l.
\end{equation}
\end{defi}

Une d\'ecomposition (\ref{PipDickens}) est conforme \`a $\L$ si, 
et seulement si $\e^i\in\AA(\L)$ pour chaque $i$. 
Dans ce cas, on a $\L^i=\e^i\L$ et $\PP(\L)\cap\A^i=\PP(\L^i)$.

\begin{exem}
On suppose que (\ref{PipDickens}) est conforme \`a $\L$.
Pour $k\>1$, le groupe $\U_k(\L)$ admet une d\'e\-com\-po\-sition
d'Iwahori rela\-tivement \`a $(\M,\P)$, pour tout sous-groupe 
para\-bo\-lique $\P$ de $\mult\A$ de facteur de Levi $\M$.
\end{exem}

\subsection{}
\label{DefFondaScindee}

Soit $\A$ une $\F$-alg\`ebre centrale simple, soit $\V$ un $\A$-module
\`a gauche simple et soit $\D$ l'alg\`ebre oppos\'ee \`a $\End_{\A}(\V)$.
On reprend les notations des paragraphes pr\'ec\'edents.

\begin{defi}
Une \emph{strate} de $\A$ est un quadruplet $[\L,n,r,\b]$ constitu\'e
d'une $\o_\D$-suite $\L$ de $\V$, de deux entiers $r,n$ v\'erifiant
$0\<r\<n-1$ et d'un \'el\'ement $\b\in\PP_{-n}(\L)$.
Deux strates $[\L,n,r,\b_i]$, $i\in\{1,2\}$, sont dites
\emph{\'equivalentes} si $\b_2-\b_1\in\PP_{-r}(\L)$.  
\end{defi}

\begin{exem}
\label{LesLapinsNAimentPasLesCarottes}
Soit un caract\`ere additif $\psi_\F:\F\f\mult\CC$ trivial sur
$\p_\F$ mais pas sur $\o_\F$.
Soit $[\L,n,r,\b]$ une strate de $\A$ telle que 
$\lfloor n/2\rfloor\<r$.
Le caract\`ere $\psi_\b$ de $\U_{r+1}(\L)$ d\'efini par~:
\begin{equation*}
\psi_\b:x\mapsto\psi_\F\circ\tr_{\A/\F}(\b(x-1)),
\end{equation*}
ne d\'epend que de la classe d'\'equivalence de $[\L,n,r,\b]$.
\end{exem}

\'Etant donn\'ee une strate $[\L,n,r,\b]$ de $\A$, 
on note $\E$ la $\F$-alg\`ebre engendr\'ee par $\b$. 

\begin{defi}
La strate $[\L,n,r,\b]$ est dite {\it pure} si $\E$ est un corps, 
si la suite $\L$ est normalis\'ee par $\mult\E$ et si
$\v_{\L}(\b)=-n$.
\end{defi}

Soit $[\L,n,r,\b]$ une strate pure.
On note $\B$ le commutant de $\E$ dans $\A$.
Pour tout entier $k\in\ZZ$, on pose :
\begin{equation*}
\mathfrak{n}_k(\b,\L)=\{x\in\PP_0(\L)\ |\ \b x-x\b\in\PP_k(\L)\}.
\end{equation*}
Le plus petit entier $k\>\v_{\L}(\b)$
pour lequel le r\'eseau $\mathfrak{n}_{k+1}(\b,\L)$ est inclus dans 
$\PP_0(\L)\cap\B+\PP_1(\L)$ est not\'e $k_0(\b,\L)$ et porte le nom
d'\emph{exposant critique} de la strate $[\L,n,r,\b]$.  
(Il s'agit de la convention adopt\'ee dans \cite{St4} : dans le cas
o\`u $\E=\F$, on a $k_0(\b,\L)=-n$, tandis qu'avec la convention
d'usage dans \cite{BK,VS1}, on aurait $k_0(\b,\L)=-\infty$.)

\begin{defi}
\label{stratepure}
La strate $[\L,n,r,\b]$ est dite \emph{simple} si elle est pure et si
$r\<-k_0(\b,\L)-1$.  
\end{defi}

\subsection{}

Soit $[\L,n,r,\b]$ une strate simple de $\A$.
On pose $q=-k_0(\b,\L)$. 

\begin{defi}
\label{WeltAnschauung}
\begin{enumerate}
\item[(i)]
Si $q=n$, on dit que $\b$ est {\it minimal} sur $\F$.
\item[(ii)]
Si $q\<n-1$, une {\it approximation de $\b$ relativement \`a $\L$} 
est un \'el\'ement $\g\in\A$ tel que $[\L,n,q,\g]$ soit une strate
simple \'equivalente \`a $[\L,n,q,\b]$.
\end{enumerate}
\end{defi}

Selon \cite[Th\'eor\`eme 2.2]{VS3}, un \'el\'ement $\b$ qui n'est pas
minimal sur $\F$ admet des approximations relativement \`a
n'importe quelle suite de r\'eseaux $\E$-pure, et on peut m\^e\-me
choisir $\g$ de fa\c con que la sous-extension non ramifi\'ee maximale
de $\F(\g)$ soit incluse dans $\E$.
(Par exemple, si $\E/\F$ est non ramifi\'ee, on peut choisir
$\g\in\E$.) 
Ceci permet d'avoir acc\`es \`a la machinerie, d\'evelopp\'ee dans
\cite{BK}, des constructions par r\'ecurrence sur l'exposant critique.

\begin{prop}
\label{KingLear}
Soit $\V=\V^1\oplus\ldots\oplus\V^l$ une d\'ecomposition de $\V$
en sous-$\E\otimes_{\F}\D$-modules, qui soit conforme \`a $\L$.
\begin{enumerate}
\item[(i)]
Pour tout $1\<i\<l$, la strate $[\L^i,n,r,\e^i\b]$
est une strate simple de $\A^i$.
\item[(ii)]
Si $\b$ n'est pas minimal sur $\F$, il existe une approximation de
$\b$ relativement \`a $\L$ commutant aux $\e^i$.
\end{enumerate}
\end{prop}

\begin{proof}
Pour le (i), voir \cite[Proposition 2.28]{VS1}.
Pour le (ii), c'est ce que dit la preuve de 
\cite[Th\'eor\`eme 2.2]{VS3} (voir aussi {\it ibid.}, \S1.3).
\end{proof}

\subsection{}
\label{PairesSimplesEtRealisations}

Soit $(k,\b)$ une paire simple sur $\F$ au sens de \cite{BH1},
et soit $\E=\F(\b)$.
On appelle {\it donn\'ee admissible} pour la paire $(k,\b)$
un quintuplet $(\A,\ii,\V,\L,m)$ constitu\'e d'une $\F$-alg\`ebre 
centrale simple $\A$, d'un plongement de $\F$-alg\`ebres $\ii:\E\f\A$,
d'un $\A$-module \`a gauche simple $\V$, d'une $\o_{\D}$-suite 
$\E$-pure $\L$ de r\'eseaux de $\V$ et d'un entier $m$ v\'erifiant~:
\begin{equation*}
\left\lfloor\frac{m}{e(\L|\o_{\E})}\right\rfloor=k.
\end{equation*}
On pose $n=-v_{\L}(\ii\b)$.
D'apr\`es \cite[\S2.3.3]{VS1}, la strate $[\L,n,m,\ii\b]$ est une
strate simple de $\A$. 

\begin{defi}
La strate simple $[\L,n,m,\ii\b]$ est appel\'e une {\it r\'ealisation} 
de la paire simple $(k,\b)$. 
\end{defi}


\section{Caract\`eres simples} 
\label{CSsuites}

Soit $\A$ une $\F$-alg\`ebre centrale simple. 
Dans cette section, on associe \`a toute strate simple $[\L,n,m,\b]$
de $\A$ un ensemble $\Cc(\L,m,\b)$ de ca\-ract\`eres simples,
jouissant de propri\'et\'es remarquables de transfert et
d'entrelacement. 
Cette cons\-truc\-tion g\'en\'eralise les cons\-truc\-tions
effectu\'ees dans \cite{BK2,VS1,St4} lorsque $\A$ est
d\'eploy\'ee (\ie lorsque $\D=\F$) et celles effectu\'ees dans 
\cite{Grabitz,VS1} lors\-que  $\L$ est une $\Oo_\D$-cha\^\i ne
de r\'eseaux.

L'id\'ee suivie est la m\^eme que dans \cite{BK2} : on rajoute \`a
$\L$ une $\Oo_\D$-cha\^\i ne $\L^{\pt}$ d'un $\D$-espace vectoriel
ad\'equat $\V^{\pt}$, et on d\'efinit $\Cc(\L,m,\b)$ comme la projection 
de $\Cc(\L\oplus\L^{\pt},m,\b)$ sur le facteur $\Aut_\D(\V)$.

On fixe une fois pour toutes un caract\`ere additif
$\psi_\F:\F\f\mult\CC$ trivial sur $\p_\F$ mais pas sur $\o_\F$.
Toutes les constructions de caract\`eres simples d\'ependent
du choix de $\psi_\F$.

\subsection{} 
\label{TempsAnciens}

Soit $[\L,n,0,\b]$ une strate simple de $\A$ et soit $q=-k_0(\b,\L)$.
Dans ce paragraphe et le suivant, on suppose que $\L$ est stricte. 
Dans \cite{VS1}, on associe \`a une telle strate les objets 
suivants~:
\begin{itemize}
\item[(i)] 
Deux sous-$\o_\F$-ordres $\JJ(\b,\L)$ et $\HH(\b,\L)$ de $\AA(\L)$.
Ils ne d\'ependent que de la classe d'\'equivalence de $[\L,n,0,\b]$.
Chacun d'eux est filtr\'e par une suite d\'ecroissante d'id\'eaux
bilat\`eres~:
\begin{eqnarray*}
\JJ^k(\b,\L)&=&\JJ(\b,\L)\cap\PP_k(\L),\\
\HH^k(\b,\L)&=&\HH(\b,\L)\cap\PP_k(\L), \quad k\>0.
\end{eqnarray*}
Ces id\'eaux bilat\`eres sont en particulier
(avec les notations du \S\ref{PureteDeMoeurs})
des sous-$\BB(\L)$-bimodules de $\A$. 
On note $\J(\b,\L)$ (resp. $\H(\b,\L)$) le groupe multiplicatif de 
$\JJ(\b,\L)$ (resp. de $\HH(\b,\L)$).
De fa\c con similaire, chacun d'eux est filtr\'e par une suite
d\'ecroissante de sous-groupes ouverts compacts~:
\begin{eqnarray*}
\J^k(\b,\L)&=&\J(\b,\L)\cap\U_k(\L),\\
\H^k(\b,\L)&=&\H(\b,\L)\cap\U_k(\L), \quad k\>0.
\end{eqnarray*}
\item[(ii)] 
Pour tout entier $0\<m\<q-1$, un ensemble fini $\Cc(\L,m,\b)$ de
carac\-t\`e\-res de $\H^{m+1}(\b,\L)$ appel\'es 
{\it caract\`eres simples de niveau $m$}. 
Ces ca\-rac\-t\`eres v\'erifient une propri\'et\'e de fonctorialit\'e 
appel\'ee propri\'et\'e de {\it transfert}, dont l'\'etablissement
\'etait l'un des principaux objectifs de \cite{VS1}.
On reviendra largement dessus plus bas.
\end{itemize}

\subsection{}
\label{DelireMystique}

Dans ce paragraphe, on d\'ecrit, pour $\L$ stricte, le comportement de
$\Cc(\L,m,\b)$ par passage \`a une strate simple \'equivalente, puis
par augmentation du niveau. 
On pose $r=\lfloor q/2\rfloor+1$.

\begin{prop}
\label{GrosBoxon1}
Soit un entier $0\<l\<q$, et soit $[\L,n,l,\g]$ une strate simple
\'equivalente \`a $[\L,n,l,\b]$.
L'application $\t\mapsto\t\psi_{\g-\b}$ induit des bijections~: 
\begin{equation*}
\Cc(\L,m,\b)\f\Cc(\L,m,\g), \quad \lfloor l/2\rfloor\<m\<q-1.
\end{equation*}
\end{prop}

\begin{proof}
La preuve est analogue \`a celle de \cite{BK}, corollaire $3.3.18$ 
si $l=q$, et corollaire $3.3.20$(ii) si $l\neq q$.
Dans les deux cas, il suffit de rem\-pla\-cer \cite[Proposition 3.3.9]{BK}
par la proposition \cite[3.30]{VS1} et \cite[Proposition 2.4.11]{BK}
par \cite[Lemma 1.9]{Grabitz}.
\end{proof}

\begin{prop}
\label{GrosBoxon2}
Soit $[\L,n,m,\b]$ une strate simple de $\A$.
\begin{enumerate}
\item[(i)]
On suppose que $m\<\lfloor q/2\rfloor$.
Si deux caract\`eres de $\Cc(\L,m,\b)$ co\"\i nci\-dent sur
$\H^{\rr}(\b,\L)$, ils sont tordus l'un l'autre par un caract\`ere de 
$\U_{m+1}(\L)\cap\mult\B$ trivial sur $\U_{\rr}(\L)\cap\mult\B$ et se
factorisant par ${\rm N}_{\B/\E}$.
\item[(ii)]
On suppose que $m\<q-2$.
La restriction de $\Cc(\L,m,\b)$ \`a $\Cc(\L,m+1,\b)$
est surjective.
\end{enumerate}
\end{prop}

\begin{proof}
Le (i) est une cons\'equence de \cite[Lemme 3.24]{VS1}.
Pour le (ii), on proc\`ede par r\'ecurrence comme pour 
\cite[Corollary 3.3.21]{BK}.
Si $q=n$, on s\'epare deux cas.
Pour $\lfloor n/2\rfloor\<m$, on a~:
\begin{equation}
\label{LesHerissonsNonPlus}
\Cc(\L,m,\b)=\{\psi_{\b}\}
\end{equation}
({\it cf.} \cite[Lemme 3.23]{VS1}).
Pour $m\<\lfloor n/2\rfloor$, on utilise (i).
Ensuite, si $q\<n-1$, on choisit une strate simple $[\L,n,q,\g]$
\'equivalente \`a $[\L,n,q,\b]$ et on applique la proposition
\ref{GrosBoxon1}. 
\end{proof}

\subsection{}
\label{KarlAdam}

Soit $[\L,n,0,\b]$ une strate simple de $\A$.
Maintenant, $\L$ est une $\o_\D$-suite quelconque.
On note $e$ sa p\'eriode sur $\o_\D$ et on pose $q=-k_0(\b,\L)$.

Soit $(\V^{\pt},\L^{\pt})$ un couple constitu\'e d'un
$\E\otimes_\F\D$-module de type fini $\V^{\pt}$ et d'une $\o_\D$-suite
{\it stricte} $\E$-pure $\L^{\pt}$ de r\'eseaux de $\V^{\pt}$, de
p\'eriode $e$. 
On pose :
\begin{equation}
\label{decL}
\bar{\V}=\V\oplus\V^{\pt},
\quad
\bar{\L}=\L\oplus\L^{\pt}.
\end{equation}
La $\F$-alg\`ebre $\E$ se plonge naturellement dans 
$\bar{\A}=\End_{\D}(\bar{\V})$ et $\bar\L$ est une 
$\o_\D$-suite stricte $\E$-pure de r\'eseaux de $\bar\V$. 
La strate $[\bar\L,n,0,\b]$ est une stra\-te simple de $\bar\A$. 
(On renvoie \`a \cite[\S5]{BK2} pour ce proc\'ed\'e.)

On note $\bar\B$ le commutant de $\E$ dans $\bar\A$.
Par cons\-truc\-tion, la d\'ecomposition (\ref{decL}) est une 
d\'ecom\-po\-sition de $\bar\V$ en $\E\otimes\D$-modules, qui 
est conforme \`a $\L$, \ie que les pro\-jec\-teurs $\e:\bar\V\f\V$ et
$\e^{\pt}:\bar\V\f\V^{\pt}$ appartiennent tous deux \`a 
$\BB(\bar\L)=\AA(\bar\L)\cap\bar\B$.

On pose :
\begin{equation*}
\M=\Aut_{\D}(\V)\times\Aut_{\D}(\V^{\pt}).
\end{equation*}
Pour $k\>1$, les $\o_\F$-r\'eseaux $\JJ^{k}(\b,\bar\L)$ et
$\HH^{k}(\b,\bar\L)$ sont des sous-$\BB(\bar\L)$-bimodules 
de $\bar\A$.
D'apr\`es le lemme \ref{ProtoIwahoriX}, les groupes 
$\J^{k}(\b,\bar\L)$, $\H^{k}(\b,\bar\L)$ 
admettent cha\-cun une d\'e\-com\-po\-sition d'Iwahori 
relativement \`a $(\M,\P)$, pour tout sous-groupe parabolique 
$\P$ de $\Aut_{\D}(\bar\V)$ de facteur de Levi $\M$. 

On pose enfin $\G=\mult\A$ et $\bar\G=\bar\A^{\times}$.

\begin{theo}
\label{CoherenceDefJHC}
On suppose que $\L$ est stricte.
\begin{enumerate}
\item[(i)]
Pour $k\>0$, on a $\J^k(\b,\bar\L)\cap\G=\J^k(\b,\L)$ et 
$\H^k(\b,\bar\L)\cap\G=\H^k(\b,\L)$.
\item [(ii)]
Pour $0\<m\<q-1$ et pour $\t\in\Cc(\bar\L,m,\b)$, la paire
$(\H^{m+1}(\b,\bar\L),\t)$ est d\'ecompos\'ee par $(\M,\P)$ 
pour tout sous-groupe parabolique $\P$ de $\bar\G$ de facteur 
de Levi $\M$, et la restriction de $\t$ \`a $\H^{m+1}(\b,\L)$
est \'egale au transfert de $\t$ \`a $\Cc(\L,m,\b)$.
\end{enumerate}
\end{theo}

\begin{proof}
La preuve est analogue \`a celles de \cite[Proposition 7.1.12]{BK} et 
\cite[Proposition 7.1.19]{BK}.
Dans le cas o\`u $q\<n-1$, il suffit de choisir une approximation 
de $\b$ commutant \`a $\e$ et $\e^{\pt}$ 
({\it cf.} Proposition \ref{KingLear}(ii)) et,  
pour le (ii), on peut rem\-placer \cite[Corollary 3.3.21]{BK} par la 
proposition \ref{GrosBoxon2}. 
\end{proof}

\begin{rema}
\label{petitnom}
On a en fait un r\'esultat un peu plus g\'en\'eral que 
\ref{CoherenceDefJHC}(i). 
L'analogue de \cite[Proposition 7.1.12]{BK} est~:
\begin{eqnarray*}
\JJ^k(\b,\bar\L)\cap\A&=&\JJ^k(\b,\L),\\
\HH^k(\b,\bar\L)\cap\A&=&\HH^k(\b,\L), \quad k\>0.
\end{eqnarray*}
On v\'erifie \'egalement que
$\mathfrak{n}_k(\b,\bar\L)\cap\A=\mathfrak{n}_k(\b,\L)$,
pour $k\in\ZZ$.
\end{rema}

\subsection{}
\label{Wittg}

Soit $[\L,n,0,\b]$ une strate simple de $\A$, et soit $\bar\L$
d\'efini comme dans la section pr\'ec\'edente.
On lui associe deux $\o_\F$-ordres~: 
\begin{eqnarray*}
\label{ProtoTristesTropiques}
\JJ(\b,\L)&=&\JJ(\b,\bar\L)\cap\A,\\
\HH(\b,\L)&=&\HH(\b,\bar\L)\cap\A.
\end{eqnarray*}
De fa\c con similaire au \S\ref{TempsAnciens}, chacun d'eux est
filtr\'e par une suite d\'ecroissante d'id\'eaux bilat\`eres
$\JJ^k(\b,\L)$ et $\HH^k(\b,\L)$.
Ces id\'eaux bilat\`eres sont en particulier
des sous-$\BB(\L)$-bimodules de $\A$. 

On note $\J(\b,\L)$ (resp. $\H(\b,\L)$) le groupe multiplicatif de 
$\JJ(\b,\L)$ (resp. de $\HH(\b,\L)$), et chacun d'eux est filtr\'e par
une suite d\'ecroissante de sous-groupes ouverts compacts
$\J^k(\b,\L)$, $\H^k(\b,\L)$.
Si aucune confusion n'en r\'esulte, on notera $\J^k$ plut\^ot que
$\J^k(\b,\L)$. 
La m\^eme remarque vaut pour $\JJ^k(\b,\L)$, $\HH^k(\b,\L)$ 
et $\H^k(\b,\L)$.

\begin{prop}
\label{ventredelabaleine}
Pour $k\>0$, les groupes $\J^{k}$ et $\H^{k}$ sont 
normalis\'es par $\J$ et $\KK(\L)\cap\mult\B$.
\end{prop}

\begin{proof}
C'est une cons\'equence de \cite[Pro\-po\-sition 3.43]{VS1}.
\end{proof}

\begin{defi}
\label{defCSsuites}
Pour un entier $0\<m\<q-1$, on appelle 
{\it caract\`ere simple de niveau $m$}
attach\'e \`a la strate simple $[\L,n,0,\b]$ la restriction \`a
$\H^{m+1}(\b,\L)$ d'un caract\`ere simple de $\Cc(\bar\L,m,\b)$.
Ces caract\`eres forment un ensemble not\'e $\Cc(\L,m,\b)$.
\end{defi}

D'apr\`es le th\'eor\`eme \ref{CoherenceDefJHC},
(voir aussi la remarque \ref{petitnom}), 
cette d\'efinition co\"\i ncide avec celle de \cite[\S3.3]{VS1}
lorsque $\L$ est stricte.
Dans le cas contraire, elle d\'epend {\it a priori} de
$(\V^{\pt},\L^{\pt})$.
On va voir que ce n'est pas le cas. 

\begin{prop}
\label{Cenest}
Les groupes $\H^k(\b,\L)$ et $\J^k(\b,\L)$, ainsi que l'ensemble 
$\Cc(\L,m,\b)$, sont in\-d\'e\-pen\-dants du choix de
$(\V^{\pt},\L^{\pt})$. 
\end{prop}

\begin{proof}
Soit $(\V^{\bullet},\L^{\bullet})$ un autre couple comme au
\S\ref{KarlAdam}. 
On pose $\L^{\dag}=\L\oplus\L^{\pt}\oplus\L^{\bullet}$.
On note ${\rm R}^{\pt}$ la restriction de
$\H^{m+1}(\b,\L\oplus\L^{\pt})$ \`a $\H^{m+1}(\b,\L)$ et 
${\rm S}^{\pt}$ la restriction de $\H^{m+1}(\b,\L^{\dag})$ 
\`a $\H^{m+1}(\b,\L\oplus\L^{\pt})$.
On d\'efinit ${\rm R}^{\bullet}$ et ${\rm S}^{\bullet}$ de fa\c con
analogue en substituant $\L^{\bullet}$ \`a $\L^{\pt}$.
D'apr\`es le th\'eor\`eme \ref{CoherenceDefJHC}(ii), et puisque les
applications de transfert sont surjectives, on peut \'ecrire
$\Cc(\L,m,\b)$ comme l'image de $\Cc(\L^{\dag},m,\b)$ par 
${\rm R}^{\pt}\circ{\rm S}^{\pt}=
{\rm R}^{\bullet}\circ{\rm S}^{\bullet}$. 
Ceci prouve que $\Cc(\L,m,\b)$, et {\it a fortiori} $\H^{m+1}(\b,\L)$, 
sont ind\'ependants du choix de $(\V^{\pt},\L^{\pt})$ pour un entier 
$0\<m\<q-1$.

Pour $\H^k(\b,\L)$ et $\J^k(\b,\L)$, avec $k\>1$, on raisonne de fa\c
con analogue \`a partir du point (i) du th\'eo\-r\`eme
\ref{CoherenceDefJHC}. 
\end{proof}

\subsection{}
\label{ChangementBaseNonRamifie}

On rappelle bri\`evement le proc\'ed\'e de changement de base non
ramifi\'e d\'ecrit en d\'etail dans \cite{VS1}. 
Soit $\LL/\F$ une extension non ramifi\'ee maximale dans $\D$.
La $\LL$-alg\`ebre $\A\otimes_\F\LL$ est centrale simple d\'eploy\'ee,
et s'identifie cano\-niquement \`a $\End_{\LL}(\V)$.
On note $\L_{\LL}$ la suite $\L$ vue comme une $\Oo_\LL$-suite de
r\'eseaux de $\V$. 
On a une d\'ecomposition~:
\begin{equation*}
\E\otimes_\F\LL=\E^{1}\oplus\ldots\oplus\E^{l}
\end{equation*}
de la $\LL$-alg\`ebre $\E\otimes_\F\LL$ en une somme finie de $l\>1$
extensions de $\LL$, o\`u $l$ est le pgcd de $[\LL:\F]$ et du degr\'e
r\'esiduel de $\E/\F$.
On note $\e^i$ l'idem\-po\-tent de $\E\otimes_\F\LL$ correspondant \`a la
projection sur $\E^i$, ce qui d\'efinit une d\'e\-com\-po\-si\-tion de
$\V$ en une somme de $\E\otimes_\F\LL$-modules $\V^i=\e^i\V$, qui est
con\-for\-me \`a $\L_{\LL}$. 
On pose $\L^i_{\LL}=\e^i\L_{\LL}$ et $\b^i=\e^i\b$.
Le r\'esultat suivant vient de \cite{VS1} ({\it cf.} \S2.3.4)~:

\begin{theo}
\label{UnterSchugung}
$[\L^i_{\LL},n,m,\b^i]$ 
est une strate simple de $\End_{\LL}(\V^i)$.
\end{theo}

\begin{rema}
Signalons une erreur dans la preuve qui en est donn\'ee dans
\cite{VS1} ({\it cf.} Th\'eor\`eme 2.30). 
Avec les notations de {\it loc. cit.}, on n'a pas en g\'en\'eral 
\'egalit\'e entre $k_0({\rm e}^i\b,{\bf\L}^i)$ et
$k_0(\e^i\b,\bar\L^i)$ 
(le corollaire $2.21$ ne s'applique pas \`a $\e^i\b$). 
Par contre, puisque le r\'eseau 
$\mathfrak{n}_k(\e^i\b,\bar\L^i)\cap{\bf A}^i$ 
est \'egal \`a $\mathfrak{n}_k(\e^i\b,\e^i{\bf\L})$ 
({\it cf.} Proposition 2.20), 
on a $k_0({\rm e}^i\b,{\bf\L}^i)\<k_0(\e^i\b,\bar\L^i)$,
ce qui suffit pour conclure.
Voici un cas o\`u l'\'egalit\'e n'a pas lieu.
Choisissons par exemple $\E/\F$ non ramifi\'ee, $\b$ non minimal sur
$\F$ ({\it i.e.} $q\neq n$) et $\LL$ contenant $\E$.
Alors chaque $\E^i$ est de degr\'e $1$ sur $\LL$, \ie que
$\b^i$ est scalaire, donc minimal sur $\LL$.
\end{rema}

\subsection{}
\label{TransDEP}

Dans ce paragraphe, on \'etend le transfert aux
caract\`eres simples attach\'es \`a une suite de r\'eseaux 
({\it cf.} \S\ref{Wittg}). 

\begin{prop}
\label{transfertCSsuites}
La restriction de $\H^{m+1}(\b,\bar\L)$ \`a $\H^{m+1}(\b,\L)$ induit 
une bijection de $\Cc(\bar\L,m,\b)$ sur $\Cc(\L,m,\b)$.
\end{prop}

\begin{proof}
La preuve est analogue \`a celle de \cite[Th\'eor\`eme 3.12]{VS1}.
Il suffit de remplacer \cite[Proposition 3.2.4]{BK} par 
(\ref{LesHerissonsNonPlus}) et \cite[Proposition 3.2.5]{BK} par la
proposition \ref{GrosBoxon2}(i) et, dans le cas o\`u $q\<n-1$, 
de choisir une approximation de $\b$ commutant \`a $\e$ et $\e^{\pt}$ 
({\it cf.} Proposition \ref{KingLear}(ii)). 
\end{proof}

Soit $(k,\b)$ une paire simple sur $\F$ et soit $\E=\F(\b)$ 
({\it cf.} \S\ref{PairesSimplesEtRealisations}). 
On note $[k,\b]$ la strate simple de $\End_\F(\E)$ correspondant \`a
$(k,\b)$, \ie correspondant \`a l'unique ordre h\'er\'editaire de
$\End_\F(\E)$ normalis\'e par $\mult\E$ 
(voir par exemple \cite[\S2.3.3]{VS1}), et on note $\Cc(k,\b)$
l'ensemble des caract\`eres simples (de niveau $k$) attach\'e \`a la
strate simple $[k,\b]$.

Soit $[\L,n,m,\ii\b]$ une r\'ealisation de $(k,\b)$ dans $\A$.
Lorsque la suite $\L$ est stricte, on a une bijection canonique de
$\Cc(k,\b)$ dans $\Cc(\L,m,\ii\b)$ appel\'ee application de transfert 
({\it cf.} \cite[\S3.3.3]{VS1}) et not\'ee $\tau_{\L,m,\ii\b}$.
Lorsque $\L$ est quelconque, on pose la d\'efinition suivante.

\begin{defi}
\label{datc}
L'{\it application de transfert} de $\Cc(k,\b)$ \`a $\Cc(\L,m,\ii\b)$ 
est l'application bijective compos\'ee de $\tau_{\bar\L,m,\ii\b}$ avec 
la restriction de $\H^{m+1}(\ii\b,\bar\L)$ \`a $\H^{m+1}(\ii\b,\L)$.
On la note $\tau_{\L,m,\ii\b}$.
\end{defi}

D'apr\`es le th\'eor\`eme \ref{CoherenceDefJHC}, cette d\'efinition 
du transfert co\"\i ncide avec celle de \cite[\S3.3]{VS1} lorsque $\L$
est stricte.
Dans le cas contraire, elle d\'epend {\it a priori} de 
$(\V^{\pt},\L^{\pt})$.
Un raisonnement analogue \`a celui de la preuve de la proposition 
\ref{Cenest} montre que ce n'est pas le cas.

Si $[\L,n,m,\ii\b]$ et $[\L',n',m',\ii'\b]$ sont deux r\'ealisations 
d'une m\^eme paire, on d\'efinit une application de transfert entre 
$\Cc(\L,m,\ii\b)$ et $\Cc(\L',m',\ii'\b)$ en composant 
$\tau_{\L,m,\ii\b}^{-1}$ avec $\tau_{\L',m',\ii'\b}$.
On a une propri\'et\'e de transitivit\'e \'evidente.

\begin{exem}
Dans la situation du \S\ref{KarlAdam}, l'application de transfert de 
$\Cc(\bar\L,m,\b)$ \`a $\Cc(\L,m,\b)$ est la restriction de 
$\H^{m+1}(\b,\bar\L)$ \`a $\H^{m+1}(\b,\L)$.
\end{exem}

Voici une premi\`ere propri\'et\'e du transfert (qu'on peut qualifier 
de transfert {\it interne}, dans la mesure o\`u on change la suite
$\L$, mais pas le groupe $\G$).

\begin{theo}
\label{BOUlouLOU}
Soient $[\L,n,m,\b]$ et $[\L',n',m',\b]$ deux strates simples
de $\A$.
On suppose que~:
\begin{equation*}
\Big\lfloor\frac{m}{e(\L|\o_\E)}\Big\rfloor=
\Big\lfloor\frac{m'}{e(\L'|\o_\E)}\Big\rfloor.
\end{equation*}
Alors pour tout $\t\in\Cc(\L,m,\b)$, le transfert de $\t$ \`a
$\Cc(\L',m',\b)$ co\"\i ncide avec $\t$ sur 
$\H^{m+1}(\b,\L)\cap\H^{m'+1}(\b,\L')$.
\end{theo}

\begin{rema}
Dans le cas o\`u $\L$ et $\L'$ sont strictes, ce r\'esultat, quoique
annonc\'e dans \cite{VS1} ({\it cf.} Exemple 3.54) et d\'ej\`a
utilis\'e dans \cite{VS2}, n'y est pas d\'emontr\'e.
On donne ici une d\'emonstration dans le cas g\'en\'eral,
lorsque $\L$ et $\L'$ sont quelconques.
\end{rema}

\begin{proof}
La preuve s'effectue en quatre \'etapes.
\begin{enumerate}
\item[(i)]
On choisit un caract\`ere simple $\t\in\Cc(\L,m,\b)$ et on note
$\t'$ son transfert \`a $\Cc(\L',m',\b)$.
On fixe un couple $(\V^{\pt},\L^{\pt})$ pour $\L$ et un couple
$(\V'^{\pt},\L'^{\pt})$ pour $\L'$ ({\it cf.} \S\ref{KarlAdam}),
de telle sorte que $\V'^{\pt}=\V^{\pt}$.
On pose $\bar\L=\L\oplus\L^{\pt}$ et $\bar\L'=\L'\oplus\L'^{\pt}$.
On note $\bar\t$ (resp. $\bar\t'$) le caract\`ere simple de
$\Cc(\bar\L,m,\b)$ (resp. de $\Cc(\bar\L',m,\b)$) qui 
prolonge $\t$ (resp. $\t'$).
Ainsi $\bar\t'$ est le transfert de $\bar\t$.
\item[(ii)]
On suppose momentan\'ement que $\A$ est d\'eploy\'ee. 
Alors $\bar\t$ et $\bar\t'$ co\"\i ncident sur l'inter\-section
$\H^{m+1}(\b,\bar\L)\cap\H^{m'+1}(\b,\bar\L')$
d'apr\`es \cite[Theorem 3.6.1]{BK}.
Par restriction \`a $\G$, les caract\`eres $\t$ 
et $\t'$ co\"\i ncident sur 
$\H^{m+1}(\b,\L)\cap\H^{m'+1}(\b,\L')$.
\item[(iii)]
La $\F$-alg\`ebre $\A$ est \`a nouveau quelconque, mais on suppose
maintenant que $\L$ et $\L'$ sont strictes.
On choisit une extension non ramifi\'ee $\LL/\F$ maximale dans $\D$
({\it cf.} \S\ref{ChangementBaseNonRamifie}).
Pour chaque $1\<i\<l$, on a une strate simple
$[\L^{i}_\LL,n,m,\b^{i}]$ de $\End_{\LL}(\V^i)$.
On fixe un caract\`ere additif $\psi_\LL:\LL\f\mult\CC$ trivial sur
$\p_\LL$ mais pas sur $\o_\LL$, dont la restriction \`a $\F$ est
$\psi_\F$. 
On choisit un caract\`ere simple $\t^i\in\Cc(\L^{i}_\LL,m,\b^{i})$
de telle sorte que la famille $\{\t^i\}$ d\'efinisse un caract\`ere
quasi-simple dont la restriction \`a $\H^{m+1}(\b,\L)$ soit $\t$
({\it cf.} \cite{VS1}, \S3.2.4).
On note $\t'^{i}$ le transfert de $\t^i$ \`a $\Cc(\L'^{i},m',\b^i)$.
D'apr\`es \cite[Th\'eor\`eme 3.53]{VS1}, la famille $\{\t'^{i}\}$
d\'efinit un caract\`ere quasi-simple dont la restriction \`a
$\H^{m+1}(\b,\L')$ est $\t'$.
D'apr\`es (ii), les caract\`eres $\t^{i}$ et $\t'^{i}$ co\"\i ncident
sur $\H^{m+1}(\b^{i},\L^{i})\cap\H^{m'+1}(\b^{i},\L'^{i})$.
Ceci prouve le th\'eor\`eme dans le cas o\`u $\L$ et $\L'$ sont
strictes. 
\item[(iv)]
On revient maintenant au cas g\'en\'eral.
Les caract\`eres $\bar\t$ et $\bar\t'$ co\"\i ncident sur 
l'inter\-section $\H^{m+1}(\b,\bar\L)\cap\H^{m'+1}(\b,\bar\L')$
d'apr\`es (iii). 
Par restriction \`a $\G$, les caract\`eres $\t$ 
et $\t'$ co\"\i ncident sur 
$\H^{m+1}(\b,\L)\cap\H^{m'+1}(\b,\L')$.
\end{enumerate}
Ceci termine la preuve du th\'eor\`eme \ref{BOUlouLOU}.
\end{proof}

\subsection{}
\label{TransDEPpri}

Dans ce paragraphe, on d\'ecrit le comportement de $\Cc(\L,m,\b)$ par
passage \`a une strate simple \'equivalente, puis par augmentation du
niveau.
On pose $\rr=\lfloor q/2\rfloor+1$.

\begin{prop}
\label{GrosBoxon1L}
Soit un entier $0\<l\<q$, et soit $[\L,n,l,\g]$ une strate simple
\'equivalente \`a $[\L,n,l,\b]$.
L'application $\t\mapsto\t\psi_{\g-\b}$ induit des bijections~: 
\begin{equation*}
\Cc(\L,m,\b)\f\Cc(\L,m,\g), \quad \lfloor l/2\rfloor\<m\<q-1.
\end{equation*}
\end{prop}

\begin{proof}
Par transfert de $\L$ \`a $\bar\L$, on se ram\`ene au cas o\`u 
$\L$ est stricte, puis on applique la proposition \ref{GrosBoxon1}.
\end{proof}

\begin{prop}
\label{GrosBoxon2L}
Soit $[\L,n,m,\b]$ une strate simple de $\A$.
\begin{enumerate}
\item[(i)]
On suppose que $m\<\lfloor q/2\rfloor$.
Si deux caract\`eres de $\Cc(\L,m,\b)$ co\"\i nci\-dent sur
$\H^{\rr}(\b,\L)$, ils sont tordus l'un l'autre par un caract\`ere de 
$\U_{m+1}(\L)\cap\mult\B$ trivial sur $\U_{\rr}(\L)\cap\mult\B$ et se
factorisant par ${\rm N}_{\B/\E}$.
\item[(ii)]
On suppose que $m\<q-2$.
La restriction de $\Cc(\L,m,\b)$ \`a $\Cc(\L,m+1,\b)$
est surjective.
\end{enumerate}
\end{prop}

\begin{proof}
Par transfert de $\L$ \`a $\bar\L$, on se ram\`ene au cas o\`u 
$\L$ est stricte, puis on applique la proposition \ref{GrosBoxon2}.
\end{proof}

On termine par le r\'esultat suivant.

\begin{theo}
\label{PaireDecomposeeTheta}
Soit 
$\V=\V^1\oplus\ldots\oplus\V^l$ une d\'ecomposition de $\V$
en sous-$\E\otimes_{\F}\D$-modules, qui soit conforme \`a $\L$.
Soit $\M$ le sous-groupe de Levi de $\G$ correspondant.
Soit un entier $1\<i\<l$.
\begin{enumerate}
\item[(i)]
Pour $k\>0$, on a~:
\begin{equation*}
\J^k(\b,\L)\cap\Aut_{\D}(\V^i)=\J^k(\b,\L^i),
\quad
\H^k(\b,\L)\cap\Aut_{\D}(\V^i)=\H^k(\b,\L^i).
\end{equation*}
\item[(ii)]
Pour $0\<m\<q-1$ et pour $\t\in\Cc(\L,m,\b)$, la paire
$(\H^{m+1}(\b,\L),\t)$ est d\'ecompos\'ee par $(\M,\P)$ 
pour tout sous-groupe parabolique $\P$ de $\G$ de facteur 
de Levi $\M$, et la restriction de $\t$ \`a $\H^{m+1}(\b,\L^i)$
est \'egale au transfert de $\t$ \`a $\Cc(\L^{i},m,\b)$.
\end{enumerate}
\end{theo}

\begin{proof}
On fixe un couple $(\V^{\pt},\L^{\pt})$ comme au \S\ref{KarlAdam}, 
et on applique le th\'eor\`eme \ref{CoherenceDefJHC} \`a 
la d\'ecomposition 
$\L\oplus\L^{\pt}\oplus\L^{\pt}=(\L^i\oplus\L^{\pt})\oplus
((\bigoplus_{j\neq i}\L^j)\oplus\L^{\pt})$.
D'abord, on a~:
\begin{equation*}
\J^k(\b,\L\oplus\L^{\pt}\oplus\L^{\pt})\cap\Aut_{\D}(\V^i\oplus\V^{\pt})
=\J^k(\b,\L^i\oplus\L^{\pt}).
\end{equation*}
Si on projette sur $\Aut_{\D}(\V^i)$, on obtient l'\'egalit\'e
voulue pour $\J^k$.
Avec un rai\-son\-nement analogue, on obtient l'\'egalit\'e pour $\H^k$.
Ensuite, soit $\t$ un caract\`ere simple de $\Cc(\L,m,\b)$, et soit 
$\tilde\t\in\Cc(\L\oplus\L^{\pt}\oplus\L^{\pt},m,\b)$ prolongeant $\t$.
La res\-triction de $\tilde\t$ \`a $\H^{m+1}(\b,\L^i\oplus\L^{\pt})$
est \'egale au transfert de $\t$ \`a $\Cc(\L^i\oplus\L^{\pt},m,\b)$.
Si on restreint \`a $\H^{m+1}(\b,\L^i)$, on obtient l'\'egalit\'e
voulue entre restriction de $\t$ \`a $\H^{m+1}(\b,\L^i)$ et 
transfert de $\t$ \`a $\Cc(\L^i,m,\b)$.
Enfin, il reste \`a prouver que $\t$ est trivial sur chaque 
sous-groupe de la forme $1+\Hom(\V^i,\V^j)$, avec $i\neq j$.
C'est vrai pour $\tilde\t$ sur 
$1+\Hom(\V^i\oplus\V^{\pt},(\bigoplus_{j\neq i}\V^j)\oplus\V^{\pt})$,
et on obtient le r\'esultat par restriction \`a 
$1+\Hom(\V^i,\V^j)$.
\end{proof}

\subsection{}
\label{metasuivantpre}

Soit $[\L,n,0,\b]$ une strate simple de $\A$ et soit 
$q=-k_0(\b,\L)$.
On pose $\ss=\lceil q/2\rceil$.
Pour $k\>1$, on pose :
\begin{equation*}
\label{defMM}
\MM_k(\b,\L)=
\PP_{k}(\L)\cap\mathfrak{n}_{-q+k}(\b,\L)+
\JJ^{\ss}(\b,\L) 
\end{equation*}
et on pose $\Omega_k(\b,\L)=1+\MM_k(\b,\L)$, qui est un sous-groupe 
ouvert compact de $\U_{k_{0}}(\L)$, avec 
$k_{0}=\min\{k,\ss\}$. 
Souvent, on notera simplement $\Om_k(\L)$, ou m\^eme $\Om_k$.
Remarquer que la convention de notation est diff\'erente de celle
utilis\'ee dans \cite{VS1,St4} 
(dont le $\MM_{k}$ correspond \`a notre $\MM_{q-k}$).

Soit $(\V^{\pt},\L^{\pt})$ un couple comme au \S\ref{KarlAdam}, 
dont on reprend les notations. 
L'objectif de ce paragraphe et du suivant est de prouver le r\'esultat
suivant. 

\begin{prop}
\label{WiganPier}
Soit $k\>1$.
On a :
\begin{equation*}
\Om_k(\b,\bar\L)\bar\B^{\times}\Om_k(\b,\bar\L)\cap\G
=\Om_k(\b,\L)\B^{\times}\Om_k(\b,\L).
\end{equation*}
\end{prop}

La d\'emonstration se fait, comme dans \cite[\S3.1]{VS1} 
({\it cf.} \cite[Lemme 3.7]{VS1}), par changement de base 
non ramifi\'e.

\begin{lemm}
\label{C3}
Soit $k\>1$.
Pour tout $b\in\B^{\times}$, on a~:
\begin{equation*}
\Om_k(\L)b\Om_k(\L)\cap\B^{\times}
=\left(\Om_k(\L)\cap\B^{\times}\right)b
\left(\Om_k(\L)\cap\B^{\times}\right).
\end{equation*}
\end{lemm}

\begin{proof}
Si l'on pose $k_0=\min\{k,\ss\}$, c'est une cons\'equence de~: 
\begin{equation*}
\U_{k_0}(\L)\cap\mult\B\subset\Om_k(\L)\subset\U_{k_0}(\L)
\end{equation*}
et de la propri\'et\'e d'intersection simple 
\cite[Corollaire 3.3]{VS1}. 
La majoration est imm\'ediate.
Pour la minoration, il faut remarquer d'une part que 
$\BB(\L)$ est inclus dans $\JJ(\b,\L)$, ce qui implique 
$\PP_{k_0}(\L)\cap\B\subset\JJ^{\ss}(\b,\L)$ dans le cas
o\`u $k\>\ss$, et d'autre part que $\PP_{k_0}(\L)\cap\B$
est inclus dans $\PP_{k}(\L)\cap\mathfrak{n}_{-q+k}(\b,\L)$
dans le cas o\`u $k\<\ss$.
\end{proof}

\begin{lemm}
\label{C4}
On a $\Om_k(\b,\bar\L)\cap\G=\Om_k(\b,\L)$ pour $k\>1$.
\end{lemm}

\begin{proof}
On v\'erifie que~:
\begin{equation*}
\PP_{k}(\bar\L)\cap\mathfrak{n}_{-q+k}(\b,\bar\L)\cap\A=
\PP_{k}(\L)\cap\mathfrak{n}_{-q+k}(\b,\L).
\end{equation*}
Le r\'esultat est une cons\'equence de la d\'efinition 
de $\JJ^{k}(\b,\L)$ et du fait que, si $\R,\SS$ sont des 
sous-$\BB(\bar\L)$-bimodules de $\bar\A$, on a 
$(\R+\SS)\cap\A=\R\cap\A+\SS\cap\A$.
\end{proof}

\subsection{}
\label{metasuivant}

Soit $\F^{\sharp}$ une extension finie non ramifi\'ee de $\F$, 
de degr\'e premier au degr\'e r\'esiduel de $\E/\F$ et \`a la
dimension de $\D$ sur $\F$, et de groupe de Galois not\'e $\Gg$.
On pose~:
\begin{equation*}
\A^{\sharp}=\A\otimes_\F\F^{\sharp},
\quad
\V^{\sharp}=\V\otimes_{\F}\F^{\sharp},
\quad
\L^{\sharp}=\L\otimes_{\o_\F}\o_{\F^{\sharp}}
\quad
\D^{\sharp}=\D\otimes_{\F}\F^{\sharp}.
\end{equation*}
Ainsi $\F^{\sharp}[\b]$ est un corps, $\D^{\sharp}$ est une
$\F^{\sharp}$-alg\`ebre \`a division, $\V^{\sharp}$ est un  
$\A^{\sharp}$-module simple et $\A^{\sharp}$ s'identifie 
naturellement \`a $\End_{\D^{\sharp}}(\V^{\sharp})$.

\begin{prop}
\label{LaPossibiliteDUneIle}
La strate $[\L^{\sharp},n,r,\b\otimes1]$ de $\A^{\sharp}$ est simple.
\end{prop}

\begin{proof}
La d\'emonstration suit formellement \cite[\S2]{VS1} (voir notamment
la proposition 2.9 et les corollaires 2.10 et 2.11), le fait que
$\F^{\sharp}/\F$ d\'eploie $\A$ n'y jouant aucun r\^ole.
\end{proof}

On suppose en outre que l'extension $\F^{\sharp}/\F$ est non triviale,
et on choisit dans $\F^{\sharp}$ une racine de l'unit\'e $\xi$, 
non triviale et d'ordre premier \`a la carac\-t\'e\-ris\-ti\-que 
r\'esiduelle $p$ de $\F$.
On note $\Delta$ le groupe cyclique engendr\'e par
$\xi\cdot\e+\e^{\pt}$ et $\M^{\sharp}$ son centralisateur dans 
$\Aut_{\D^{\sharp}}(\bar\V^{\sharp})$, \ie le groupe des points 
fixes de $\Aut_{\D^{\sharp}}(\bar\V^{\sharp})$ par $\Delta$.
Si on identifie $\A$ \`a la $\F$-alg\`ebre $\A^{\sharp\Gg}$ des 
$\Gg$-invariants de $\A^{\sharp}$, on a $\M^{\sharp\Gg}=\M$.

\begin{rema}
Ce proc\'ed\'e permet de calculer l'intersection de certaines 
parties de $\bar\G$ avec $\M$ par des m\'ethodes de descente
comme en \cite[\S2.4]{VS1}.
Le changement de base est n\'ecessaire dans le cas o\`u le 
corps r\'esiduel de $\F$ n'a que deux \'el\'ements.
\end{rema}

\begin{proof}[Preuve de la proposition \ref{WiganPier}]
On note $\bar\B^{\sharp}$ le commutant de $\E$ dans
$\bar\A^{\sharp}$.
Si on applique le lemme \ref{C3} \`a la strate simple
$[\bar\L^{\sharp},n,0,\b]$, alors, compte tenu de 
\cite[Proposition 2.36]{VS1} et de \cite[Lemme 2.35]{VS1}, 
on obtient~:
\begin{equation*}
\Om_k(\bar\L^{\sharp})\bar\B^{\sharp\times}
\Om_k(\bar\L^{\sharp})\cap\M^{\sharp}=
\left(\Om_k(\bar\L^{\sharp})\cap\M^{\sharp}\right)
\left(\bar\B^{\sharp\times}\cap\M^{\sharp}\right)
\left(\Om_k(\bar\L^{\sharp})\cap\M^{\sharp}\right)
\end{equation*}
puis, en projetant sur $\A^{\sharp\times}$~:
\begin{equation}
\label{SaloperieDeFormule}
\Om_k(\bar\L^{\sharp})\bar\B^{\sharp\times}
\Om_k(\bar\L^{\sharp})\cap\A^{\sharp\times}=
\Om_k(\L^{\sharp})\B^{\sharp\times}\Om_k(\L^{\sharp}).
\end{equation}
Il reste \`a calculer les points fixes de 
(\ref{SaloperieDeFormule}) par $\Gg$.
Pour le membre de gauche, on applique \cite[Lemme 2.35]{VS1} 
en tenant compte de \cite[Proposition 2.41]{VS1} et du lemme
\ref{C3} appliqu\'e \`a la strate simple $[\bar\L^{\sharp},n,0,\b]$.
On obtient :
\begin{equation*}
\Om_k(\bar\L^{\sharp})\bar\B^{\sharp\times}
\Om_k(\bar\L^{\sharp})\cap\G=
\Om_k(\bar\L)\bar\B^{\times}\Om_k(\bar\L)\cap\G.
\end{equation*}
Pour le membre de droite, on applique \cite[Lemme 2.35]{VS1} 
en tenant compte de \cite[Proposition 2.41]{VS1} et du lemme
\ref{C3} appliqu\'e \`a la strate simple $[\L^{\sharp},n,0,\b]$.
On obtient :
\begin{equation*}
\Om_k(\L^{\sharp})\B^{\sharp\times}\Om_k(\L^{\sharp})\cap\G=
\Om_k(\L)\B^{\times}\Om_k(\L),
\end{equation*}
ce qui termine la d\'emonstration.
\end{proof}

\subsection{}

On fixe un entier $0\<m\<q-1$ et un caract\`ere simple
$\t\in\Cc(\L,m,\b)$.
On rappelle que, si $\K$ est un sous-groupe de $\G$, et si $\chi$ est
un caract\`ere de $\K$, l'{\it entrelacement} de $\chi$ dans $\G$,
not\'e $\I_{\G}(\chi)$, est l'ensemble des \'el\'ements $g\in\G$ pour 
les\-quels $\chi$ et son caract\`ere conjugu\'e $\chi^g$ 
co\"\i ncident sur $\K\cap g^{-1}\K g$.

\begin{theo}
\label{entrelacementCSsuites}
On a $\I_{\G}(\t)=\Om_{q-m}(\b,\L)\mult\B\Om_{q-m}(\b,\L)$.
\end{theo}

\begin{proof}
La preuve est analogue \`a celle de \cite[Proposition 2.5]{St4}.
Soit $\bar\t\in\Cc(\bar\L,m,\b)$ le caract\`ere prolongeant $\t$.
D'apr\`es \cite[Th\'eor\`eme 3.50]{VS1}, \cite[Lemma 2.1]{St4}
et la proposition \ref{PaireDecomposeeTheta} appliqu\'ee \`a
(\ref{decL}), on a~: 
\begin{equation*}
\begin{split}
\I_{\bar\G}(\bar\t_{|\H^{m+1}(\b,\bar\L)\cap\M})\cap\M&=
\I_{\bar\G}(\bar\t)\cap\M\\
&=\Om_{q-m}(\bar\L)\bar\B^{\times}\Om_{q-m}(\bar\L)\cap\M.
\end{split}
\end{equation*}
On applique la proposition \ref{WiganPier}, puis on projette 
sur $\G$, ce qui donne~:
\begin{equation*}
\I_{\bar\G}(\bar\t_{|\H^{m+1}(\b,\bar\L)\cap\M})\cap\G=
\Om_{q-m}(\b,\L)\mult\B\Om_{q-m}(\b,\L).
\end{equation*}
Le membre de gauche vaut $\I_{\G}(\t)$, ce qui termine la 
d\'emonstration.
\end{proof}

\begin{prop}
\label{normBKSSTcs}
Tout ca\-rac\-t\`e\-re simple de $\Cc(\L,m,\b)$
est normalis\'e par $(\KK(\L)\cap\B^{\times})\Om_{q-m}(\b,\L)$.
\end{prop}

\begin{proof}
Soit $\bar\t\in\Cc(\bar\L,m,\b)$ le caract\`ere prolongeant $\t$.
D'apr\`es le th\'eo\-r\`e\-me \ref{entrelacementCSsuites}
et \cite[Th\'eor\`eme 3.50]{VS1}, le caract\`ere $\bar\t$
est normalis\'e par le grou\-pe 
$(\KK(\bar\L)\cap\bar\B^{\times})\Om_{q-m}(\bar\L)$
et $\t$ est entrelac\'e par $(\KK(\L)\cap\B^{\times})\Om_{q-m}(\L)$.
\end{proof}

\subsection{}
\label{CorMod}
\def\fr#1{\smash{\mathop{\f}\limits^{#1}}}

Soit $[\L,n,0,\b]$ une strate simple de $\A$.
On pose $\tilde\A=\End_{\F}(\V)$, qu'on identifie \`a 
$\A\otimes_{\F}\End_{\A}(\V)$.

\begin{defi}
Une {\it corestriction mod\'er\'ee sur $\A$ relative \`a $\E/\F$} 
est un homomorphisme de $\B$-bimodules $s:\A\f\B$ tel que 
$\tilde{s}=s\otimes{\rm id}_{\End_{\A}(\V)}$
soit une corestriction mod\'er\'ee sur $\tilde\A$ relative \`a 
$\E/\F$ au sens de \cite[Definition 1.3.3]{BK}.
\end{defi}

Bien entendu, lorsque $\A$ est d\'eploy\'ee sur $\F$, 
cette d\'efinition co\"\i ncide avec celle de \cite{BK},
puisque, dans ce cas, on a $\End_{\A}(\V)=\F$.

D'apr\`es \cite[Lemma 4.2.1]{Br4}, on a un moyen de construire 
des corestrictions mo\-d\'e\-r\'ees sur $\A$ relatives \`a $\E/\F$, 
ce qui prouve qu'il en existe.
On a les propri\'et\'es suivantes. 

\begin{prop}
\label{Piniwi}
Soit $s$ une corestriction mod\'er\'ee sur $\A$ relative 
\`a $\E/\F$.
\begin{enumerate}
\item[(i)]
Si $s'$ est une corestriction mod\'er\'ee sur $\A$ relative
\`a $\E/\F$, il existe $u\in\mult\o_\E$ tel que $s'=us$.
\item[(ii)]
Soit $\V=\V^1\oplus\V^2$ une d\'ecomposition de $\V$ en 
sous-$\E\otimes_\F\D$-modules, qui soit conforme \`a $\L$.
Pour $i\in\{1,2\}$, la restriction $s_i$ de $s$ \`a $\A^i$ est une
co\-res\-tric\-tion mo\-d\'e\-r\'ee sur $\A^i$ relative \`a $\E/\F$.
\end{enumerate}
\end{prop}

\begin{proof}
Les deux sont vrais lorsque la $\F$-alg\`ebre $\A$ est d\'eploy\'ee. 
Pour (i), il existe donc $u\in\mult\o_\E$ tel que
$\tilde{s}'=u\tilde{s}$, ce qui implique $s'=us$.
Pour (ii), la restriction $\tilde{s}_i$ de $\tilde{s}$ \`a
$\tilde{\A}^i$ est une corestriction mo\-d\'e\-r\'ee sur
$\tilde{\A}^i$ relative \`a $\E/\F$ \'egale \`a 
$s_i\otimes{\rm id}_{\End_{\A^i}(\V^i)}$.
\end{proof}

Soit $\V=\V^1\oplus\V^2$ une d\'ecomposition de $\V$ en 
sous-$\E\otimes_\F\D$-modules, qui soit conforme \`a $\L$.
On note $\M$ le sous-groupe de Levi correspondant.
Si $\mathfrak{l}$ est un sous-$\o_\F$-r\'e\-seau de $\A$, on pose
$\mathfrak{l}^{ij}=\mathfrak{l}\cap\A^{ij}$ pour $i,j\in\{1,2\}$.

Soit $s$ une corestriction mod\'er\'ee sur $\A$ relative 
\`a $\E/\F$.
Pour $x\in\A$, on pose $a_\b(x)=\b x-x\b$.
On pose $\psi_\A=\psi_\F\circ\tr_{\A/\F}$.
Pour toute partie $\R$ de $\A$, on note~:
\begin{equation*}
\R^*=\{a\in\A\ |\ \psi_\A(ax)=1,\ x\in\R\}
\end{equation*}
le dual de $\R$ relativement \`a $\psi_\A$.

\begin{prop}
\label{ExactSequences}
Pour $0\<m\<q-1$, la suite~:
\begin{equation}
\label{setilE}
0\f\PP_{q-m}(\L)\cap\B\f\MM_{q-m}(\b,\L)
\fr{a_\b}(\HH^{m+1}(\b,\L))^*\fr{s}\PP_{-m}(\L)\cap\B\f0
\end{equation}
est exacte.
Si on d\'esigne cette suite par 
$0\f\mathfrak{l}_{1}\f\mathfrak{l}_{2}
\f\mathfrak{l}_{3}\f\mathfrak{l}_{4}\f0$, 
alors la suite~:
\begin{equation}
\label{setilEij}
0
\f h^{-1}\mathfrak{l}^{ij}_{1}h+\mathfrak{l}^{ij}_{1}
\f h^{-1}\mathfrak{l}^{ij}_{2}h+\mathfrak{l}^{ij}_{2}
\f h^{-1}\mathfrak{l}^{ij}_{3}h+\mathfrak{l}^{ij}_{3}
\f h^{-1}\mathfrak{l}^{ij}_{4}h+\mathfrak{l}^{ij}_{4}
\f 0
\end{equation}
est exacte pour tout $h\in\mult\B\cap\M$ et tous $i,j\in\{1,2\}$.
\end{prop}

\begin{proof}
On note $\tilde\L$ la $\o_\F$-suite de $\V$ sous-ja\-cente \`a $\L$.
On note $\tilde\B$ le commutant de $\E$ dans $\tilde\A$,
on note $*$ la dualit\'e relativement \`a 
$\psi_{\F}\circ\tr_{\tilde{\A}/\F}$,
on note $\tilde{a}_{\b}$
l'application $x\mapsto \b x-x\b$ d\'efinie sur $\tilde\A$ 
et on pose $\tilde{s}=s\otimes{\rm id}_{\End_{\A}(\V)}$. 
La strate $[\tilde\L,n,m,\b]$ est une strate simple 
de $\tilde\A$ de m\^eme exposant critique que $[\L,n,m,\b]$
({\it cf.} \cite[2.23]{VS1}).
On lui applique \cite[Lemma 6.3]{BK2}.
La suite~:
\begin{equation}
\label{setil}
0
\f\PP_{q-m}(\tilde\L)\cap\tilde\B\f\MM_{q-m}(\b,\tilde\L)
\fr{\tilde{a}_\b}(\HH^{m+1}(\b,\tilde\L))^*
\fr{\tilde{s}}\PP_{-m}(\tilde\L)\cap\tilde\B
\f0
\end{equation}
est exacte.
Si on d\'esigne cette suite par 
$0\f\tilde{\mathfrak{l}}_{1}\f\tilde{\mathfrak{l}}_{2}
\f\tilde{\mathfrak{l}}_{3}\f\tilde{\mathfrak{l}}_{4}\f0$, 
alors la suite~:
\begin{equation}
\label{setilij}
0
\f h^{-1}\tilde{\mathfrak{l}}^{ij}_{1}h+\tilde{\mathfrak{l}}^{ij}_{1}
\f h^{-1}\tilde{\mathfrak{l}}^{ij}_{2}h+\tilde{\mathfrak{l}}^{ij}_{2}
\f h^{-1}\tilde{\mathfrak{l}}^{ij}_{3}h+\tilde{\mathfrak{l}}^{ij}_{3}
\f h^{-1}\tilde{\mathfrak{l}}^{ij}_{4}h+\tilde{\mathfrak{l}}^{ij}_{4}
\f 0
\end{equation}
est exacte pour tout $h\in\B^{\times}\cap\M$ et tous $i,j\in\{1,2\}$.

\begin{lemm}
\label{lk}
On a $\tilde{\mathfrak{l}}_{k}\cap\A=\mathfrak{l}_{k}$
pour $1\<k\<4$.
\end{lemm}

\begin{proof}
Pour $k\in\{1,4\}$, c'est imm\'ediat.
Pour les autres cas, on choisit un couple $(\V^{\pt},\L^{\pt})$ comme
au \S\ref{KarlAdam}. 
Par d\'efinition ({\it cf.} \cite[(54) et (64)]{VS1} et 
\S2.3), on a~:
\begin{equation*}
\HH^{m+1}(\b,\L)=\HH^{m+1}(\b,\tilde\L\oplus\tilde\L^{\pt})\cap\A
\end{equation*}
et le membre de droite est \'egal \`a $\HH^{m+1}(\b,\tilde\L)\cap\A$
d'apr\`es (\ref{petitnom}).
On en d\'eduit que $\tilde{\mathfrak{l}}_{3}\cap\A=\mathfrak{l}_{3}$
\`a l'aide de \cite[Lemme 2.45]{VS1}.
Un raisonnement analogue permet d'obtenir
$\tilde{\mathfrak{l}}_{2}\cap\A=\mathfrak{l}_{2}$.
\end{proof}

Pour terminer la d\'emonstration de la proposition \ref{ExactSequences}, 
il reste \`a v\'erifier que 
l'exac\-titude de (\ref{setilE}) et (\ref{setilEij}) est 
con\-ser\-v\'ee par restriction \`a $\A$.
Ceci se fait, comme dans \cite[\S4]{Br4}, en choisissant une 
extension non ramifi\'ee $\LL/\F$ maximale dans $\D$, et en 
appliquant successivement le foncteur des $\LL$-invariants puis 
le foncteur des $\Gal(\LL/\F)$-invariants
(voir aussi \cite[\S2.4]{VS1}).
\end{proof}

De fa\c con analogue, on d\'emontre \`a partir de 
\cite[Corollary 1.4.10]{BK}~: 

\begin{prop}
\label{Cor1410}
Pour $k\in\ZZ$, la suite~:
\begin{equation*}
\label{setilERR}
0\f\PP_{q+k}(\L)\cap\B\f\PP_{q+k}(\L)\cap\mathfrak{n}_k(\b,\L)
\fr{a_\b}\PP_{k}(\L)\fr{s}\PP_{k}(\L)\cap\B\f0
\end{equation*}
est exacte.
\end{prop}

\subsection{}

On \'etablit une propri\'et\'e de non-d\'eg\'en\'erescence des
caract\`eres simples, qui g\'en\'eralise \cite[Theorem 3.4.1]{BK}.
Pour $x,y\in\G$, on note $[x,y]$ le commutateur de $x$ et $y$.

\begin{lemm}
\label{QuelleHorreurInsoutenableALInfiniRepetee}
Soit $\t\in\Cc(\L,m,\b)$ avec $\lfloor q/2\rfloor\<m\<q-1$.
Soient deux entiers $k,l\>1$ tels que $k+l\>m+1$ et $k+2l\>q+1$.  
On suppose qu'on est dans l'une des situations suivantes :
\begin{enumerate}
\item
$x\in1+\PP_{k}(\L)\cap\mathfrak{n}_{k-q}(\b,\L)$ 
et $y\in1+\PP_{l}(\L)\cap\mathfrak{n}_{l-q}(\b,\L)$.
\item
$x\in1+\PP_{k}(\L)\cap\mathfrak{n}_{k-q}(\L)$ et $y\in\J^{l}(\b,\L)$.
\item
$x\in\J^{k}(\b,\L)$ et $y\in\J^{l}(\b,\L)$.
\end{enumerate}
Alors $[x,y]\in\H^{m+1}(\b,\L)$, et on a 
$\t([x,y])=\psi_{x^{-1}\b x-\b}(y)$.
\end{lemm}

\begin{proof}
Par transfert de $\L$ \`a $\bar\L$ ({\it cf.} \S\ref{KarlAdam}),
on se ram\`ene au cas o\`u $\L$ est stricte.
Puis on applique \cite[Lemmes 3.25--3.27]{VS1}.
\end{proof}

\begin{prop}
\label{bfsndcs}
Soit $\t\in\Cc(\L,m,\b)$ avec $m\<\lfloor q/2\rfloor$.
L'appli\-ca\-tion~:
\begin{equation*}
(x,y)\mapsto\t([x,y]), 
\quad x,y\in\J^{m+1},
\end{equation*}
induit une forme altern\'ee non d\'eg\'en\'er\'ee :
\begin{equation*}
\boldsymbol{k}_{\t}:\J^{m+1}/\H^{m+1}\times\J^{m+1}/\H^{m+1}\f\mult\CC. 
\end{equation*}
\end{prop}

\begin{proof}
La preuve est analogue \`a celle de \cite[Proposition 3.9]{VS1}.
Il suffit de remplacer \cite[Theorem 3.4.1]{BK} par 
\cite[Th\'eor\`eme 3.52]{VS1} et \cite[Proposition 3.2.12]{BK}
par le lemme \ref{QuelleHorreurInsoutenableALInfiniRepetee}.
\end{proof}


\section{Le processus de raffinement}
\label{praf}

Soit $\A$ une $\F$-alg\`ebre centrale simple, soit $\G$ son groupe
multiplicatif et soit $\pi$ une repr\'esentation irr\'eductible 
de niveau non nul de $\G$.
Dans cette section, on prouve que trois cas seulement peuvent se 
produire ({\it cf.} Th\'eor\`eme \ref{??})~:
{(a)} ou bien il existe une strate simple $[\L,n,0,\b]$ de $\A$,
avec $\L$ stricte, telle que $\pi$ contienne un caract\`ere 
simple $\t\in\Cc(\L,0,\b)$, 
{(b)} ou bien $\pi$ contient une strate scind\'ee 
({\it cf.} D\'efinition \ref{stratefonda}), 
{(c)} ou bien $\pi$ contient un caract\`ere scind\'e 
({\it cf.} D\'efinition \ref{CarScinde}).

\subsection{}
\label{Diotime}

On fixe un $\A$-module \`a gauche simple $\V$.
L'ensemble $\Seq(\V,\Oo_\D)$ des $\o_\D$-suites de r\'eseaux de $\V$ 
est muni d'une structure affine provenant de celle de l'immeuble de 
$\G$ sur $\F$, qu'on peut d\'ecrire de la fa\c con suivante 
(voir \cite{Br3}).
Dans ce paragraphe, on note $n$ la dimension de $\V$ sur $\D$.

\begin{rema}
\label{WARNING}
Sauf mention explicite du contraire, dans ce qui suit, 
les suites de r\'eseaux
sont des $\o_\D$-suites, et on note $e(\L)$ la p\'eriode sur 
$\o_\D$ d'une suite $\L$.
\end{rema}

\`A chaque base $b$ de $\V$ sur $\D$ correspond
d'une part un isomorphisme $\V\simeq\D^{n}$ de $\D$-espaces vectoriels
\`a droite, d'autre part l'ensemble $\Seq^{b}(\V,\Oo_\D)$ des $\o_\D$-suites 
de r\'eseaux d\'ecompos\'ees par $b$.
Pour chaque entier $i$, on a une fonction affine 
$a_i:\Seq^{b}(\V,\Oo_\D)\f\RR$ telle que, pour 
$\L\in\Seq^{b}(\V,\Oo_\D)$ et $k\in\RR$, on ait~:
\begin{equation}
\label{defLk}
\L_k=\bigoplus\limits_{i=1}^n\p_\D^{\lceil k/e(\L)-a_i(\L)\rceil}.
\end{equation}
Cette base permet d'identifier les $\F$-alg\`ebres $\A$ et $\Mat_{n}(\D)$, 
donc de faire de $\A$ un $\D$-espace vectoriel \`a droite par transport 
de structure.
Pour chaque couple d'entiers $(i,j)$, on pose $\a_{ij}=a_i-a_j$.
Pour $\L\in\Seq^{b}(\V,\Oo_\D)$, la $\o_\F$-suite $\PP(\L)$
est munie d'une structure de $\o_\D$-suite et, pour $k\in\RR$, on a~:
\begin{equation}
\label{decarf}
\PP_k(\L)=
\bigoplus\limits_{i,j=1}^n\p_D^{\lceil k/e(\L)-\a_{ij}(\L)\rceil}.
\end{equation}
En d'autres termes, un \'el\'ement $a\in\A$ appartient \`a $\PP_k(\L)$
si et seulement si, pour chaque $1\<i,j\<n$, la valuation normalis\'ee 
de $a_{ij}$ est sup\'erieure ou \'egale \`a $k/e(\L)-\a_{ij}(\L)$.


\subsection{}

On \'etablit une liste de lemmes techniques.
La remarque \ref{WARNING} vaut toujours.

\begin{lemm}
\label{UATruth}
Soit $\L\in\Seq(\V,\Oo_\D)$ et soit un entier $m\in\ZZ$. 
Il existe une suite $\L'\in\Seq(\V,\Oo_\D)$ stricte et 
un entier $m'\in\ZZ$ tels qu'on ait~:
\begin{equation}
\label{Rouze}
\PP_{-m}(\L)\subset\PP_{-m'}(\L'),
\quad 
\frac{m'}{e(\L')}\<\frac{m}{e(\L)}.
\end{equation}
\end{lemm}

\begin{proof}
La d\'emonstration  est tr\`es proche de celle de 
\cite[Proposition 2.3]{HM}.
On se contente d'en donner les grandes lignes, 
et on renvoie \`a {\it loc. cit.} pour les d\'etails.
D'abord, on remarque que, pour $i\in\ZZ$, l'application 
naturelle~: 
\begin{equation}
\label{quotientadd}
\PP_{-m}(\L)/\PP_{-m+1}(\L)\f\bigoplus\limits_{l=0}^{e(\L)-1}
\Hom_{k_\D}(\L_{i+l}/\L_{i+l+1},\L_{i+l-m}/\L_{i+l-m+1}) 
\end{equation}
est un isomorphisme de $k_\D$-espaces vectoriels.
Donc  si $\L_i\neq\L_{i+1}$, on a l'\'egalit\'e 
$\PP_{-m}(\L)\L_i=\L_{i-m}$.
En d'autres termes, la partie $\PP_{-m}(\L)$ est {\it taut} par 
rapport \`a $\L$ au sens de \cite{HM},  
\ie que que $\PP_{-m}(\L)$ op\`ere sur l'ensemble 
$\Ll=\{\L_i\ |\ i\in\ZZ\}$.
Comme $\L$ n'est pas n\'ecessairement stricte, 
$\PP_{-m}(\L)$ n'est pas n\'ecessairement 
{\it completely taut}.
On note $\Ll^{'}$ la plus grande partie de $\Ll$
sur laquelle l'action de $\PP_{-m}(\L)$ est bijective,
et on choisit une $\o_\D$-suite stricte $\L'$ de classe de
translation $\Ll^{'}$.
Soit $m'\in\ZZ$ l'entier d\'efini par~:
\begin{equation}
\label{preuvehk}
\PP_{-m}(\L)\L'_i=\L'_{i-m'},
\quad i\in\ZZ.
\end{equation}
En raisonnant comme dans la preuve de \cite[Proposition 2.3]{HM}, 
on obtient~:
\begin{equation*}
\frac{m'}{e(\L')}\<\frac{m}{e(\L)}.
\end{equation*}
D'apr\`es (\ref{preuvehk}), on a l'inclusion cherch\'ee.
\end{proof}

\begin{rema}
(\ref{Rouze}) reste vrai si l'on remplace $m'$ par 
$me(\L')/e(\L)\in\QQ$.
\end{rema}

\begin{lemm}
\label{cotcot}
Soient $\L,\L'\in\Seq(\V,\Oo_\D)$, et soient $m,m'\in\QQ$.
On pose~:
\begin{equation}
\label{NotSegment}
\L(t)=(1-t)\L+t\L',
\quad
e(t)=e(\L(t)),
\quad
\frac{m(t)}{e(t)}=(1-t)\frac{m}{e(\L)}+t\frac{m'}{e(\L')}
\end{equation}
pour $t\in[0,1]$ rationnel.
Alors~:
\begin{equation*}
\PP_{-m}(\L)\cap\PP_{-m'}(\L')\subset\PP_{-m(t)}(\L(t)).
\end{equation*}
\end{lemm}

\begin{proof}
On choisit une base de $\V$ sur $\D$ d\'ecomposant $\L$ et $\L'$.
Pour $1\<i,j\<n$, on pose~:
\begin{equation}
\label{Rhoiji}
\rho_{ij}(t)=-\frac{m(t)}{e(t)}-\a_{ij}(\L(t)),
\end{equation}
qui est une application affine.
Compte tenu de (\ref{decarf}), il suffit de prouver que~:
\begin{equation*}
\max(\lceil\rho_{ij}(0)\rceil,\lceil\rho_{ij}(1)\rceil)
\>
\lceil(1-t)\rho_{ij}(0)+t\rho_{ij}(1)\rceil,
\end{equation*}
ce qui est imm\'ediat.
\end{proof}

\begin{coro}
\label{cotcotcoot}
Sous les hypoth\`eses de \ref{cotcot}, on suppose que~:
\begin{equation*}
\PP_{-m}(\L)\subset\PP_{-m'}(\L').
\end{equation*}
Alors, pour tous $0\<s\<t\<1$ rationnels, on a~:
\begin{equation*}
\PP_{-m(s)}(\L(s))\subset\PP_{-m(t)}(\L(t)).
\end{equation*}
\end{coro}

\begin{proof}
L'hypoth\`ese implique que $\rho_{ij}(0)\>\rho_{ij}(1)$,
\ie que chacune des fonctions affines $\rho_{ij}$ est 
d\'ecroissante.
Le r\'esultat s'ensuit.
\end{proof}

\begin{coro}
\label{corcotcot}
Soient $\L_1,\ldots,\L_r\in\Seq(\V,\Oo_\D)$ et soit $m\in\QQ$.
On pose~:
\begin{equation*}
\L'=\frac{1}{r}\sum\limits_{i=1}^{r}\L_i,
\quad
\frac{m'}{e(\L')}=\frac{1}{r}\sum\limits_{i=1}^{r}\frac{m}{e(\L_i)}.
\end{equation*}
Alors~:
\begin{equation*}
\bigcap\limits_{i=1}^{r}\PP_{m}(\L_i)\subset\PP_{m'}(\L').
\end{equation*}
\end{coro}

\begin{proof}
Par r\'ecurrence sur $r$ \`a partir du lemme \ref{cotcot}.
\end{proof}

\begin{lemm}
\label{FinsFaciles}
Soient $\L,\L'\in\Seq(\V,\Oo_\D)$, et soient $m,m'\in\ZZ$.
Il existe des couples 
$(\L_{0},m_{0}),\ldots,(\L_{l},m_{l})\in\Seq(\V,\Oo_\D)\times\QQ$ 
tels que~:
\begin{equation*}
\label{inclu}
\PP_{\lfloor m_{{k+1}}\rfloor+1}(\L_{{k+1}})\subset\PP_{m_{k}}(\L_{k}),
\quad
0\<k\<l-1,
\end{equation*}
et tels que $(\L_0,m_0)=(\L,m)$ et $(\L_l,m_l)=(\L',m')$. 
\end{lemm}

\begin{proof}
On choisit une base de $\V$ sur $\D$ d\'ecomposant $\L$ et $\L'$, 
et on reprend les notations (\ref{NotSegment}) et (\ref{Rhoiji}).
Chaque $\r_{ij}$ est une application affine, de sorte que
$t\mapsto\lceil\r_{ij}(t)\rceil$ ne prend qu'un nombre fini de valeurs.
Donc, compte tenu de (\ref{decarf}), l'application~:
\begin{equation}
\label{appliARP}
t\mapsto\PP_{m(t)}(\L(t))
\end{equation}
ne prend qu'un nombre fini de valeurs sur l'intervalle 
rationnel $[0,1]\cap\QQ$.
On choisit une suite strictement croissante $t_0,\ldots,t_l$ de 
cet intervalle telle que les r\'eseaux $\PP_{m(t_k)}(\L(t_k))$
d\'ecrivent les valeurs successives prises par (\ref{appliARP}).
On pose $\L_k=\L(t_k)$ et $m_k=m(t_k)$.
Pour prouver l'inclusion, il suffit de montrer que~:
\begin{equation*}
\lceil\r_{ij}(t_k)\rceil\<\lfloor\r_{ij}(t_{k+1})\rfloor+1
\end{equation*}
pour $1\<i,j\<m$ et $0\<k\<l-1$, ce qui est imm\'ediat. 
\end{proof}

\begin{lemm}
\label{FinsFacilesEpsilon}
Soient $\L,\L'\in\Seq(\V,\Oo_\D)$, soient $m,m'\in\ZZ$
et soit $\varepsilon\in\RR_+^{\times}$.
Il existe des couples 
$(\L_{0},m_{0}),\ldots,(\L_{l},m_{l})\in\Seq(\V,\Oo_\D)\times\QQ$ 
tels que~:
\begin{equation*}
\PP_{m_{k+1}}(\L_{k+1})\subset\PP_{m_{k}-e(\L_k)\varepsilon}(\L_{k}),
\quad
0\<k\<l-1,
\end{equation*}
et tels que $(\L_0,m_0)=(\L,m)$ et $(\L_l,m_l)=(\L',m')$. 
\end{lemm}

\begin{proof}
On choisit une base de $\V$ sur $\D$ d\'ecomposant $\L$ et $\L'$, 
et on reprend les notations (\ref{NotSegment}) et (\ref{Rhoiji}).
Compte tenu de (\ref{decarf}), l'inclusion~:
\begin{equation*}
\PP_{m(t)}(\L(t))\subset\PP_{m(s)-e(s)\varepsilon}(\L(s))
\end{equation*}
a lieu, pour $s,t\in[0,1]$ rationnels, si et seulement si on a
$\lceil\r_{ij}(t)\rceil\>\lceil\r_{ij}(s)-\varepsilon\rceil$.
Puisque $\rho_{ij}$ est affine, il suffit de choisir un entier 
$l\>1$ suffisamment grand et de poser $t_k=k/l$, pour $0\<k\<l$,
puis $\L_k=\L(t_k)$ et $m_k=m(t_k)$.
\end{proof}

\subsection{}
\label{Parmenide}

Soit $[\L,m,m-1,b]$ une strate de $\A$.
On pose $e=e(\L|\o_\F)$ et on note
$g$ le plus grand diviseur commun \`a $e$ et $m$.
On choisit une uniformisante $\varpi_{\F}$ de $\F$, et on pose
$y_b=\varpi_{\F}^{m/g}b^{e/g}$, que l'on consid\`ere comme un
\'el\'ement de $\End_{\F}(\V)$. 
Son polyn\^ome caract\'eristique est \`a  coefficients dans $\o_\F$,
et la r\'eduction modulo $\p_\F$ de celui-ci est appel\'ee le 
{\it polyn\^ome caract\'eristique} de la strate.
On le note $\varphi_b$.
Il est \`a coefficients dans $k_\F$
({\it cf.} \cite[\S2.2]{St3}).

\begin{defi}
\label{stratefonda}
La strate $[\L,m,m-1,b]$ est dite \emph{fondamentale} 
(resp. \emph{scind\'ee}) si son polyn\^ome caract\'eristique 
$\varphi_b\in k_\F[\X]$
n'est pas une puis\-sance de $\X$ (resp. a au moins deux facteurs 
irr\'eductibles distincts).
\end{defi}

\begin{rema}
Le polyn\^ome $\varphi_b$ d\'epend de l'uniformisante choisie,
mais pas les notions de strate fondamentale et de strate scind\'ee. 
\end{rema}

Les r\'esultats suivants g\'en\'eralisent respectivement 
\cite[Proposition 1.2.2]{Br4} et \cite[Theorem 1.2.5]{Br4}.

\begin{prop}
\label{RafNonFonda}
Soit $[\L,m,m-1,b]$ une strate non fondamentale de $\A$.
Il existe un entier $m'\in\ZZ$ et une $\o_\D$-suite stricte 
$\L'$ tels que~:
\begin{equation*}
b+\PP_{1-m}(\L)\subset\PP_{-m'}(\L'),
\quad 
\frac{m'}{e(\L')}<\frac{m}{e(\L)}.
\end{equation*}
\end{prop}

\begin{proof}
On proc\`ede par changement de base non ramifi\'e.
On choisit un couple $(\V^{\pt},\L^{\pt})$ comme au \S\ref{KarlAdam} 
et une extension non ramifi\'ee $\F^{\sharp}$ de $\F$ comme 
au \S\ref{metasuivant}, dont on reprend les notations.
On pose $\bar{b}=b\oplus0$, de sor\-te que la strate 
$[\bar\L^{\sharp},m,m-1,\bar{b}]$ de $\bar\A^{\sharp}$ 
est non fondamentale.
On peut donc ap\-pli\-quer \cite[Proposition 1.2.2]{Br4}.
On en tire une $\o_{\D^{\sharp}}$-suite stricte $\Ll$ 
de $\bar\V^{\sharp}$ et un entier $k\in\ZZ$ tels que~:
\begin{equation}
\label{momentchocolat}
\bar{b}+\PP_{1-m}(\bar\L^{\sharp})\subset\PP_{-k}(\Ll),
\quad 
\frac{k}{e(\Ll|\o_{\D^{\sharp}})}<\frac{m}{e(\L)}.
\end{equation}
Puisque le membre de gauche est stable par le groupe fini 
$\Gg\ltimes\Delta$, la relation (\ref{momentchocolat}) est
toujours valable si l'on remplace $\Ll$ par un de ses
conjugu\'es par ce groupe.
Si on applique le corollaire \ref{corcotcot} \`a la famille
des conjugu\'es de $\Ll$, on obtient~:
\begin{equation}
\label{TempoK}
\bar{b}+\PP_{1-m}(\bar\L^{\sharp})\subset\PP_{-k'}(\Ll'),
\quad 
\frac{k'}{e(\Ll'|\o_{\D^{\sharp}})}
=\frac{k}{e(\Ll|\o_{\D^{\sharp}})}<\frac{m}{e(\L)},
\end{equation}
o\`u $\Ll'$ d\'esigne la moyenne de $\Ll$ relativement \`a 
$\Gg\ltimes\Delta$, \ie l'iso\-bary\-centre des conjugu\'es 
de $\Ll$.
La suite $\Ll'$ est \`a la fois invariante par $\Gg$ et 
invariante par $\Delta$, \ie que $\Delta$ est contenu dans 
$\U(\Ll')$, donc que $\Ll'$ est d\'e\-com\-po\-s\'ee par la 
d\'ecomposition $\bar\V^{\sharp}=\V^{\sharp}\oplus\V^{\pt\sharp}$.
En projetant (\ref{TempoK}) sur $\A$, on obtient~:
\begin{equation*}
{b}+\PP_{1-m}(\L)\subset\PP_{-k'}(\Ll'\cap\V),
\quad\frac{k'}{e(\Ll'\cap\V)}<\frac{m}{e(\L)}.
\end{equation*}
Enfin, on applique le lemme \ref{UATruth} \`a la suite $\Ll'\cap\V$ :
il existe une $\o_\D$-suite stricte $\L'$ de $\V$ et un entier
$m'\in\ZZ$ v\'e\-ri\-fiant l'inclusion et l'in\'egalit\'e voulues.
\end{proof}

\begin{prop}
\label{RafSimlple}
Soit $[\L,m,m-1,b]$ une strate fondamentale non scin\-d\'ee de $\A$.
Il existe une strate simple $[\L',m',m'-1,b']$ avec $\L'$ stricte, 
v\'erifiant~:
\begin{equation*}
b+\PP_{1-m}(\L)\subset b'+\PP_{1-m'}(\L'),
\quad 
\frac{m'}{e(\L')}=\frac{m}{e(\L)}.
\end{equation*}
\end{prop}

\begin{proof}
On note $\tilde\L$ la $\o_\F$-suite sous-jacente \`a $\L$.
L'\'el\'ement carac\-t\'eristique $y_b$ 
est dans $\AA(\tilde\L)$, et sa r\'eduction modulo 
$\mathfrak{P}(\tilde\L)$ est inver\-sible, puisque son polyn\^ome 
caract\'eristique $\varphi_b$ est une puissance d'un poly\-n\^o\-me 
irr\'eductible distinct de $\X$.
D'apr\`es (\ref{quotientadd}), on en d\'eduit que $y_b\in\U(\tilde\L)$, 
puis que $b\in\KK(\L)$.

Soit $\Ll$ une $\o_\D$-suite stricte de $\V$ telle que
$\PP_0(\Ll)=\PP_0(\L)$. 
En particulier, $b$ normalise $\Ll$.
D'apr\`es \cite[Lemma 2.1.9(i)]{Br4}, si on pose 
$k=-\v_{\Ll}(b)$, la strate $[\Ll,k,k-1,b]$ est 
fondamentale non scind\'ee, et on a :
\begin{equation}
\label{WhiteRabbit}
\PP_{-m}(\L)=\PP_{-k}(\Ll),
\quad
\PP_{1-m}(\L)=\PP_{1-k}(\Ll).
\end{equation}
D'apr\`es \cite[Theorem 1.2.5]{Br4}, il existe une strate 
simple $[\L',m',m'-1,b']$ avec $\L'$ stricte, telle que~:
\begin{equation*}
b+\PP_{1-k}(\Ll)\subset b'+\PP_{1-m'}(\L').
\end{equation*}
Compte tenu de (\ref{WhiteRabbit}), ceci met fin \`a la 
d\'emonstration.
\end{proof}

\subsection{}

Soit $[\L,n,0,\b]$ une strate simple de $\A$.
On pose $q=-k_0(\b,\L)$ et on fixe un entier $1\<m\<q-1$.
On pose $r=\lfloor q/2\rfloor+1$ et $s=\lceil q/2\rceil$.

\begin{lemm}
\label{IndTau}
Soit $m_0=\min\{m,s\}$ et soit $\vartheta$ un caract\`ere de $\H^{m_0}$
dont la restriction \`a $\H^{m+1}$ est dans $\Cc(\b,m,\L)$.
Il existe une unique repr\'esentation irr\'eductible $\tau$ de
$\U_{m_0}(\L)$ dont la restriction \`a $\H^{m_0}$ contient $\vartheta$.
\end{lemm}

\begin{proof}
La d\'emonstration est similaire \`a celle de \cite[Lemma 8.1.8]{BK}.
D'apr\`es la
proposition~\ref{QuelleHorreurInsoutenableALInfiniRepetee}, la forme
altern\'ee~:
\begin{equation*}
(x,y)\mapsto\vartheta([x,y]),
\quad x,y\in\J^{m_0},
\end{equation*}
ne d\'epend que de la restriction de $\vartheta$
\`a $\H^{m+1}$.
D'apr\`es la proposition \ref{bfsndcs}, c'est donc une forme
non d\'eg\'en\'er\'ee.
Il existe donc une unique repr\'esentation irr\'eductible $\mu$
de $\J^{m_0}$ dont la restriction \`a $\H^{m_0}$ contient $\vartheta$.
Il reste \`a prouver que l'entrelacement de $\mu$ dans $\U_{m_0}(\L)$ est
contenu dans $\J^{m_0}$, et l'induite~:
\begin{equation*}
\tau=\Ind_{\J^{m_0}}^{\U_{m_0}(\L)}(\mu)
\end{equation*}
sera irr\'eductible.
L'entrelacement de $\mu$ est contenu dans~:
\begin{equation*}
\I_{\G}(\t)=\Om_{q-m}(\b,\L)\mult\B\Om_{q-m}(\b,\L)
\end{equation*}
par le th\'eor\`eme \ref{entrelacementCSsuites}. 
Il faut donc montrer que~:
\begin{equation*}
\Om_{q-m}(\b,\L)\cap \U_{m_0}(\L) \subset \J^{m_0}(\b,\L),
\end{equation*}
ce qui revient \`a montrer~:
\begin{equation}\label{contwithm0}
\PP_{q-m}(\L)\cap\mathfrak{n}_{-m}(\b,\L)\cap\PP_{m_0}(\L)\subset
\JJ^{m_0}(\b,\L).
\end{equation}
On continue par r\'ecurrence sur $k_0(\b,\L)$. 
Si $m<\rr$, l'\'equation \eqref{contwithm0} est impliqu\'ee 
par~\cite[Proposition 3.1.10(i)]{BK} dans le cas d\'eploy\'e 
$\D=\F$, d'o\`u le cas g\'en\'eral en appliquant le lemme~\ref{lk}.

On suppose donc que $m\>\rr$ (et $m_0=\ss$). 
Si $\b$ est minimal sur $\F$, alors $\JJ^{m_0}(\b,\L)=\PP_{m_0}(\L)$
et il n'y a rien \`a d\'emontrer.
Sinon, soit $[\L,n,q,\g]$ une strate simple \'equivalente \`a
$[\L,n,q,\b]$, soit $\B_1$ le commutant de $\g$ dans $\A$ et soit
$q_1=-k_0(\g,\L)$.
Alors on a~:
\begin{eqnarray*}
\PP_{q-m}(\L)\cap\mathfrak{n}_{-m}(\b,\L)&=&
\PP_{q-m}(\L)\cap\mathfrak{n}_{-m}(\g,\L)\\ &=&
\PP_{q-m}(\L)\cap\B_1+\PP_{q_1-m}(\L)\cap\mathfrak{n}_{-m}(\g,\L).
\end{eqnarray*}
On a $\PP_{q-m}(\L)\cap\B_\g\subset\JJ^{m_0}(\g,\L)$ et, puisque
$m_0<\rr_1=\lfloor q_1/2\rfloor+1$, on a~:
\begin{equation*}
\PP_{q_1-m}(\L)\cap\mathfrak{n}_{-m}(\g,\L)\cap\PP_{m_0}
\subset\JJ^{m_0}(\g,\L)
\end{equation*}
par r\'ecurrence. 
Puisque $\JJ^{m_0}(\g,\L)=\JJ^{m_0}(\b,\L)$, la preuve est termin\'ee.
\end{proof}

\begin{lemm}
\label{TransT}
Soit $m_0=\min\{m,s\}$ et soit $\vartheta$ un caract\`ere de $\H^{m_0}$
dont la restriction \`a $\H^{m+1}$ est dans $\Cc(\L,m,\b)$.
Soit $\K$ un sous-groupe ouvert compact de $\U_{m_0}(\L)$, et soit
$\rho$ une re\-pr\'e\-sen\-tation irr\'eductible de $\K$ dont la
restriction \`a $\H^{m_0}\cap\K$ contient $\vartheta_{|\H^{m_0}\cap\K}$.
Alors toute repr\'esentation irr\'eductible de $\G$ contenant
$\vartheta$ contient $\rho$.
\end{lemm}

\begin{proof}
La preuve est similaire \`a celle de \cite[Proposition 8.1.7]{BK}.
Il suf\-fit de remplacer \cite[Lemma 8.1.8]{BK} par le lemme
\ref{IndTau}. 
\end{proof}

\subsection{}
\label{LeSoleilBrilleAlaPlage}

Soit $\pi$ une repr\'esentation irr\'eductible de $\G$ et soit $s$ une 
corestriction mod\'er\'ee sur $\A$ relative \`a $\E/\F$.

\begin{prop}\label{Transfer}
Soient $\tilde\t\in\Cc(\L,\lceil m\rceil-1,\b)$ et
$c\in\PP_{\lceil-m\rceil}(\L)$ tels que $\pi$ contienne le caract\`ere 
$\vartheta=\tilde\t\psi_c$ de $\H^{\lceil m\rceil}(\b,\L)$.
Soit $\L'$ une $\o_{\D}$-suite $\E$-pure de p\'eriode $e(\L')$ et soit
$m'=me(\L')/e(\L)$. 
Soit
$\a'\in\PP_{\lceil-m\rceil}(\L)\cap\PP_{\lceil-m'\rceil}(\L')\cap\B$
tel que~: 
\begin{equation*}
s(c)+\PP_{1-\lceil m\rceil}(\L)\cap\B\subset
\a'+ \PP_{1-\lceil m'\rceil}(\L')\cap\B.
\end{equation*}
Alors il existe $\tilde\t'\in\Cc(\L',\lceil m'\rceil-1,\b)$ et
$c'\in\PP_{\lceil-m'\rceil}(\L')$ tels que $s(c')=\a'$ et tels que 
$\pi$ contienne le caract\`ere $\vartheta'=\tilde\t'\psi_{c'}$ de
$\H^{\lceil m'\rceil}(\b,\L')$. 
Si $\a'=0$, on peut choisir $c'=0$.
\end{prop}

\begin{proof} 
Soit $m_0=\min\{\lceil m\rceil,s\}$. 
Puisque 
$\H^{m_0}(\b,\L)/\H^{\lceil m\rceil}(\b,\L)$ est ab\'elien, $\pi$
contient un caract\`ere $\tilde\vartheta$ de $\H^{m_0}(\b,\L)$ qui
prolonge $\vartheta$. 
En prolongeant $\tilde\t$ en un caract\`ere simple
de $\H^{m_0}(\b,\L)$, que l'on note aussi $\tilde\t$, et en changeant
$c$ dans sa classe modulo $\PP_{1-\lceil m\rceil}(\L)$, on a encore
$\tilde\vartheta=\tilde\t\psi_c$.

Nous d\'emontrons la proposition dans un premier temps sous
l'hypoth\`ese suppl\'ementaire~:
\begin{equation*}
\H^{\lceil m'\rceil}(\b,\L') \subset \U^{m_0}(\L). \leqno{(\H)}
\end{equation*}

On commence par prouver le lemme suivant.

\begin{lemm}
\label{lemmeaajouteri}
Pour tous $k,k'\in\ZZ$, on a~:
\begin{equation*}
s(\PP_k(\L)\cap\PP_{k'}(\L'))=s(\PP_k(\L))\cap s(\PP_{k'}(\L')).
\end{equation*}
\end{lemm}

\begin{proof}
Dans le cas o\`u $\A$ est d\'eploy\'ee sur $\F$ et o\`u $\L$, $\L'$
sont strictes, c'est une cons\'equence de \cite[\S1.3]{BK} et du fait
que $\PP_k(\L)$ et $\PP_{k'}(\L')$ sont des r\'eseaux $\E$-exacts.
Dans le cas g\'en\'eral, on choisit un couple $(\V^{\pt},\L^{\pt})$ 
comme au \S\ref{KarlAdam}.
L'\'egalit\'e est valable pour les $\o_\F$-suites sous-jacentes \`a
$\bar\L$ et \`a $\bar\L'$, appliqu\'ee avec une corestriction
mod\'er\'ee de $\End_\F(\V\oplus\V^{\pt})$ dont la restriction \`a
$\A$ est $s$.
On obtient le lemme \ref{lemmeaajouteri} par projection sur $\A$.
\end{proof}

Puisque $\a'$ appartient \`a 
$\PP_{\lceil-m\rceil}(\L)\cap\PP_{\lceil-m'\rceil}(\L')\cap\B$, 
il existe, d'apr\`es le lemme~\ref{lemmeaajouteri}, un \'el\'ement 
$c'\in\PP_{\lceil-m\rceil}(\L)\cap\PP_{\lceil-m'\rceil}(\L')$ tel que
$s(c')=\a'$, et on peut prendre $c'=0$ si $\a'=0$. 
On pose $\d=c'-c\in\PP_{\lceil-m\rceil}(\L)$.

\begin{lemm}\label{Approx}
Il existe $x\in\PP_{q-\lfloor m\rfloor}(\L)
\cap\mathfrak{n}_{-\lfloor m\rfloor}(\b,\L)$ tel que~:
\begin{equation*}
\d-(1+x)^{-1}a_{\b}(x)-(1+x)^{-1}(cx-xc)\in\PP_{1-\lceil m'\rceil}(\L').
\end{equation*}
\end{lemm}

\begin{proof} 
On d\'emontre par r\'ecurrence que, pour $t\>0$, il existe 
$x_t\in\PP_{q-\lfloor m\rfloor}(\L)
\cap\mathfrak{n}_{-\lfloor m\rfloor}(\b,\L)$ tel que~:
\begin{equation}\label{Inductivestep}
\d-(1+x_t)^{-1}a_{\b}(x_t)-(1+x_t)^{-1}(cx_t-x_tc)\in
\PP_{1-\lceil m'\rceil}(\L')+\PP_{t-\lfloor m\rfloor}(\L).
\end{equation}
Le lemme s'ensuit puisque 
$\PP_{t-\lfloor m\rfloor}(\L)\subset\PP_{1-\lceil m'\rceil}(\L')$ pour
$t$ suffisamment grand.
Puisque $\d\in\PP_{-\lfloor m\rfloor}(\L)$,
on peut prendre $x_0=0$. 
Supposons donc que $t\>0$ et qu'on ait
trouv\'e $x_t\in\PP_{q-\lfloor m\rfloor}(\L)\cap
\mathfrak{n}_{-\lfloor m\rfloor}(\b,\L)$ tel que
\eqref{Inductivestep} soit satisfaite. 
Il existe alors $\d_t\in\PP_{t-\lfloor m\rfloor}(\L)$ tel que~:
\begin{equation*}
\d-(1+x_t)^{-1}a_{\b}(x_t)-(1+x)^{-1}(cx_t-x_tc)\in
\d_t+\PP_{1-\lceil m'\rceil}(\L').
\end{equation*}
On impose $\d_0=\d$ dans le cas $t=0$ et, dans ce cas~:
\begin{equation*}
s(\d_0)=\a'-s(c)\in\PP_{1-\lceil m'\rceil}(\L')\cap\B.
\end{equation*}
\'Egalement, pour $t>0$, on a~:
\begin{equation*}
s(\d_t)\in\PP_{t-\lfloor m\rfloor}(\L)\cap\B
\subset \PP_{1-\lceil m\rceil}(\L)\cap\B
\subset \PP_{1-\lceil m'\rceil}(\L')\cap\B.
\end{equation*}
Donc $\d_t\in\PP_{t-\lfloor m\rfloor}(\L)\cap
\left(\PP_{1-\lceil m'\rceil}(\L')+a_{\b}(\A)\right)$.

\begin{lemm}
\label{lemmeaajouterii}
Pour tous $k,k'\in\ZZ$, on a~:
\begin{equation*}
\PP_k(\L)\cap\left(\PP_{k'}(\L')+a_\b(\A)\right)=
\PP_k(\L)\cap\PP_{k'}(\L')+\PP_k(\L)\cap a_\b(\A).
\end{equation*}
\end{lemm}

\begin{proof}
Dans le cas o\`u $\A$ est d\'eploy\'ee sur $\F$ et o\`u $\L$, $\L'$
sont strictes, c'est une cons\'equence de \cite[8.1.13]{BK}.
On traite le cas g\'en\'eral comme au lemme \ref{lemmeaajouteri}.
\end{proof}

D'apr\`es le lemme~\ref{lemmeaajouterii}, on a~:
\begin{equation*}
\PP_{t-\lfloor m\rfloor}(\L)\cap
\left(\PP_{1-\lceil m'\rceil}(\L')+a_{\b}(\A)\right)
= \PP_{t-\lfloor m\rfloor}(\L)\cap\PP_{1-\lceil m'\rceil}(\L')
+ \PP_{t-\lfloor m\rfloor}(\L)\cap a_{\b}(\A).
\end{equation*}
D'apr\`es la proposition \ref{Cor1410}, il existe 
$y_t\in \PP_{q+t-\lfloor m\rfloor}(\L)\cap
\mathfrak{n}_{t-\lfloor m\rfloor}(\b,\L)$ tel que~:
\begin{equation*}
\d_t-a_{\b}(y_t) \in
\PP_{t-\lfloor m\rfloor}(\L)\cap\PP_{1-\lceil m'\rceil}(\L'),
\end{equation*}
et $x_{t+1}=x_t+y_t$ est comme il faut. 
\end{proof}

Revenons \`a la d\'emonstration de la proposition \ref{Transfer}
sous l'hypoth\`ese $(\H)$.
Soit $x$ comme dans le lemme \ref{Approx}.
Par le lemme~\ref{QuelleHorreurInsoutenableALInfiniRepetee},
l'\'el\'ement $1+x$ normalise $\H^{m_0}(\b,\L)$ et~:
\begin{equation*}
\tilde\vartheta^{1+x}=\tilde\vartheta \psi_{(1+x)^{-1}a_{\b}(x)}
\psi_{(1+x)^{-1}(cx-xc)}.
\end{equation*}
Soit $\tilde\t'$ le caract\'ere simple dans
$\Cc(\L',\lceil m'\rceil-1,\b)$ qui co\"incide avec $\tilde\t$ sur
$\H^{\lceil m'\rceil}(\b,\L')\cap\H^{m_0}(\b,\L)$. 
Le lemme \ref{Approx} implique que, comme caract\`ere du groupe 
$\H^{\lceil m'\rceil}(\b,\L')\cap\H^{m_0}(\b,\L)$, on a~:
\begin{equation*}
\tilde\vartheta^{1+x} = \tilde\t'\psi_{c'}.
\end{equation*}
Puisque l'hypoth\`ese $(\H)$ est satisfaite, on
d\'eduit du lemme \ref{TransT} que $\pi$ contient
$\vartheta'=\tilde\t'\psi_{c'}$.

Traitons maintenant le cas g\'en\'eral. 
Supposons d'abord que $m\<q/2$.
D'apr\`es le lemme \ref{FinsFaciles}, il existe une famille
finie $(\L_i)_{0\<i\<l}$ v\'erifiant~:
\begin{equation*}
\U_{\lceil m_{i}\rceil}(\L_{i})
\supset
\U_{\lfloor m_{{i+1}}\rfloor+1}(\L_{{i+1}}),
\quad
0\<i\<l-1,
\end{equation*}
avec $(\L_0,m_0)=(\L,m)$, $(\L_l,m_l)=(\L',m')$ et $m_i=me_i/e$,
o\`u $e_i$ est la p\'eriode de $\L_i$.
Par le lemme~\ref{cotcot}, on a
$s(c)\in\PP_{-\lfloor m_i\rfloor}(\L_i)\cap\B$, pour chaque $i$, 
et~:
\begin{equation*}
s(c)+\PP_{1-\lceil m_i\rceil}(\L_i)\cap\B \subset
s(c)+\PP_{1-\lceil m_{i+1}\rceil}(\L_{i+1})\cap\B,
\quad
0\<i\<l-1.
\end{equation*}
Par dualit\'e, on obtient~:
\begin{equation*}
\U_{\lceil m_{i}\rceil}(\L_{i})\cap\B
\supset
\U_{\lceil m_{{i+1}}\rceil}(\L_{{i+1}})\cap\B.
\end{equation*}
Soit $q_i=-k_0(\b,\L_i)=qe_i/e$. 
Puisque $m\<q/2$, on a $m_i\< q_i/2$ et donc~:
\begin{equation*}
\H^{\lceil m_{{i+1}}\rceil}(\b,\L_{{i+1}})
=
\left(\U_{\lceil m_{{i+1}}\rceil}(\L_{{i+1}})\cap\B\right)
\H^{\lfloor m_{{i+1}}\rfloor+1}(\b,\L_{{i+1}})
\subset
\U_{\lceil m_{i}\rceil}(\L_{i}).
\end{equation*}
Appliquant le cas o\`u l'hypoth\`ese $(\H)$ est v\'erifi\'ee, on
voit que, pour chaque $i$, il existe
$\tilde\t_i\in\Cc(\L_i,\lceil m_i\rceil-1,\b)$ et
$c_i\in\PP_{\lceil-m_i\rceil}(\L_i)$ tels que $s(c_i)=s(c)$ et
$\pi$ contienne le caract\`ere $\vartheta_i=\tilde\t_i\psi_{c_i}$ de
$\H^{\lceil m_i\rceil}(\b,\L_i)$. 
\`A la derni\`ere \'etape, on peut remplacer $s(c)$ par $\a'$ et on 
en d\'eduit le r\'esultat. 

Finalement, supposons que $m>q/2$ et fixons $\varepsilon>0$ tel que
$e\varepsilon<m-q/2$. 
Par le lemme~\ref{FinsFacilesEpsilon},
il existe une famille finie $(\L_i)_{0\<i\<l}$ v\'erifiant~:
\begin{equation*}
\U_{\lceil m_{i+1}\rceil}(\L_{i+1})
\subset
\U_{\lceil m_{i}-e_i\varepsilon\rceil}(\L_{i}),
\quad
0\<i\<l-1,
\end{equation*}
avec $(\L_0,m_0)=(\L,m)$, $(\L_l,m_l)=(\L',m')$ et $m_i=me_i/e$,
o\`u $e_i$ est la p\'eriode de $\L_i$. 
On a~:
\begin{equation*}
m_i-e_i\varepsilon = \frac{e_i}e (m-e\varepsilon) >
\frac{e_i}e \frac q2 = \frac{q_i}2,
\end{equation*}
donc $\lceil m_{i}-e_i\varepsilon\rceil\> \lceil q_i/2\rceil=s_i$ et~:
\begin{equation*}
\H^{\lceil m_{{i+1}}\rceil}(\b,\L_{{i+1}})
\subset
\U_{\lceil m_{{i+1}}\rceil}(\L_{{i+1}})
\subset
\U_{\lceil m_{i}-e_i\varepsilon\rceil}(\L_{i})
\subset
\U_{\lceil s_i\rceil}(\L_{i}).
\end{equation*}
La d\'emonstration se termine maintenant comme dans le cas pr\'ec\'edent.
\end{proof}

\subsection{}
\label{SuperSection}

Soit $\pi$ une repr\'esentation irr\'eductible de niveau non nul de $\G$.

\begin{prop}
\label{Dicho}
On est dans l'un des deux cas suivants~:
\begin{enumerate}
\item
Il existe une strate scind\'ee $[\L,n,n-1,b]$ de $\A$, 
avec $n\>1$ et $\L$ stricte, telle que la restriction de $\pi$ \`a
$\U_n(\L)$ contienne $\psi_b$ ;
\item
Il existe une strate simple $[\L,n,m,\b]$ de $\A$,
avec $\L$ stricte, et un ca\-rac\-t\`ere simple $\t\in\Cc(\L,m,\b)$, 
tels que la restriction de $\pi$ \`a $\H^{m+1}(\b,\L)$ contienne
$\t$.
\end{enumerate}
\end{prop}

\begin{rema}
Dans le cas (1), on dit que $\pi$ contient une strate scind\'ee.
\end{rema}

\begin{proof}
D'abord, d'apr\`es \cite[Theorem 1.2.1(i)]{Br4}, la repr\'esentation
$\pi$ contient une strate fondamentale $[\L,n,n-1,b]$ de $\A$, avec
$\L$ stricte, \ie que la restriction de $\pi$ \`a $\U_n(\L)$ contient
$\psi_b$. 
Si elle est scind\'ee, alors on est dans le premier cas.
Si elle ne l'est pas, alors, d'apr\`es \cite[Theorem 1.2.4]{Br4},
la repr\'esentation $\pi$ contient une strate simple $[\L,n,m,\b]$ de
$\A$, avec $\L$ stricte, \ie que la restriction de $\pi$ \`a
$\H^{m+1}(\b,\L)$ contient un caract\`ere simple $\t\in\Cc(\L,m,\b)$.
\end{proof}

Il est commode d'introduire la d\'efinition suivante.
Soit $[\L,n,m,\b]$ une strate simple de $\A$ avec $m\>1$. 
Soit $\V_{\E}$ un $\B$-module \`a gauche simple, 
soit $\D_\E$ l'alg\`ebre oppos\'ee \`a $\End_{\B}(\V_\E)$,
soit $\Ga$ une $\o_{\D_{\E}}$-suite de r\'eseaux de $\V_{\E}$ 
v\'erifiant la condition du th\'eor\`eme \ref{CNNBr1} et
soit $s$ une corestriction mod\'er\'ee sur $\A$ 
relativement \`a $\E/\F$ ({\it cf.} \S\ref{CorMod}).

\begin{defi}
\label{StrateDerivee!}
Une {\it strate d\'eriv\'ee} de $[\L,n,m,\b]$ est une strate de $\B$
de la forme $[\Ga,m,m-1,s(c)]$ avec $c\in\PP_{-m}(\L)$.
\end{defi}

On d\'efinit maintenant les caract\`eres scind\'es de $\G$.
C'est ce qui est appel\'e, dans \cite{BK}, les types scind\'es 
de niveau $(x,y)$, avec $x>y>0$. 

\begin{defi}
\label{CarScinde}
Un couple $(\K,\vartheta)$ est un {\it caract\`ere scind\'e} de $\G$
s'il existe une strate simple $[\L,n,m,\b]$ de $\A$,
avec $m\>1$ et $\L$ stricte, un caract\`ere simple
$\t\in\Cc(\L,m-1,\b)$ et $c\in\PP_{-m}(\L)$ tels que~:
\begin{itemize}
\item[(i)]
On a $\K=\H^{m}(\b,\L)$ et $\vartheta=\t\psi_c$.
\item[(ii)]
La strate d\'eriv\'ee $[\Ga,m,m-1,s(c)]$ est scind\'ee, pour n'importe
quelle corestriction mod\'er\'ee $s$ sur $\A$ 
relativement \`a $\E/\F$.
\end{itemize}
\end{defi}

La distinction entre strate scind\'ee et caract\`ere scind\'e est
assez superficielle~: c'est \`a peu pr\`es la m\^eme que celle qu'on 
fait entre types simples de niveau $0$ et de niveau $>0$ 
({\it cf.} \cite{VS3}).

On est maintenant en mesure de formuler le r\'esultat principal de
cette section.

\begin{theo}
\label{??}
Soit $\pi$ une repr\'esentation irr\'eductible de niveau non nul de
$\G$.
Alors~: 
\begin{enumerate}
\item
ou bien il existe une strate simple $[\L,n,0,\b]$ de $\A$,
avec $\L$ stricte, telle que $\pi$ contienne un caract\`ere 
simple $\t\in\Cc(\L,0,\b)$ ;
\item
ou bien $\pi$ contient une strate scind\'ee.
\item
ou bien $\pi$ contient un caract\`ere scind\'e.
\end{enumerate}
\end{theo}

Le reste du paragraphe est consacr\'e \`a la d\'emonstration du th\'eor\`eme
\ref{??}. 
D'apr\`es la proposition \ref{Dicho}, si $\pi$ ne contient pas
de strate scind\'ee, il existe une strate simple $[\L,n,m,\b]$
de $\A$, avec $\L$ stricte, et un caract\`ere simple
$\t\in\Cc(\L,m,\b)$, tels que la restriction de $\pi$ \`a
$\H^{m+1}(\b,\L)$ contienne $\t$.
On choisit $[\L,n,m,\b]$ et $\t$ tels que le rapport $m/e(\L|\o_{\F})$
soit minimal.
Si $m=0$, alors on est dans le cas $(1)$ du th\'eor\`eme \ref{??}.
Dans toute la suite du \S, on suppose que $m\>1$. 
On fixe un $\B$-module \`a gauche simple $\V_{\E}$ et
on note $\D_\E$ l'alg\`ebre oppos\'ee \`a $\End_{\B}(\V_\E)$.
On fixe une $\o_{\D_{\E}}$-suite $\Ga$ de r\'eseaux de $\V_{\E}$ 
v\'erifiant la condition du th\'eor\`eme \ref{CNNBr1}.
On fixe un caract\`ere $\vartheta$ de $\H^{m}(\b,\L)$ contenu 
dans $\pi$ et pro\-longeant $\t$, un caract\`ere simple  
$\tilde\t\in\Cc(\L,m-1,\b)$ prolongeant $\t$ et un \'el\'e\-ment 
$c\in\PP_{-m}(\L)$ tel que $\vartheta=\tilde\t\psi_c$.
Enfin, on fixe une corestriction mod\'er\'ee $s$ sur $\A$ 
relativement \`a $\E/\F$.
Il s'agit de prouver que la strate d\'eriv\'ee $[\Ga,m,m-1,s(c)]$ 
est scind\'ee.

\begin{prop}
\label{THderFONDA}
La strate $[\Ga,m,m-1,s(c)]$ est fondamendale.
\end{prop}

\begin{rema}
Pour simplifier les notations, on calcule toutes les p\'e\-rio\-des 
de suites de r\'eseaux sur $\o_\F$ et, si $\L$ est une suite de
r\'eseaux, on note $e(\L)$ pour $e(\L|\o_{\F})$.
Cette remarque se substitue donc \`a la remarque \ref{WARNING}.
Pour une $\o_{\D}$-suite de r\'eseaux $\L$, on a 
$e(\L|\o_{\F})=e(\L|\o_{\D})d$, o\`u $d$ d\'esigne le
degr\'e r\'eduit de $\D$ sur $\F$.
Pour une $\o_{\D_\E}$-suite de r\'eseaux $\Ga$, on a 
$e(\Ga|\o_{\F})=e(\Ga|\o_{\E})e_{\E/\F}$, o\`u $e_{\E/\F}$ d\'esigne
l'indice de ramification de $\E/\F$.
\end{rema}

\begin{proof}
On raisonne par l'absurde.
D'apr\`es la proposition \ref{RafNonFonda}, il existe 
un entier $k'\in\ZZ$ et une $\o_{\D_\E}$-suite stricte 
$\Ga'$ de $\V_{\E}$ tels que~:
\begin{equation}
\label{Arenaire1}
s(c)+\PP_{1-m}(\L)\cap\B
=s(c)+\PB_{1-m}(\Ga)
\subset\PB_{-k'}(\Ga')
\end{equation}
et~:
\begin{equation}
\label{Arenaire2}
\frac{k'}{e(\Ga')}<\frac{m}{e(\Ga)}.
\end{equation}
On fixe un couple $(\epsilo',\L')$ correspondant \`a $\Ga'$ par le
th\'eor\`eme \ref{Edescent} et on pose $m'=me(\L')/e(\L)$. 

\begin{lemm}
On a~:
\begin{equation}
\label{Arenaire3}
s(c)+\PP_{1-m}(\L)\cap\B\subset\PP_{1-\lceil m'\rceil}(\L')\cap\B
\end{equation}
et~:
\begin{equation}
\label{Arenaire4}
\frac{\lceil m'\rceil-1}{e(\L')}<\frac{m}{e(\L)}.
\end{equation}
\end{lemm}

\begin{proof}
\`A partir de (\ref{Arenaire2}) et de la remarque
\ref{Tetragrammaton}, on \'ecrit~:
\begin{equation*}
\epsilo'k'<\frac{e(\L')}{e(\L)}m=m'.
\end{equation*}
Compte tenu du fait que $\epsilo'k'$ est entier, on en d\'eduit~:
\begin{equation*}
-\epsilo'k'\>1-\lceil m'\rceil,
\end{equation*}
ce qui, avec (\ref{Arenaire1}), donne (\ref{Arenaire3}).
Ensuite, on \'ecrit $m'=\lceil m'\rceil-1+a/e(\L)$, avec $1\<a\<e(\L)$.
On obtient~:
\begin{equation*}
\frac{\lceil m'\rceil-1}{e(\L')}
=\frac{m'}{e(\L')}-\frac{a}{e(\L)e(\L')}
<\frac{m}{e(\L)},
\end{equation*}
ce qui termine la d\'emonstration.
\end{proof}

D'apr\`es le lemme \ref{Transfer} appliqu\'e avec $\a'=0$, 
il existe un caract\`ere simple de $\Cc(\L',\lceil m'\rceil-1,\b)$ 
contenu dans $\pi$, ce qui contredit la minimalit\'e de $m/e(\L)$.
\end{proof}

\begin{prop}
\label{THder}
La strate $[\Ga,m,m-1,s(c)]$ est scind\'ee.
\end{prop}

\begin{proof}
On raisonne par l'absurde.
D'apr\`es la proposition \ref{RafSimlple}, il existe un 
entier $k'\in\ZZ$, une $\o_{\D_{\E}}$-suite stricte $\Ga'$ de
$\V_{\E}$ et $\a'\in\B$ tels que~:
\begin{equation*}
s(c)+\PP_{1-m}(\L)\cap\B\subset\a'+\PB_{1-k'}(\Ga')
\end{equation*}
et tels que la strate $[\Ga',k',k'-1,\a']$ soit simple.
On fixe un couple $(\epsilo',\L')$ correspondant \`a $\Ga'$ par le
th\'eor\`eme \ref{Edescent} et on pose $m'=\epsilo'k'\in\ZZ$.
On a donc~:
\begin{equation*}
s(c)+\PP_{1-m}(\L)\cap\B\subset\a'+\PP_{1-m'}(\L')\cap\B.
\end{equation*}
Avant d'appliquer \`a nouveau le lemme \ref{Transfer}, on a besoin des 
deux lemmes suivants.

\begin{lemm}
\label{con1}
On a $m/e(\L)=m'/e(\L')$.
\end{lemm}

\begin{proof}
Puisque la strate $[\Ga,m,m-1,s(c)]$ est fondamentale, 
on a $s(c)\in\KK(\Ga)$.
Ensuite, on a $s(c)-\a'\in\PP_{1-k'}(\Ga')$ et
$\a'\in\PP_{-k'}(\Ga')$, \ie que les strates 
$[\Ga',k',k'-1,s(c)]$ et $[\Ga',k',k'-1,\a']$
sont \'equivalentes.
Elles ont donc le m\^eme polyn\^ome caract\'eristique, 
de sorte que $[\Ga',k',k'-1,s(c)]$ est fondamentale.
Ainsi $s(c)\in\KK(\Ga')$.
On en d\'eduit que $m/e(\Ga)=k'/e(\Ga')$,
ce qui implique l'\'egalit\'e voulue.
\end{proof}

\begin{lemm}
\label{con2}
On a $\PP_{-m'}(\L')\cap\B\subset\PP_{-m}(\L)\cap\B$.
\end{lemm}

\begin{proof}
En prenant le dual de l'inclusion~:
\begin{equation*}
\PP_{1-m}(\L)\cap\B\subset\PP_{1-m'}(\L')\cap\B,
\end{equation*}
on obtient $\PP_{m'}(\L')\cap\B\subset\PP_{m}(\L)\cap\B$.
Puis, en multipliant par $s(c)^2$ et en tenant compte du lemme
\ref{con1}, on obtient l'inclusion voulue.
\end{proof}

D'apr\`es le lemme \ref{Transfer}, il existe un 
caract\`ere simple $\tilde\t'\in\Cc(\L',m'-1,\b)$ et
$c'\in\PP_{-m'}(\L')$ tels que $s(c')=\a'$ et que
$\pi$ contienne le caract\`ere $\vartheta'=\tilde\t'\psi_{c'}$.
Pour terminer la preuve de la proposition \ref{THder}, on a besoin 
du r\'esultat suivant, analogue de \cite[Theorem 2.2.8]{BK}.

\begin{prop}
\label{TheoremeEnSuspens}
La strate $[\L',n',m'-1,\b+c']$ est \'equivalente \`a une strate
simple.
\end{prop}

On reporte la preuve de la proposition \ref{TheoremeEnSuspens} au
paragraphe suivant.
En at\-ten\-dant, on termine la preuve de la proposition \ref{THder}.
On choisit une strate simple $[\L',n',m'-1,\b']$ \'equivalente \`a
$[\L',n',m'-1,\b+c']$.
La strate $[\L',n',m',\b']$ est \'equivalente \`a la strate
simple $[\L',n',m',\b]$.
On applique la proposition \ref{GrosBoxon1L}.
On a une bijection~:
\begin{equation*}
\Cc(\L',m'-1,\b)\f\Cc(\L',m'-1,\b')
\end{equation*}
envoyant $\tilde\t'$ sur $\vartheta'\psi_{\b'-\b-c'}$.
Mais $\b'-\b-c'\in\PP_{1-m'}(\L')$.
Ceci implique que $\vartheta'\in\Cc(\L',m'-1,\b')$, et contredit  
la minimalit\'e de $m/e(\L|\o_\F)$.
\end{proof}

\subsection{}

Dans ce paragraphe, on d\'emontre la proposition
\ref{TheoremeEnSuspens}. 
D'apr\`es \cite[Th\'eo\-r\`eme 2.2]{VS3}, il suffit de prouver que 
$[\L',n',m'-1,\b+c']$ est \'equivalente \`a une strate pure.

On pose $\K=\E(\a')$ et on note $\C$ le commutant de $\K$ dans $\A$.
On fixe un $\K\otimes_\F\D$-module \`a droite simple $\SS$ et on pose
$\A(\K)=\End_\D(\SS)$. 
On note $\D_\K$ le commutant de $\K$ dans $\A(\K)$.
C'est une $\K$-alg\`ebre \`a division.
On choisit une d\'ecomposition de $\V$ en somme de
$\K\otimes_\F\D$-modules qui soit conforme \`a $\L'$.
On en d\'eduit un plongement de $\F$-alg\`ebres $\iota:\A(\K)\f\A$ 
et un isomorphisme de $(\A(\K),\C)$-bimodules 
({\it cf.} \cite[\S1.3]{VS3})~:
\begin{equation}
\label{BonVieuxTempsDesWEDecompositions}
\A(\K)\otimes_{\D_\K}\C\f\A.
\end{equation}

On a besoin du lemme suivant.

\begin{lemm}
Soit $s'$ une corestriction mod\'er\'ee sur $\A(\K)$ relativement
\`a $\E/\F$.
Alors $s'\otimes{\rm id}_{\C}$ est une corestriction mod\'er\'ee sur
$\A$ relativement \`a $\E/\F$.
\end{lemm}

\begin{proof}
La preuve est analogue \`a celle de \cite[Proposition 1.3.9]{BK},
compte tenu de \cite[Lemmas 4.2.1--4.2.2]{Br4}.
\end{proof}

On note $s'$ la corestriction mod\'er\'ee sur $\A(\K)$ relativement
\`a $\E/\F$ telle que $s'\otimes{\rm id}_{\C}$ corresponde
\`a $s$ {\it via} (\ref{BonVieuxTempsDesWEDecompositions}).
On note $\AA(\K)$ l'unique ordre h\'er\'editaire de $\A(\K)$
normalis\'e par $\mult\K$ et $\mathfrak{P}(\K)$ son radical de
Jacobson.
On note $e$ le rapport de $e(\L'|\Oo_\D)$ sur $e(\AA(\K)|\Oo_\D)$.
On pose $n''=n'/e$ et $m''=m'/e$.
Ce sont des entiers, \'egaux respectivement \`a $-\v_{\AA(\K)}(\b)$ et
\`a $-\v_{\AA(\K)}(\a')$.
Soit enfin un \'el\'ement $c_0\in\mathfrak{P}(\K)^{-m''}$ tel que
$s'(c_0)=\a'$.

\begin{lemm}
\label{PlaPla}
La strate $[\AA(\K),n'',m''-1,\b+c_0]$ est pure.
\end{lemm}

\begin{proof}
La preuve est tr\`es proche de celle de \cite[Proposition 2.2.3]{BK}.
D'apr\`es le th\'eor\`eme \ref{CNNBr1}, le normalisateur de $\AA(\K)$
dans $\D_\K$ est \'egal au normalisateur de l'unique $\o_{\K}$-ordre
de $\D_\K$, qui est $\o_{\D_\K}$. 
Donc $\AA(\K)$ est normalis\'e par $\mult\D_{\K}$.
On note $\B(\K)$ le commutant de $\E$ dans $\A(\K)$. 

Soit $x\in\A(\K)^{\times}$ commutant \`a $\b+c_0$ et soit $t\in\ZZ$ le
plus grand entier tel que $x\in\mathfrak{P}(\K)^{t}$.
En raisonnant comme dans le preuve de \cite[Proposition 2.2.3]{BK} et
en rempla\c cant \cite[Corollary 1.4.10]{BK} par la proposition
\ref{Cor1410}, on obtient~:
\begin{equation*}
x\in\left(\mathfrak{P}(\K)^{t}\cap\D_\K+\mathfrak{P}(\K)^{t+1}\right)
\backslash\mathfrak{P}(\K)^{t+1}\subset\KK(\AA(\K)).
\end{equation*}
Donc le centralisateur de $\F[\b+c_0]$ dans $\A(\K)^{\times}$ est compact
modulo le centre.
Ainsi la $\F$-alg\`ebre $\F[\b+c_0]$ est un corps dont le groupe
multiplicatif est contenu dans $\KK(\AA(\K))$, ce qui termine la
d\'emonstration du lemme \ref{PlaPla}.
\end{proof}

Ainsi l'image par $\iota$ de la strate pure
$[\AA(\K),n'',m''-1,\b+c_0]$ est une strate pure
$[\L',n',m'-1,\b+c_0]$.
On a $s(\b+c')=s(\b+c_0)$. 
D'apr\`es la proposition \ref{Cor1410}, il existe un \'el\'ement 
$y\in\PP_{q'-m'}(\L)\cap\mathfrak{n}_{-m'}(\L)$
tel que~:
\begin{equation*}
c'-c_0\equiv a_\b(y)\mod{\PP_{1-m'}(\L')}.
\end{equation*}
On en d\'eduit que $[\L',n',m'-1,\b+c']$ est \'equivalente \`a la 
conjugu\'ee de la strate $[\L',n',m'-1,\b+c_0]$ par $1+y$, ce qui 
met fin \`a la fois \`a la preuve de la proposition
\ref{TheoremeEnSuspens} et \`a celle du th\'eor\`eme \ref{??}.


\section{Modules de Jacquet}
\label{AymeJacquet}

Soit $\pi$ une representation irr\'eductible de niveau non nul de
$\G$.
Par le th\'eor\`eme~\ref{??}, on sait que $\pi$ contient soit un
caract\`ere simple d'un groupe $\H^1(\b,\L)$ avec $\L$ stricte, 
soit un caract\`ere scind\'e, soit une strate scind\'ee.
Le but principal de cette section est de d\'emontrer le th\'eor\`eme
suivant~:

\begin{theo}\label{scsimple}
Soit $\pi$ une representation irr\'eductible supercuspidale
de niveau non nul de $\G$.
Il existe une strate simple $[\L,n,0,\b]$, avec $\L$ stricte, 
et un caract\`ere simple $\t\in\Cc(\L,0,\b)$ tels que 
$\pi_{|\H^1(\b,\L)}$ contienne $\t$. 
\end{theo}

L'id\'ee est d'utiliser la notion de paire couvrante pour d\'emontrer
que, dans le cas o\`u la repr\'esentation $\pi$ contient un
caract\`ere scind\'e ou une strate scind\'ee, elle a un module de
Jacquet non nul.

\subsection{}
\label{S41}

Soit $[\L,n,m,\b]$ une strate simple de $\A$ avec $m\>1$ et avec $\L$
stricte. 
Soit $\V=\V^1\oplus \V^2$ une d\'ecomposition de $\V$ en
$\E\otimes_\F\D$-modules, qui est conforme \`a $\L$, et soit~:
\begin{equation*}
\M=\Aut_\D(\V^1)\times\Aut_\D(\V^2),
\end{equation*}
qui est un sous-groupe de Levi de $\G$.
Soient $\N=1+\A^{12}$ et $\N^-=1+\A^{21}$, et soient $\P=\M\N$ et
$\P^-=\M\N^-$.
Donc $\P$ est un sous-groupe parabolique avec facteur de Levi $\M$,
et $\P^-$ est le sous-groupe parabolique oppos\'e.

Soit $\B$ le commutant de $\E$ dans $\A$. 
On fixe un $\B$-module \`a gauche simple $\V_{\E}$ et on note
$\D_{\E}$ l'alg\`ebre oppos\'ee \`a $\End_{\B}(\V_{\E})$. 
Soit $\Ga$ une $\Oo_{\D_\E}$-suite de r\'eseaux telle que
$\PP_k(\L)\cap\B=\PB_k(\Ga)$ pour $k\in\ZZ$, donn\'ee par le
th\'eor\`eme~\ref{CNNBr1}. 
La d\'ecomposition $\V=\V^1\oplus\V^2$ correspond \`a une
d\'ecomposition $\V_\E=\V_\E^1\oplus\V_\E^2$ telle que~:
\begin{equation*}
\M\cap\B=\Aut_{\D_\E}(\V^1_\E)\times\Aut_{\D_\E}(\V^2_\E).
\end{equation*}
C'est une d\'ecomposition conforme \`a $\Ga$, et on pose 
$\Ga^i=\Ga\cap\V^i_\E$.
On fixe aussi une corestriction mod\'er\'ee $s:\A\to\B$. 
Alors les restrictions $s_i=s_{|\A^i}:\A^i\to\B^i$ sont aussi des
corestrictions mod\'er\'ees.

Soient $c_i\in \A^{i}\cap\PP_{-m}(\L)$.
On pose $c=c_1+c_2\in \A$.
On suppose que la strate $[\Ga,m,m-1,s(c)]$ de $\B$ est (fondamentale)
scin\-d\'ee par la d\'ecomposition, \ie que les polyn\^omes
caracteristiques des strates $[\Ga^i,m,m-1,s(c_i)]$, pour $i\in\{1,2\}$, 
sont premier entre eux.
On suppose aussi que $s(c_1)$ normalise $\Ga^1$ et que
$\v_{\Ga^1}(s(c_1))=-m$. 

\begin{lemm}\label{bornescinde}
On a~:
\begin{equation*}
\I_{\B^\times}({\psi_c}_{|\U_m(\L)\cap\B})\subset
(\U_1(\L)\cap\B)\cdot(\M\cap\B)\cdot(\U_1(\L)\cap\B).
\end{equation*}
\end{lemm}

\begin{proof}
L'entrelacement de ${\psi_c}_{|\U_m(\L)\cap\B}$ est le m\^eme que
l'en\-tre\-la\-ce\-ment de la strate $[\Ga,m,m-1,s(c)]$.
La d\'emonstration est alors presque identique \`a celle
de~\cite[Theorem 4.9]{St1} : il suffit de 
remplacer~\cite[Lemma 4.11]{St1} par~\cite[Lemma 2.3.4]{Br4}.
\end{proof}

Soit maintenant $\t$ un caract\`ere simple dans $\Cc(\L,m-1,\b)$. 
On consid\`ere le caract\`ere $\xi=\t\psi_c$ de $\H^m(\b,\L)$.
On pose~: 
\begin{equation*}
\Om=(\U_1(\L)\cap\B)\Om_{q-m+1}(\b,\L).
\end{equation*}
D'apr\`es la proposition \ref{normBKSSTcs}, le caract\`ere $\t$ est
normalis\'e par $\Om$. 
Puisque $\Om\subset\U_1(\L)$ normalise $\psi_c$, le caract\`ere $\xi$
est lui aussi normalis\'e par $\Om$.

\begin{theo}\label{entrelacementscinde}
On a $\I_\G(\xi)\subset\Om\M\Om$.
\end{theo}

\begin{proof}
Si $g\in\G$ entrelace $\xi$, alors $g$ entrelace {\it a fortiori}
sa restriction $\xi_{|\H^{m+1}(\b,\L)}=\t_{|\H^{m+1}(\b,\L)}$. 
D'apr\`es le th\'eor\`eme \ref{entrelacementCSsuites}, l'\'el\'ement
$g$ appartient donc \`a $\Om_{q-m}(\b,\L)\B^\times\Om_{q-m}(\b,\L)$.

On \'ecrit $\g=(1+x)t(1+y)^{-1}$, avec $x,y\in\MM_{q-m}(\b,\L)$ et
$t\in\B^\times$. 
Puisque $1+x$ est dans
$(1+\PP_{q-m}(\L)\cap\mathfrak{n}_{-m}(\b,\L))\J^s(\b,\L)$, on a, 
d'apr\`es la proposition 
\ref{QuelleHorreurInsoutenableALInfiniRepetee}~:
\begin{equation}
\label{conjugx}
\xi^{1+x} = \xi\psi_{(1+x)^{-1}\b(1+x)-\b} =\xi\psi_{a_\b(x)}.
\end{equation}
L'\'el\'ement $t$ entrelace $\xi^{1+x}$ et $\xi^{1+y}$, donc
leur restriction au groupe
$\H^m(\b,\L)\cap\B=\U_m(\L)\cap\B$.
Puisque les
restriction de $\psi_{a_\b(x)}$ et $\psi_{a_\b(y)}$
\`a $\U_m(\L)\cap\B$ sont triviales, et
puisque $t$ entrelace certainement $\t_{|\U_m(\L)\cap\B}$, on voit que
$t$ entrelace aussi le caract\`ere ${\psi_c}_{|\U_m(\L)\cap\B}$.
D'apr\`es le lemme~\ref{bornescinde}, on a~:
\begin{equation*}
t\in (\U_1(\L)\cap\B) (\M\cap\B) (\U_1(\L)\cap\B).
\end{equation*}
Puisque $\U_1(\L)\cap\B$ normalise $\xi$ et le groupe
$\Om_{q-m}(\b,\L)$, on peut donc supposer que $t\in\M\cap\B$.

Par la d\'ecomposition d'Iwahori de $\Om_{q-m}(\b,\L)$ par rapport \`a
$(\M,\P)$, on \'ecrit $1+x=n_x^- m_x n_x$, avec
$n_x^-\in\Om_{q-m}(\b,\L)\cap\N^-$, $m_x\in\Om_{q-m}(\b,\L)\cap\M$
et $n_x=1+x_n\in\Om_{q-m}(\b,\L)\cap\N$. 
On \'ecrit $1+y=n_y^- m_yn_y$ de la m\^eme mani\`ere.

On pose~:
\begin{equation*}
\HH_- = \begin{pmatrix} \HH^{m+1}(\b,\L) & \HH^{m+1}(\b,\L) \\
\HH^{m}(\b,\L) & \HH^{m+1}(\b,\L) \end{pmatrix}
\end{equation*}
et $\H_-=1+\HH_-=\H^{m+1}(\b,\L)(\H^m(\b,\L)\cap\N^-)$.
L'\'el\'ement $t$ entrelace les restrictions ${\xi^{1+x}}_{|\H_-}$ et
${\xi^{1+y}}_{|\H_-}$ et, puisque $n_x^- m_x$ normalise $\xi_{|\H^-}$,
on a l'\'egalit\'e ${\xi^{1+x}}_{|\H_-}={\xi^{n_x}}_{|\H_-}$.
Puisque $\xi_{|\H_-}=\t_{|\H^-}$ est entrelac\'e
par $t$, on d\'eduit de \eqref{conjugx} que $t$ entrelace les caract\`eres
$\psi_{a_\b(x_n)}$ et $\psi_{a_\b(y_n)}$, \ie que~:
\begin{equation*}
t^{-1}a_\b(x_n)t \equiv a_\b(y_n) \mod{t^{-1}\HH_-^*t+\HH_-^*}.
\end{equation*}
Cette \'equivalence est certainement satisfaite dans tous les blocs sauf
peut-\^etre le bloc $\A^{12}$, o\`u la condition est~:
\begin{equation*}
a_\b(t^{-1}x_nt-y_n) \in
(t^{-1}(\HH^m(\b,\L))^*t+(\HH^m(\b,\L))^*)\cap\A^{12}.
\end{equation*}
D'apr\`es le lemme~\ref{ExactSequences}, il existe $x_n',y_n'$ dans
$\MM_{q-m+1}(\b,\L)\cap A^{12}$ tels que~:
\begin{equation*}
a_\b(t^{-1}x_n't-y_n')=a_\b(t^{-1}x_nt-y_n).
\end{equation*}
Donc $(t^{-1}x_n't-y_n')-(t^{-1}x_nt-y_n)$ appartient au $(1,2)$-bloc
de~:
\begin{equation*}
\left(t^{-1}\MM_{q-m}(\b,\L)t + \MM_{q-m}(\b,\L)\right) \cap\B
\subset t^{-1}\left(\PP_1(\L)\cap\B\right)t + \left(\PP_1(\L)\cap\B\right).
\end{equation*}
Il existe donc
$x_n'',y_n''\in\left(\PP_1(\L)\cap\B+\MM_{q-m+1}\right)\cap\A^{12}$
tels que $t^{-1}x_n''t-y_n''= t^{-1}x_nt-y_n$ et donc~:
\begin{equation*}
(1+x_n'')t(1+y_n'')^{-1}= n_xtn_y^{-1}.
\end{equation*}
Puisque $(1+x_n''),(1+y_n'')\in\Om\cap\N$ qui normalise $\xi$, et
puisque le groupe $\Om_{q-m}(\b,\L)$ normalise $\Om$, on peut 
supposer que $g=n_x^- m_x t (n_y^- m_y)^{-1}$. 
De la m\^e\-me mani\`ere, en regardant
la restriction de $\xi$ \`a $\H_+=\H^{m+1}(\b,\L)(\H^m(\b,\L)\cap\N)$,
on se ram\`ene au cas $g=m_x t m_y^{-1} \in\M$, ce qui d\'emontre le
th\'eor\`eme.
\end{proof}

On pose $\K=\H^m(\b,\L)(\Om\cap\N)$, qui est un sous-groupe ouvert
compact de $\G$, puisque $\Om$ normalise $\H^m(\b,\L)$. 
Comme $\Om$ normalise aussi $\xi$, on obtient~:

\begin{lemm}
\label{LemmeSansNom}
Il existe un unique caract\`ere $\tilde\xi$ de $\K$ qui est trivial
sur $\K\cap\N$ et qui prolonge $\xi$.
\end{lemm}

\begin{coro}
Si $g\in\N$ entrelace $\tilde\xi$, alors $g\in\K\cap\N$.
\end{coro}

\begin{proof}
Supposons que $g\in\N$ entrelace $\xi$. D'apr\`es le th\'eor\`eme
\ref{entrelacementscinde}, il existe $m\in\M$ et $\g_1,\g_2\in\Om$
tels que $g=\g_1^{-1} m \g_2$. Puisque $\Om$ poss\`ede un d\'ecomposition
d'Iwahori par rapport \`a $(\M,\P)$, on peut \'ecrire
$\g_i=\g_i^-\g_i^\M\g_i^+$, pour $i=1,2$, avec $\g_i^-\in\Om\cap\N^-$,
$\g_i^\M\in\Om\cap\M$ et $\g_i^+\in\Om\cap\N$. 
On a donc~:
\begin{equation*}
\g_1^+ g (\g_2^+)^{-1} = (\g_1^-\g_1^\M)^{-1} m \g_2^-\g_2^\M
\in \N\cap\P^-=\{1\},
\end{equation*}
d'o\`u on d\'eduit que $g=(\g_1^+)^{-1}\g_2^+\in\Om\cap\N=\K\cap\N$.
\end{proof}

\subsection{}

Pour $i\in\{1,2\}$, on suppose donn\'es d'une part un sous-groupe
ouvert $\tilde\K_i$ de $\U(\L)\G^i$ qui contient et normalise le
groupe $\H^m(\b,\L)\cap\G^i$, d'autre part une repr\'esentation
irr\'eductible $\varrho_i$ de $\tilde\K_i$ dont la restriction \`a
$\tilde\K_i\cap\K$ est multiple de
$\xi_{|\tilde\K_i\cap\H^m(\b,\L)}$.

\begin{coro}\label{cover}
\begin{enumerate}
\item[(i)] 
L'ensemble $\tilde\K=(\tilde\K_1\times\tilde\K_2)\cdot\K$ est un groupe.
\item[(ii)] 
Il existe une unique repr\'esentation irr\'eductible $\varrho$ de
$\tilde\K$ telle que les res\-trictions $\varrho_{|\tilde\K\cap\N}$ et
$\varrho_{|\tilde\K\cap\N^-}$ soient triviales, et que
$\varrho_{|\tilde\K\cap\M}\simeq\varrho_1\otimes\varrho_2$.
\item[(iii)] 
La paire $(\tilde\K,\varrho)$ est une paire couvrante de 
$(\tilde\K_1\times\tilde\K_2,\varrho_1\otimes\varrho_2)$.
\end{enumerate}
\end{coro}

\begin{proof} Avec l'\'el\'ement fortement $(\P,\tilde\K)$-positif~:
\begin{equation*}
\zeta=\begin{pmatrix} \w_\F&0\\0&1\end{pmatrix},
\end{equation*}
o\`u $\w_\F$ d\'esigne une uniformisante de $\F$, la d\'emonstration
est identique \`a celle de \cite[Corollary 6.6]{BK2}.
\end{proof}

\subsection{}
\label{S43}

On rappelle que $r=\lfloor q/2\rfloor+1$ et $s=\lceil q/2\rceil$.
On d\'efinit trois entiers $m_0=\max\{q-m,m\}$, $m_\ss=\max\{q-m+1,\ss\}$ et
$m_\rr=\max\{q-m,\rr-1\}$. Soit $\varepsilon=m_\ss-m_\rr$. 
On pose~:
\begin{equation*}
\K_l = \begin{cases}
\H^m(\b,\L)(\H^{l-\varepsilon+1}(\b,\L)\cap\N), &
m_\ss\< l< m_0+\epsilon,\\
\H^m(\b,\L)(\Om_{l+1}(\b,\L)\cap\N), &
q-m\< l < m_\ss, \\
\K_{q-m}(\U_{l+1}(\L)\cap\B\cap\N),&0\< l< q-m,
\end{cases}
\end{equation*}
et~:
\begin{equation*}
\Xi_l = \begin{cases}
((\U_{m-l+\varepsilon}(\L)\cap\B)\Om^{q-l+\varepsilon}(\b,\L))\cap\N^-, &
m_\ss\< l\< m_0+\varepsilon,\\
((\U_{m-l}(\L)\cap\B)\H^{q-l}(\b,\L))\cap\N^-, &
q-m<l<m_\ss, \\
(\U_{m-l}(\L)\cap\B)\cap\N^-, &0< l\< q-m.
\end{cases}
\end{equation*}
On a le lemme suivant.

\begin{lemm}\label{conjxi}
Pour $0< l< m+\varepsilon$, le groupe $\Xi_l$ agit transitivement sur les
caract\`eres de $\K_{l-1}$ qui prolongent $\tilde\xi_{|\K_l}$
\end{lemm}

\begin{proof}
On consid\'ere le cas $m_\ss\<l<m+\varepsilon$, les autres cas \'etant
similaires. 
Comme le quotient $\K_{l-1}/\K_l$ est ab\'elien, tout
caract\`ere de $\K_{l-1}$ qui \'etend $\tilde\xi_{|\K_l}$ est de la
forme $\tilde\xi_{|\K_{l-1}}\psi_b$, pour
$b\in\HH^{l-\varepsilon+1}(\b,\L)^*\cap\A^{21}$. Donc
$s(b)\in\PP_{-l+\varepsilon}(\L)\cap\B^{21}$.

D'apr\`es \cite[Lemma 2.3.8]{Br4}, il existe
$x\in\PP_{m-l+\varepsilon}(\L)\cap\B^{21}$ tel qu'on ait
$s(b)=s(c)x-xs(c)=s(a_c(x))$. 
Pour $h\in\K_{l-1}$, le commutateur $[1+x,h]$ appartient \`a 
$\H^{m}(\b,\L)$ et, comme $x\in\B$, la proposition
\ref{QuelleHorreurInsoutenableALInfiniRepetee} implique~:
\begin{equation*}
\tilde\xi([1+x,h])=\psi_c([1+x,h])=\psi_{a_c(x)}(h).
\end{equation*}

D'apr\`es le lemme~\ref{ExactSequences} il existe un \'el\'ement 
$y\in\MM_{q-l+\varepsilon}\cap\A^{21}$ tel que $a_\b(y)=b-a_c(x)$. 
D'apr\'es la proposition 
\ref{QuelleHorreurInsoutenableALInfiniRepetee}, pour $h\in\K_{l-1}$,
le commutateur $[1+y,h]$ appartient \`a $\H^{m+1}(\b,\L)$ et~:
\begin{equation*}
\tilde\xi([1+y,h]) =\t([1+y,h])=\psi_{(1+y)^{-1}a_\b(y)}(h).
\end{equation*}
On a $ya_\b(y)=0$ et $\psi_{a_c(x)}^{1+y}=\psi_{a_c(x)}$, car
$y\in\PP_1(\L)$.
On obtient donc~:
\begin{equation*}
\tilde\xi_{|\K_{l-1}}^{(1+x)(1+y)}
=\tilde\xi_{|\K_{l-1}}\psi_{a_\b(y)}\psi_{a_c(x)}
=\tilde\xi_{|\K_{l-1}}\psi_b,
\end{equation*}
d'o\`u le lemme.
\end{proof}

\begin{coro}\label{nonsc}
Soit $\pi$ une repr\'esentation lisse de $\G$ qui contient le 
ca\-rac\-t\`ere $\xi$ de $\H^m(\b,\L)$. 
Alors $\pi$ n'est pas supercuspidale.
\end{coro}

\begin{proof} 
En appliquant le lemme \ref{conjxi}, on voit que $\pi$
contient aussi le caract\`ere $\tilde\xi$ de $\K=\K_{0}$. Par le
corollaire \ref{cover}, la paire $(\K,\tilde\xi)$ est une paire
couvrante de $(\K\cap\M,\xi_{|\H^m(\b,\L)\cap\M})$. D'apr\`es
\cite[Theorem 7.9]{BK1}, la composante isotypique
$\pi_\N^{\xi_{|\H^m(\b,\L)\cap\M}}$ du module de Jacquet de $\pi$ par
rapport \`a $\P$ est non nulle.
\end{proof}

\subsection{}

Soit $\pi$ une repr\'esentation lisse de niveau non nul de $\G$ qui ne 
contient aucun caract\`ere simple d'un groupe $\H^1(\b,\L)$ avec $\L$
stricte.
D'apr\`es le th\'eor\`eme~\ref{??}, la repr\'esentation $\pi$ contient
alors un caract\`ere scind\'e ou une strate scind\'ee.
Le cas de la strate scind\'ee est d\'ej\`a r\'egl\'e par 
\cite[Theorem 1.2.3]{Br4}.
On suppose donc qu'il existe une strate simple $[\L,n,0,\b]$ avec $\L$
stricte, un caract\`ere simple $\t\in\Cc(\L,m,\b)$, et un
$c\in\PP_{-m}(\L)$, tels que $\pi$ contienne le caract\`ere 
$\vartheta=\t\psi_c$ et que la strate $[\Ga,m,m-1,s(c)]$ soit
scind\'ee, o\`u, avec les notations habituelles, $\Ga$ est une
$\Oo_{\D_\E}$-suite de r\'eseaux de $\V_\E$ qui correspond \`a $\L$
par le th\'eor\`eme~\ref{CNNBr1}.

\begin{prop}
\label{scindecomp}
Il existe une d\'ecomposition $\V_\E=\V_\E^1\oplus \V_\E^2$
conforme \`a $\Ga$ telle que :
\begin{enumerate}
\item[(i)] 
la strate $[\Ga,m,m-1,s(c)]$ est scind\'ee par cette d\'ecomposition ;
\item[(ii)] 
l'\'el\'ement $s(c)^1\in\B^1$ normalise $\Ga^1$ et 
$\v_{\Ga^1}(s(c)^1)=-m$. 
\end{enumerate}
\end{prop}

\begin{proof} 
Bien s\^ur, $\B^1$ et $\Ga^1$ d\'esignent respectivement 
la $\E$-alg\`ebre $\End_{\D_\E}(\V_\E^1)$ et la suite 
$\Ga\cap\V_\E^1$.
Dans le cas o\`u $\Ga$ est stricte, c'est 
\cite[Proposition 2.2.1]{Br1}.
La d\'emonstration dans le cas g\'en\'eral est identique. 
\end{proof}

Soit $\V=\V^1\oplus\V^2$ la d\'ecomposition qui correspond \`a celle
de $\V_\E$ donn\'ee par la proposition~\ref{scindecomp}. 
Comme $s(c)$ stabilise la d\'ecomposition $\V_\E=\V_\E^1\oplus\V_\E^2$,
on a $s(\e^i c\e^j)=0$, pour $i\ne j$, o\`u $\e^i$ est le projecteur
sur $\V^i$.
D'apr\`es le lem\-me \ref{ExactSequences}, il existe
$x\in\MM_{q-m}$ tel que $\e^ia_\b(x)\e^j=-\e^i c\e^j$, pour $i\ne j$.
D'apr\'es la
proposition~\ref{QuelleHorreurInsoutenableALInfiniRepetee}, 
l'\'el\'ement $1+x$ normalise $\H^{m}(\b,\L)$ et~:
\begin{equation*}
\vartheta^{1+x} = \t\psi_{(1+x)^{-1}\b(1+x)-\b} \psi_c
= \t\psi_{c'},
\end{equation*}
o\`u $c'=c+a_\b(x)$. 
Rempla\c cant $c$ par $c'$, on voit qu'on est dans la situation des
\S\S\ref{S41}--\ref{S43}.
Par le corollaire~\ref{nonsc}, on conclut que $\pi$ n'est pas
supercuspidale, ce qui termine la d\'emonstration du
th\'eor\`eme~\ref{scsimple}.


\section{Le niveau z\'ero}
\label{YoupiNiveauZero}

Soit $\A$ une $\F$-alg\`ebre centrale simple et soit $\G$ son groupe
multiplicatif.
D\'esormais, toute les strates sont relatives \`a une suite de
r\'eseaux stricte.
On peut donc remplacer le langage des suites de r\'eseaux par celui
des ordres h\'er\'editaires.
Par commodit\'e, on fixe tout de m\^eme, comme d'habitude, un
$\A$-module \`a gauche simple $\V$, et on note $\D$ l'alg\`ebre
oppos\'ee \`a $\End_{\A}(\V)$.

Dans cette section, on prouve que toute repr\'esentation
irr\'eductible supercuspidale de niveau non nul de $\G$ contient un
type simple maximal au sens de \cite{VS3}.

\subsection{}

Soit $[\AA,n,0,\b]$ une strate simple de $\A$.
Soit~:
\begin{equation}
\label{DecEDdeRef}
\V=\V^1\oplus\ldots\oplus\V^l
\end{equation}
une d\'ecomposition de $\V$
en sous-$\E\otimes_{\F}\D$-modules, qui soit conforme \`a $\AA$.
Soit $\M$ le sous-groupe de Levi de $\G$ correspondant et soit $\P$ 
un sous-groupe para\-bo\-lique de $\G$ de facteur de Levi $\M$. 
On \'ecrit $\P=\M\N$, o\`u $\N$ est le radical unipotent de $\P$.
On note $\N^-$ le radical unipotent du sous-groupe para\-bo\-lique
oppos\'e \`a $\P$.
On pose~:
\begin{eqnarray*}
\H^1_{\P}&=&\H^1(\b,\AA)\left(\J^1(\b,\AA)\cap\N\right),\\
\J^1_{\P}&=&\H^1(\b,\AA)\left(\J^1(\b,\AA)\cap\P\right).
\end{eqnarray*}
Ce sont des sous-groupes ouverts compacts de $\J^1(\b,\AA)$ contenant
$\H^1(\b,\AA)$.
Pour simplifier les notations, on note $\H^1$, $\J^1$ 
respectivement pour les groupes $\H^1(\b,\AA)$, $\J^1(\b,\AA)$. 
Soit $\t\in\Cc(\AA,0,\b)$ un caract\`ere simple.
On note $\t_{\P}$ le caract\`ere de $\H^1_{\P}$ d\'efini par
$\t_\P(hu)=\t(h)$, pour $h\in\H^1$ et $u\in\J^1\cap\N$.

\subsection{}

Soit $\B$ le commutant de $\E$ dans $\A$ et soit $\V_{\E}$ un
$\B$-module \`a gauche simple.
On note $\D_{\E}$ l'alg\`ebre oppos\'ee \`a $\End_{\B}(\V_{\E})$
et $m_\E$ la dimension de $\V_\E$ sur $\D_{\E}$.
Soient $\e^i$ les idempotents de $\BB=\AA\cap\B$ d\'efinis par la
d\'ecomposition (\ref{DecEDdeRef}) et soit $n_i$ la dimension de
$\e^i\V_\E$.

\begin{defi}
\label{DefSubordonnee!}
La d\'ecomposition (\ref{DecEDdeRef}) est dite
{\it subordonn\'ee \`a} $\BB$ 
s'il exis\-te un isomorphisme de $\E$-alg\`ebres
$\Psi:\B\f\Mat_{m_\E}(\D_\E)$ tel que~:
\begin{enumerate}
\item[(i)]
Pour chaque $1\<i\<l$, l'idempotent $\Psi(\e^i)$ est \'egal \`a~: 
\begin{equation*}
\I^i={\rm diag}(0,\ldots,{\rm Id}_{n_i},\ldots,0),
\end{equation*}
o\`u la matrice identit\'e ${\rm Id}_{n_i}\in\Mat_{n_i}(\D_\E)$
appara\^\i t \`a la $i$-i\`eme place.
\item[(ii)]
L'ordre h\'er\'editaire $\Psi(\BB)$ est la sous-$\o_\E$-alg\`ebre de
$\Mat_{m_E}(\o_{\D_\E})$ cons\-ti\-tu\'ee des matrices dont la r\'eduction
modulo $\p_{\D_\E}$ est triangulaire sup\'erieure par blocs de
taille $(n_1,\ldots,n_l)$.
\end{enumerate}
\end{defi}

\begin{rema}
Si (\ref{DecEDdeRef}) est subordonn\'ee \`a $\BB$, alors $l$ est
\'egal \`a la p\'eriode de $\BB$.
Les $\I^i$ d\'e\-fi\-nis\-sent une d\'e\-com\-po\-si\-tion de
$\D_\E^{m_\E}$ conforme \`a $\Psi(\BB)$, et l'ordre
$\I^i\Psi(\BB)\I^i$ est un ordre maximal de $\Mat_{n_i}(\D_\E)$ 
\'egal \`a $\Mat_{n_i}(\o_{\D_\E})$.
\end{rema}

\subsection{}

On suppose d\'esormais que la d\'ecomposition (\ref{DecEDdeRef}) est 
subordonn\'ee \`a $\BB$.
On note $e$ la p\'eriode de $\BB$.

\begin{prop}
\label{Bitoniau}
\begin{enumerate}
\item[(i)]
On a des d\'ecompositions d'Iwahori~:
\begin{eqnarray*}
\J^1_{\P}&=&(\H^1\cap\N^-)\cdot(\J^1\cap\M)\cdot(\J^1\cap\N),\\
\H^1_{\P}&=&(\H^1\cap\N^-)\cdot(\H^1\cap\M)\cdot(\J^1\cap\N).
\end{eqnarray*}
\item[(ii)]
On a des isomorphismes de groupes~:
\begin{equation*}
\J^1_{\P}/\H^1_{\P}
\simeq\J^1\cap\M/\H^1\cap\M
\simeq\prod\limits_{i=1}^e\J^1(\b,\AA^i)/\H^1(\b,\AA^i).
\end{equation*}
\item[(iii)]
L'appli\-ca\-tion $(x,y)\mapsto\t_\P([x,y])$ d\'efinit un 
espace symplectique non d\'eg\'en\'er\'e 
$(\J^1_{\P}/\H^1_{\P},\boldsymbol{k}_{\t_\P})$ isomorphe \`a la somme
directe des espaces symplectiques
$(\J^1(\b,\AA^i)/\H^1(\b,\AA^i),\boldsymbol{k}_{\t^i})$, o\`u $\t^i$
est le transfert de $\t$ \`a $\Cc(\AA^i,0,\b)$.
\end{enumerate}
\end{prop}

\begin{proof}
Voir \cite[Proposition 7.2.3]{BK}.
Cela d\'ecoule des d\'e\-com\-po\-si\-tions d'Iwahori existant pour
$\J^1$ et $\H^1$, du th\'eor\`eme \ref{PaireDecomposeeTheta} et enfin
de la proposition \ref{bfsndcs}.
\end{proof}

\begin{prop}
\label{Ent1}
On a $\I_\G(\t_\P)=\J^1_\P\mult\B\J^1_\P$.
\end{prop}

\begin{proof}
Il suffit de v\'erifier que $\I_\G(\t_\P)$ contient $\mult\B$.
Soit~: 
\begin{equation}
\label{DecFineSub}
\V=\W^1\oplus\ldots\oplus\W^{l_0}
\end{equation}
une d\'ecomposition de $\V$ en $\E\otimes_\F\D$-modules {\it simples}, 
qui soit plus fine que (\ref{DecEDdeRef}) et conforme \`a $\AA$.
Soit $\M_0\subset\M$ le stabilisateur de (\ref{DecFineSub}) et soit
$\P_0=\M_0\N_0$ un sous-groupe parabolique de sous-groupe de Levi
$\M_0$ et contenu dans $\P$.
Alors~: 
\begin{equation*}
\U=(\U(\BB)\cap\P_0)\U_1(\BB)
\end{equation*}
est un sous-groupe d'Iwahori de $\mult\B$.
Puisque $\P_0$ est inclus $\P$, le groupe $\U$ est inclus dans 
$(\U(\BB)\cap\P)\J^1_\P$, qui normalise $\t_\P$.
D'apr\`es la d\'ecomposition de Bruhat de $\mult\B$ en doubles classes
modulo $\U$, il suffit donc de montrer que tout \'el\'ement du
normalisateur de $\M_0$ dans $\mult\B$ entrelace $\t_\P$.
Soit donc $y$ dans ce normalisateur.
Le groupe $\H^1_\P$ a des d\'ecompositions d'Iwahori relativement \`a
$(\M_0,\P_0)$ et \`a $(\M_0,\P_0^y)$.
Le groupe $\H^1_\P\cap y\H^1_\P y^{-1}$ admet donc lui-m\^eme une
d\'ecomposition d'Iwahori relativement \`a $(\M_0,\P_0)$. 
Puisque $\t_\P$ est trivial sur $\H^1_\P\cap\N_0$ et sur
$\H^1_\P\cap\N_0^y$, chacun des deux caract\`eres $\t_\P$ et
${}^y\t_\P$ est trivial sur $\H^1_\P\cap y\H^1_\P y^{-1}\cap\N_0$.
On a un r\'esultat analogue pour le radical unipotent oppos\'e
$\N_0^-$. 
Il reste donc \`a v\'erifier que $y$ entrelace~:
\begin{equation*}
\t_{\P|\H^1_\P\cap\M_0}=\t_{|\H^1\cap\M_0},
\end{equation*}
ce qui est le cas puisque $y$ entrelace $\t$.
\end{proof}

\subsection{}

On rappelle ({\it cf.} \cite[\S2.2]{VS2}) qu'il existe une
repr\'esentation irr\'eductible $\n$ de $\J^1$, uni\-que \`a
isomorphisme pr\`es, dont la restriction \`a $\H^1$ contient 
$\t$.
Elle est normalis\'ee par $(\KK(\AA)\cap\mult\B)\J$, son entrelacement  
vaut $\I_\G(\n)=\J^1\mult\B\J^1$ et, pour tout $y\in\mult\B$, on a  
$\dim\Hom_{\J^1\cap(\J^{1})^y}(\n,\n^y)=1$.

\begin{prop}
\label{Ent2}
Il existe une repr\'esentation irr\'eductible $\n_\P$ de $\J^1_\P$, 
uni\-que \`a isomorphisme pr\`es, dont la restriction \`a $\H^1_\P$
contient $\t_\P$.
En outre~:
\begin{enumerate}
\item[(i)]
Les repr\'esentations $\Ind_{\J^1_\P}^{\J^1}(\n_\P)$ et $\n$ sont
isomorphes. 
\item[(ii)]
Pour tout $y\in\mult\B$, il existe une unique
$(\J^1_\P,\J^1_\P)$-double classe dans $\J^1y\J^1$ entrela\c cant
$\n_\P$ et $\dim\Hom_{\J^1_\P\cap{\J^{1y}_\P}}(\n_\P,\n_\P^y)=1$.
\end{enumerate}
\end{prop}

\begin{rema}
En particulier, on a $\I_\G(\n_\P)=\J^1_\P\mult\B\J^1_\P$.
\end{rema}

\begin{proof}
Pour le point (i), l'argument est identique \`a celui utilis\'e pour 
\cite[Proposition 7.2.4]{BK}.
Pour le (ii), voir \cite[Corollary 4.1.5]{BK}.
\end{proof}

\begin{rema}
Les propositions \ref{Bitoniau} \`a \ref{Ent2} sont valables pour
une d\'e\-com\-po\-si\-tion (\ref{DecEDdeRef}) quelconque, \ie
conforme mais pas n\'ecessairement su\-bor\-donn\'ee. 
Il suffit, dans la propositions \ref{Bitoniau}, de remplacer $\AA^i$
par $\L^i$, o\`u $\L$ est une $\o_\D$-cha\^\i ne d\'efinissant $\AA$.
\end{rema}

\subsection{}

On suppose que la d\'ecomposition (\ref{DecEDdeRef}) est 
subordonn\'ee \`a $\BB$. 
On rappelle ({\it cf.} \cite[\S2.4]{VS2}) qu'une {\it $\b$-extension}
de $\n$ est une repr\'esentation de $\J=\J(\b,\AA)$ prolongeant $\n$
dont l'entrelacement contient $\mult\B$.
On pose~:
\begin{equation*}
\J_{\P}=\H^1(\b,\AA)\left(\J(\b,\AA)\cap\P\right).
\end{equation*}
On fixe une $\b$-extension $\k$ de $\n$ et on note $\k_\P$ la
repr\'esentation de $\J_\P$ sur les $(\J\cap\N)$-invariants de $\k$.

\begin{prop}
\begin{enumerate}
\item[(i)]
On a $\I_\G(\k_\P)=\J_\P\mult\B\J_\P$.
\item[(ii)]
Les repr\'esentations $\Ind_{\J_\P}^{\J}(\k_\P)$ et $\k$ sont
isomorphes.
\end{enumerate}
\end{prop}

\begin{proof}
Le point (ii) d\'ecoule directement du fait que la restriction de
$\k_\P$ \`a $\J^1_\P$ est \'egale \`a $\n_\P$. 
Traitons le point (i).
D'apr\`es (ii), pour chaque \'e\-l\'e\-ment $y\in\mult\B$, il existe une
unique double 
classe $\J_\P x\J_\P$ dans $\J y\J$ entrela\c cant $\k_\P$.
Par d\'ecomposition d'Iwahori, on peut supposer que $x$ appartient
\`a~: 
\begin{equation*}
(\J\cap\N^-)y(\J\cap\N^-)=(\J^1\cap\N^-)y(\J^1\cap\N^-),
\end{equation*}
donc \`a $\J^1y\J^1$.
Puisque $\k_\P$ prolonge $\n_\P$, l'\'el\'ement $x$ entrelace
$\n_\P$. 
D'apr\`es les propositions \ref{Ent1} et \ref{Ent2}, la double classe
$\J^1_\P y\J^1_\P$ est la seule double classe dans $\J^1y\J^1$ qui
entrelace $\n_\P$.
Ainsi $x\in\J^1_\P y\J^1_\P$ et $\J_\P x\J_\P=\J_\P y\J_\P$,
de sorte que $y$ entrelace $\k_\P$.
\end{proof}

\begin{rema}
Si la d\'ecomposition (\ref{DecEDdeRef}) n'est pas subordonn\'ee, 
les grou\-pes $\U(\BB)$ et $\J(\b,\AA)$ n'admettent pas, en
g\'en\'eral, de d\'e\-com\-po\-si\-tion d'Iwa\-ho\-ri relativement \`a
$(\M,\P)$. 
\end{rema}

\begin{prop}
\label{entrelacementvthP}
Soit $\xi$ une repr\'esentation irr\'eductible de $\J_\P$ triviale sur
$\J^1_\P$. 
On a $\I_\G(\k_\P\otimes\xi)=\J_\P\I_{\mult\B}(\xi)\J_\P$.
\end{prop}

\begin{proof}
La preuve est analogue \`a celle de \cite[Proposition 5.3.2]{BK}.
Il suffit de remplacer \cite[Proposition 5.1.8]{BK} par la proposition
\ref{Ent2}.
\end{proof}

\subsection{}
\label{Slieover}

Soit $\k$ une $\b$-extension de $\n$ et soit $\s$ l'inflation \`a
$\J$ d'une repr\'esentation irr\'eductible du groupe $\J/\J^1$.
On pose $\vartheta=\k\otimes\s$.
Le quotient $\J/\J^1$ est isomorphe \`a $\U(\BB)/\U_1(\BB)$, qui est le
groupe des points rationnels d'un groupe r\'eductif sur le corps fini
$k_{\D_\E}$. 

Soit $\AA'$ un ordre h\'er\'editaire $\E$-pur de $\A$ tel que
l'intersection de $\BB$ avec l'ordre $\BB'=\AA'\cap\B$ soit un 
ordre h\'er\'editaire.
On note $\t'$ le transfert de $\t$ \`a $\Cc(\AA',0,\b)$ et $\n'$
l'unique repr\'esentation irr\'eductible de $\J^1(\b,\AA')$ 
contenant $\t'$. 
On rappelle comment, dans \cite{VS3}, on associe \`a $\AA'$ une
$\b$-extension $\k'$ de $\n'$.
En proc\'edant comme dans \cite[Proposition 4.5]{VS3}, on construit
une famille finie~:
\begin{equation}
\label{appafam}
(\AA_0,\ldots,\AA_k), \quad k\>0,
\end{equation}
d'ordres $\E$-purs de $\A$, avec $\AA_0=\AA$ et $\AA_k=\AA'$, et telle
que pour tout $0\<i<k$, l'ordre $\AA_{i}$ ou bien contienne ou bien
soit contenu dans $\AA_{i+1}$.
(Il suffit de tracer, dans l'immeuble de Bruhat-Tits de $\G$, le
segment joignant $\AA$ et $\AA'$.)
Soit $\t_{i}$ le transfert de $\t$ \`a $\Cc(\AA_{i},0,\b)$ et soit
$\n_{i}$ l'unique repr\'esentation irr\'eductible de
$\J^1(\b,\AA_{i})$ contenant $\t_{i}$.
On d\'efinit par r\'ecurrence une famille finie~:
\begin{equation}
\label{kappafam}
(\k_0,\ldots,\k_k), \quad k\>0,
\end{equation}
de $\b$-extensions, en posant $\k_0=\k$ et, pour $0\<i<k$, en prenant
pour $\k_{i+1}$ l'unique $\b$-extension de $\n_{i+1}$ qui soit
coh\'erente avec $\k_i$ au sens de \cite[\S2.4.4]{VS3}.
On pose enfin $\k'=\k_k$.

\begin{defi} 
Une repr\'esentation $\vartheta'$ de $\J(\b,\AA')$ est dite
\emph{\audessus} $\vartheta$ si elle est de la forme
$\vartheta'=\k'\otimes\s'$, o\`u $\s'$ est l'inflation \`a
$\J(\b,\AA')$ d'une re\-pr\'e\-sen\-ta\-tion ir\-r\'e\-duc\-ti\-ble 
de $\J(\b,\AA')/\J^1(\b,\AA')$ telle que $\s$ et $\s'$
s'entrelacent sur $\U(\BB\cap\BB')$.
\end{defi}

\begin{rema}
\label{Oignon}
\begin{enumerate}
\item[(i)]
Par exemple, si $\BB'$ contient $\BB$, on peut consid\'erer
$\s$ comme une repr\'esentation du sous-groupe parabolique
$\U(\BB)/\U_1(\BB')$ de $\U(\BB')/\U_1(\BB')$.
La condition sur $\s'$ signifie alors que $\s'$ est l'inflation \`a
$\J(\b,\AA')$ d'une composante ir\-r\'e\-duc\-ti\-ble de l'induite de
$\s$ \`a $\U(\BB')/\U_1(\BB')$.
\item[(ii)]
Si $\BB'$ est contenu dans $\BB$, la condition sur $\s'$ signifie que
$\s'$ est l'inflation \`a $\J(\b,\AA')$ d'une composante
ir\-r\'e\-duc\-ti\-ble de la restriction de $\s$ \`a
$\U(\BB')/\U_1(\BB)$.
\end{enumerate}
\end{rema}

\begin{prop}
\label{lieover}
Soit $\pi$ une repr\'esentation irr\'eductible de $\G$ contenant
$\vartheta$. 
Alors $\pi$ contient une repr\'esentation $\vartheta'$ de
$\J(\b,\AA')$ \audessus\ $\vartheta$.
\end{prop}

\begin{proof}
On proc\`ede en trois \'etapes.
\begin{enumerate}
\item[(i)]
On consid\`ere d'abord le cas o\`u $\AA'$ ou bien contient, ou bien 
est contenu dans $\AA$. 
La d\'e\-mons\-tra\-tion est analogue \`a celle 
de~\cite[Proposition 8.3.5]{BK}~: voir \emph{op.\ cit.} p.296.
Il suffit de remplacer~\cite[(5.2.14)]{BK} par 
\cite[Proposition 2.29]{VS2} et \cite[Proposition 5.3.2]{BK} par
\cite[Lemme 4.2]{VS3}. 
Dans le cas o\`u $\AA'$ est contenu dans $\AA$, on a m\^eme un
r\'esultat plus pr\'ecis~: on voit que toute repr\'esentation de
$\J(\b,\AA')$ \audessus\ $\vartheta$ est contenue dans $\pi$.
\item[(ii)]
On consid\`ere ensuite le cas o\`u $\BB'$ ou bien contient, ou bien 
est contenu dans $\BB$. 
Dans ce cas, la famille (\ref{appafam}) peut \^etre
choisie de telle sorte que $\AA_i\cap\B=\BB$ pour $0\<i\<k-1$.
On peut donc d\'efinir $\vartheta_i=\k_i\otimes\s$.
Alors $\vartheta_i$ est \emph{l'unique} repr\'esentation de
$\J(\b,\AA_i)$ \audessus\ $\vartheta$.
C'est m\^eme l'unique re\-pr\'e\-sen\-ta\-tion de $\J(\b,\AA_i)$ 
\audessus\ $\vartheta_j$ pour tout $0\<j\<k-1$.
En appliquant (i) successivement \`a chaque paire
$\{\AA_i,\AA_{i+1}\}$ au lieu de $\{\AA,\AA'\}$, on voit que $\pi$
contient $\vartheta$ si et seulement si elle contient $\vartheta_{k-1}$. 
En appliquant encore (i) avec $\AA_{k-1}$ au lieu de $\AA$, on voit
que $\pi$ contient $\vartheta_{k-1}$ si et seulement si elle contient
une repr\'esentation $\vartheta'$ \audessus\ $\vartheta_{k-1}$, ce qui
est la m\^eme chose que d'\^etre \audessus\ $\vartheta$.
\item[(iii)]
On consid\`ere enfin le cas g\'en\'eral.
Pour se ramener \`a (ii), on passe d'abord de $\BB$ \`a $\BB\cap\BB'$,
de sorte que toute repr\'esentation de
$\J(\b,\AA\cap\AA')$ \audessus\ $\vartheta$ est contenue dans $\pi$, 
puis on passe de $\BB\cap\BB'$ \`a $\BB'$, de sorte que $\pi$ contient
une repr\'esentation $\vartheta'$ de $\J(\b,\AA')$ \audessus\
$\vartheta$.
\end{enumerate}
Ceci termine la d\'emonstration de la proposition \ref{lieover}.
\end{proof}

\begin{rema}
\label{Humito}
En particulier, si $\pi$ contient $\vartheta$, alors $\pi$ contient
aussi $\vartheta'=\k'\otimes\s$ pour tout ordre h\'er\'editaire
$\E$-pur $\AA'$ tel que $\AA'\cap\B=\BB$.
\end{rema}

\subsection{}
\label{Scusp}

Soit $\pi$ une repr\'esentation irr\'eductible de $\G$ contenant un
caract\`ere simple d'un groupe $\H^1$ -- \ie qu'il existe un couple 
$([\AA,n,0,\b],\t)$ constitu\'e d'une strate simple et d'un
caract\`ere simple $\t\in\Cc(\AA,0,\b)$ tels que la restriction de 
$\pi$ \`a $\H^1(\b,\AA)$ con\-tien\-ne $\t$. 
Parmi ces couples, on en choisit un tel que l'ordre h\'er\'editaire
$\AA$ soit minimal.  

Soit $\n$ l'unique repr\'esentation irr\'eductible de $\J^1(\b,\AA)$
contenant $\t$.
La res\-tric\-tion de $\pi$ \`a $\J^1(\b,\AA)$ contient donc $\eta$,
et la restriction de $\pi$ \`a $\J(\b,\AA)$ contient une
repr\'esentation de la forme $\vartheta=\k\otimes\s$, o\`u $\k$ est
une $\b$-extension de $\n$ et $\s$ est l'inflation \`a $\J(\b,\AA)$
d'une repr\'esentation irr\'eductible de 
$\J(\b,\AA)/\J^1(\b,\AA)\simeq\U(\BB)/\U_1(\BB)$, qu'on note encore
$\s$.

\begin{prop}\label{sigmacusp}
Dans cette situation, $\s$ est une repr\'esentation cuspidale de 
$\U(\BB)/\U_1(\BB)$.
\end{prop}

\begin{proof} 
Supposons que $\s$ n'est pas une repr\'esentation cuspidale de
$\Gg=\U(\BB)/\U_1(\BB)$. 
Il existe donc un sous-groupe parabolique propre $\Pp$ de $\Gg$, de
radical unipotent $\Uu$, tel que la restriction de $\s$ \`a $\Uu$
contienne le caract\`ere trivial. 
Il existe donc une repr\'esentation $\s'$ de $\Pp/\Uu$ telle que
$\s$ soit une composante irr\'eductible de
$\Ind_{\Pp}^{\Gg}(\s')$. 
Il existe un unique ordre h\'er\'editaire $\BB'$ contenu dans $\BB$
tel que $\Pp$ soit l'image de $\U(\BB')$ par l'application 
quotient $\U(\BB)\to\Gg$.
Le radical unipotent $\Uu$ est alors l'image de $\U_1(\BB')$.  
Par~\cite[Lemme 1.7]{VS2}, il existe un ordre h\'er\'editaire $\E$-pur
$\AA'$ contenu dans $\AA$ tel que $\AA'\cap\B=\BB'$. 
Plus pr\'ecis\'ement, $\AA'$ est strictement contenu dans $\AA$,
puisque $\BB'$ l'est dans $\BB$. 

Soit $\t'$ le transfert de $\t$ \`a $\Cc(\AA',0,\b)$, soit $\n'$
l'unique repr\'esentation ir\-r\'e\-duc\-ti\-ble de $\J^1(\b,\AA')$ 
contenant $\t$ et soit $\k'$ la $\b$-extension de $\n'$ construite 
comme en (\ref{kappafam}). 
L'inflation \`a $\J(\b,\AA')$ de la repr\'esentation $\s'$ de~:
\begin{equation*}
\Pp/\Uu\simeq\U(\BB')/\U_1(\BB')\simeq\J(\b,\AA')/\J^1(\b,\AA')
\end{equation*}
est encote not\'ee $\s'$, et on pose $\vartheta'=\k'\otimes\s'$.
Nous sommes alors dans la situation du \S\ref{Slieover} 
et $\vartheta$ est \audessus\ $\vartheta'$. 
Par la proposition~\ref{lieover}, la repr\'esentation $\pi$ contient
donc aussi $\vartheta'$ et \emph{a fortiori} $\t'$. 
Puisque $\AA'$ est strictement inclus dans $\AA$, ceci contredit la
minimalit\'e de $\AA$.
\end{proof}

\subsection{}
\label{SvthP}

On continue avec les notations du paragraphe pr\'ec\'edent |
donc $\pi$ contient une repr\'esentation de la forme
$\vartheta=\k\otimes\s$ avec $\s$ cuspidale. 
Par la proposition~\ref{lieover} (voir aussi la remarque \ref{Humito}), 
on peut changer l'ordre h\'eriditaire $\AA$ sans changer sa trace sur
$\B$.
On peut donc supposer que~: 
\begin{equation}
\label{SoundGrabitz}
\KK(\AA)\cap\mult\B=\KK(\BB).
\end{equation}
Il suffit de choisir l'ordre h\'eriditaire $\E$-pur $\AA$ associ\'e
\`a $\BB$ par le th\'eor\`eme \ref{Edescent}.
(Dans la terminologie de Grabitz~\cite{Grabitz}, un ordre principal
$\E$-pur $\AA$ v\'erifiant (\ref{SoundGrabitz}) est dit \emph{sound}.)

Soit $\V=\V^1\oplus\ldots\oplus\V^e$ une d\'ecomposition de $\V$
subordonn\'ee \`a $\BB$. 
Soit $\M$ le sous-groupe de Levi de $\G$ qui est le stabilisateur de
cette d\'ecomposition, et soit $\P=\M\N$ un sous-groupe parabolique de
$\G$ de facteur de Levi $\M$. 
On a un isomorphisme de groupes~: 
\begin{equation*}
\J_\P\cap\M\simeq\prod_{i=1}^e\J(\b,\AA_i).
\end{equation*}
Si on note $\t_i$ la restriction de $\t$ \`a $\H^1(\b,\AA_i)$, 
\ie le transfert de $\t$ \`a $\Cc(\AA_i,0,\b)$, et $\n_i$ l'unique
repr\'esentation irr\'eductible de $\J^1(\b,\AA_i)$ contenant $\t_i$,
il existe pour chaque $i$ une $\b$-extension $\k_i$ de $\n_i$ telle
que la restriction de $\k_\P$ \`a $\J_\P\cap\M$ soit \'equivalente
\`a~: 
\begin{equation*}
\k_1\otimes\ldots\otimes\k_e.
\end{equation*}
(En effet, l'entrelacement de $\k_{\P|\J_\P\cap\M}$ contient
$\mult\B\cap\M$.)
De fa\c con analogue, il existe, pour chaque entier $i$, une
repr\'esentation $\s_i$ de $\J(\b,\AA_i)$ qui est l'inflation d'une
repr\'esentation irr\'eductible cuspidale de
$\J(\b,\AA_i)/\J^1(\b,\AA_i)\simeq\U(\BB_i)/\U_1(\BB_i)$ telle que
la restriction de $\s$ \`a $\J\cap\M$ soit \'equivalente \`a~: 
\begin{equation*}
\s_1\otimes\ldots\otimes\s_e.
\end{equation*}
Si on consid\`ere $\s$ comme repr\'esentation de 
$\J_\P/\J^1_\P\simeq\J/\J^1$, on pose
$\vartheta_\P=\k_\P\otimes\s$, qui est une repr\'esentation de $\J_\P$. 
Comme dans~\cite[Proposition 7.2.17]{BK}, on a le r\'esultat suivant. 

\begin{prop}\label{vthP}
\begin{enumerate}
\item[(i)] 
$\vartheta_\P$ est irr\'eductible et 
$\vartheta\simeq\Ind_{\J_\P}^{\J}(\vartheta_\P)$.
\item[(ii)] 
Les restrictions de $\vartheta_\P$ \`a $\J_\P\cap\N$ et
$\J_\P\cap\N^-$ sont triviales, et la restriction de $\vartheta_\P$
\`a $\J_\P\cap\M$ est \'equivalente \`a 
$\vartheta_1\otimes\ldots\otimes\vartheta_e$,
o\`u $\vartheta_i=\k_i\otimes\s_i$ est irr\'eductible.
\end{enumerate}
\end{prop}

\subsection{}

On fixe une extension non ramifi\'ee $\LL/\E$ maximale dans $\D_\E$
et une uniformisante $\unif$ de $\D_\E$ normalisant $\LL$. 
Le groupe de Galois de $\LL/\E$ est engendr\'e par 
${\rm Ad}(\unif)$, l'automorphisme de conjugaison par
$\unif$.
Par r\'eduction, on identifie les groupes de Galois $\Gal(\LL/\E)$ 
et $\Gal(k_{\D_\E}/k_\E)$ et on note $\phi$ l'image de 
${\rm Ad}(\unif)$ dans $\Gal(k_{\D_\E}/k_\E)$.

On fixe un isomorphisme de $\E$-alg\`ebres
$\B\simeq\Mat_{m_\E}(\D_\E)$ v\'erifiant les conditions de la
d\'efinition \ref{DefSubordonnee!}, dont on reprend les notations. 
Cet isomorphisme induit des isomorphismes de groupes
$\mult\B\simeq\GL_{m_\E}(\D_\E)$ et~: 
\begin{equation*}
\U(\BB)/\U^1(\BB)\simeq
\GL_{n_1}(k_{\D_\E})\times\ldots\times\GL_{n_e}(k_{\D_\E}).
\end{equation*}
De cette fa\c con, le groupe $\Gal(k_{\D_\E}/k_\E)$ op\`ere sur les
repr\'esentations de $\J/\J^1\simeq\U(\BB)/\U^1(\BB)$, et notamment 
sur les $\s_i$.

\subsection{}

Dans ce paragraphe, on suppose que les $\s_i$ ne sont pas tous dans
une seule orbite sous $\Gal(k_{\D_\E}/k_\E)$.
Plus pr\'ecis\'ement, on note $\I$ l'ensemble des entiers $1\<i\<e$
tels que $\s_i$ soit \'equivalent \`a un conjugu\'e de $\s_1$ sous
$\Gal(k_{\D_\E}/k_\E)$, et on suppose que $\I$ n'est pas \'egal \`a 
$\{1,\ldots,e\}$ tout entier.
On pose~:
\begin{equation*}
\W=\bigoplus_{i\in\I}\V^i,\quad\W'=\bigoplus_{i\notin\I}\V^i.
\end{equation*}
Soit $\M'$ le stabilisateur de la d\'ecomposition $\V=\W\oplus\W'$, 
qui est un sous-groupe de Levi de $\G$, et soit $\P'=\M'\N'$ un
sous-groupe parabolique de $\G$ de facteur de Levi $\M'$. 
Soit $\P=\M\N$ un sous-groupe parabolique de $\G$ de facteur de Levi
$\M$ tel que $\P\subset\P'$.

\begin{prop} 
Dans cette situation, $(\J_\P,\vartheta_\P)$ est une
paire couvrante de $(\J_\P\cap\M',\vartheta_\P{}_{|\J_\P\cap\M'})$.
\end{prop}

\begin{proof} 
Nous allons d'abord majorer l'entrelacement de $\vartheta_\P$.
Soit $\U$ le sous-groupe d'Iwahori de $\B^\times$
contenu dans $\U(\BB)$ s'identifiant au sous-groupe d'Iwahori
standard de $\GL_{m_\E}(\D_\E)$.
Soit $\tilde\W$ le groupe de Weyl affine g\'en\'eralis\'e de
$\B^\times\simeq\GL_{m_\E}(\D_\E)$, consistu\'e des matrices
monomiales dont les coefficients non nuls sont des puissances 
de $\unif$.
Par la d\'ecomposition de Bruhat, on a
$\B^\times=\U(\BB)\tilde\W\U(\BB)$.

D'apr\`es~\cite[Proposition 1.2]{GSZ}, l'\'el\'ement $w\in\tilde\W$ 
entrelace $\s_{|\U(\BB)}$ si est seulement s'il normalise
$\s_{|\U(\BB)\cap\M}$. 
D'apr\`es la construction de $\M'$, ceci implique que $w\in\M'$,
donc $\I_{\B^\times}(\s_{|\U(\BB)})$ est inclus dans
$\U(\BB)\M'\U(\BB)$.  
En particulier, d'apr\`es  le th\'eor\`eme~\ref{entrelacementvthP},
on a~:
\begin{equation*}
\I_\G(\vartheta_\P)=\J_\P\I_{\B\times}(\s_{|\U(\BB)})\J_\P \subset
\J_\P\M'\J_\P.
\end{equation*}
Par la proposition~\ref{vthP}, et comme $\N\supset\N'$, le couple
$(\J_\P,\vartheta_\P)$ est d\'ecompos\'e au dessus de
$(\J_\P\cap\M',\vartheta_\P{}_{|\J_\P\cap\M'})$. 
La fin de la d\'emonstration est alors identique \`a celle de
\cite[Corollary 3.9(iii)]{BK2}.
\end{proof}

\begin{coro}\label{sdifferent}
La repr\'esentation $\pi$ de $\G$ n'est pas supercuspidale.
\end{coro}

\subsection{}

Dans ce paragraphe, on traite le cas o\`u chacun des $\s_i$
est \'equivalent \`a un conjugu\'e de $\s_1$ sous
$\Gal(k_{\D_\E}/k_\E)$.
En particulier, l'ordre $\BB$ est principal. 
D'apr\`es la proposition~\ref{vthP}, la repr\'esentation
$\vartheta_\P$ determine, \`a permutation circulaire pr\`es, 
le vecteur~:
\begin{equation*}
\V(\J,\vartheta) = \left( \vartheta_1,\ldots,\vartheta_e\right),
\end{equation*}
o\`u chaque $\vartheta_i$ est consid\'er\'e comme une classe
d'\'equivalence d'une repr\'esentation irr\'eductible de
$\J(\b,\AA_i)$.

\begin{prop}\label{contgalconj}
Soit $\tau$ une permutation de $\{1,\ldots,e\}$ et, pour chaque
$1\<i\<e$, soit $\g_i\in\Gal(k_{\D_\E}/k_\E)$.
Alors $\pi$ contient aussi une representation
$\vartheta'=\k\otimes\s'$ de $\J(\b,\AA)$ telle que~:
\begin{equation*}
\V(\J,\vartheta') = \left({}^{\g_1}\vartheta_{\tau(1)},\ldots,
{}^{\g_e}\vartheta_{\tau(e)}\right),
\end{equation*}
o\`u
${}^{\g_i}\vartheta_{\tau(i)}=\k_{\tau(i)}\otimes{}^{\g_i}\s_{\tau(i)}$
pour chaque $1\<i\<e$.
\end{prop}

\begin{proof}
Il suffit de consid\'erer les deux cas particuliers suivants~:
\begin{enumerate}
\item[(i)] 
D'abord, on suppose que
$\V(\J,\vartheta')=
\left({}^{\phi}\vartheta_e,\vartheta_1,\ldots,\vartheta_{e-1}\right)$.
On pose~:
\begin{equation*}
\Pi_\BB=
\begin{pmatrix}&\I_{m_\E-1}\\
\unif&\\
\end{pmatrix}\in\GL_{m_\E}(\D_\E).
\end{equation*}
L'\'el\'ement $\Pi_\BB$ normalise $\BB$ donc, par
(\ref{SoundGrabitz}), il normalise \'egalement $\J(\b,\AA)$ et la
$\b$-extension $\k$.
En particulier, $\pi$ contient la repr\'esentation
$\vartheta'=\vartheta^{\Pi_\BB}=\k\otimes\s^{\Pi_\BB}$, 
et $\s^{\Pi_\BB}$ est \'equivalent \`a 
${}^\phi\s_e\otimes\s_1\otimes\ldots\otimes\s_{e-1}$.
\item[(ii)]
Ensuite, on suppose que~:
\begin{equation*}
\V(\J,\vartheta')= 
\left(\vartheta_1,\ldots,\vartheta_{i-1},\vartheta_{i+1},
\vartheta_i,\vartheta_{i+2},\ldots,\vartheta_e\right)
\end{equation*} 
et $\s_i\not\simeq\s_{i+1}$.
Dans ce cas, la d\'emonstration est identique \`a celle
de~\cite[Proposition 8.3.4]{BK} (voir \emph{loc.\ cit.} p.297), 
quitte \`a remplacer~\cite[Proposition 8.3.5]{BK} par la
proposition~\ref{lieover}.
\end{enumerate}
Ceci termine la preuve de la proposition \ref{contgalconj}.
\end{proof}

D'apr\`es la proposition~\ref{contgalconj}, on peut supposer que
toutes les $\s_i$ sont \'equi\-va\-len\-tes.
La paire $(\J(\b,\AA),\vartheta)$ est donc un type simple au sens 
de~\cite[\S4.1]{VS3}. 

\begin{coro}\label{smeme}
La repr\'esentation $\pi$ est supercuspidale si et
seulement si $\BB$ est un ordre maximal.
\end{coro}

\begin{proof}
D'apr\`es~\cite[Th\'eor\`eme 5.6]{VS3}, le type simple
$(\J(\b,\AA),\vartheta)$ est un type pour une classe inertielle 
$[\M,\rho]_\G$, et le nombre de blocs du sous-groupe de Levi $\M$ 
est \'egal \`a la p\'eriode de $\BB$.
Ainsi $\pi$ est supercuspidale si et seulement si $\M=\G$, \ie si
et seulement si $\BB$ est un ordre maximal.
\end{proof}

\subsection{}

Le th\'eor\`eme suivant est le r\'esultat principal de cette section |
et de cet article.

\begin{theo}
\label{ExhaustionSupercuspidale}
Soit $\pi$ une repr\'esentation irr\'eductible supercuspidale de
niveau non nul de $\G$.
Alors il existe un type simple maximal $(\J,\l)$ tel que
la restriction de $\pi$ \`a $\J$ contienne $\l$.
\end{theo}

\begin{proof}
D'apr\`es le th\'eor\`eme~\ref{scsimple}, la repr\'esentation $\pi$
contient un caract\`ere simple, \ie qu'il existe une strate
simple $[\AA,n,0,\b]$ et un caract\`ere simple $\t\in\Cc(\b,0,\AA)$
tels que la restriction de $\pi$ \`a $\H^1=\H^1(\b,\AA)$ con\-tien\-ne
$\t$. 
D'apr\`es le \S\ref{Scusp}, pour une certaine choix de $\AA$ et de
$\t$, la repr\'esentation $\pi$ contient une repr\'esentation 
$\vartheta=\k\otimes\s$ de $\J=\J(\b,\AA)$, o\`u $\k$ est une
$\b$-extension de $\t$ et $\s$ l'inflation \`a $\J$ d'une
repr\'esentation irr\'eductible cuspidale de 
$\J/\J^1\simeq\U(\BB)/\U_1(\BB)$.  
Comme dans le \S\ref{SvthP}, on a~:
\begin{equation*}
\s_{|\J\cap\M}=\s_1\otimes\ldots\otimes\s_e,
\end{equation*}
o\`u $e$ est la p\'eriode de $\BB$.
Puisque $\pi$ est supercuspidale, les $\s_i$ sont toutes conjugu\'ees
sous $\Gal(k_{\D_\E}/k_\E)$ d'apr\`es le corollaire~\ref{sdifferent}. 
D'apr\`es le corollaire~\ref{smeme}, l'ordre $\BB$ est donc
ma\-xi\-mal et $(\J,\vartheta)$ est un type simple ma\-xi\-mal. 
\end{proof}

\begin{coro}
\begin{itemize}
\item[(i)]
Il existe un prolongement $\L$ de $\l$ \`a $\bar\J=\N_\G(\l)$ tel
que $\pi$ soit \'equivalente \`a l'induite compacte de $\L$ \`a $\G$. 
\item[(ii)]
Le couple $(\J,\l)$ est un type pour la classe inertielle 
$[\G,\pi]_\G$.
\end{itemize}
\end{coro}

\begin{proof}
Il s'agit de \cite[Th\'eor\`eme 5.2]{VS3}.
\end{proof}

Le th\'eor\`eme suivant r\'ecapitule tout le travail effectu\'e.

\begin{theo}
Soit $\mathfrak{s}=[\G_0^r,\pi_0^{\otimes r}]_\G$ une classe
inertielle simple de $\G$, o\`u $r$ est un diviseur de $m$ et 
$\pi_0$ une repr\'esentation irr\'eductible supercuspidale 
de $\G_0=\GL_{m/r}(\D)$.
Il existe un type simple $(\J,\l)$ qui est un type pour
$\mathfrak{s}$. 
\end{theo}

\begin{proof}
Si $\pi_0$ est de niveau z\'ero, il s'agit de \cite[Theorem 5.5]{GSZ}.
Sinon, il s'agit de \cite[Th\'eor\`eme 5.6]{VS3}, conjointement avec le
th\'eor\`eme \ref{ExhaustionSupercuspidale}.
\end{proof}

\providecommand{\bysame}{\leavevmode ---\ }
\providecommand{\og}{``}
\providecommand{\fg}{''}
\providecommand{\smfandname}{\&}
\providecommand{\smfedsname}{\'eds.}
\providecommand{\smfedname}{\'ed.}
\providecommand{\smfmastersthesisname}{M\'emoire}
\providecommand{\smfphdthesisname}{Th\`ese}


\begin{thebibliography}{10}

\bibitem{Br2}
{\scshape P.~Broussous} -- {\og Extension du formalisme de {B}ushnell et
  {K}utzko au cas d'une alg\`ebre \`a division\fg}, \emph{Proc. London Math.
  Soc. (3)} \textbf{77} (1998), no.~2, p.~292--326.

\bibitem{Br3}
{\scshape P.~Broussous {\normalfont \smfandname} B.~Lemaire} -- {\og Building
  of {${\rm GL}(m,D)$} and centralizers\fg}, \emph{Transform. Groups}
  \textbf{7} (2002), no.~1, p.~15--50.

\bibitem{Br1}
{\scshape P.~Broussous} -- {\og Hereditary orders and embeddings of local
  fields in simple algebras\fg}, \emph{J. Algebra} \textbf{204} (1998), no.~1,
  p.~324--336.

\bibitem{Br4}
\bysame , {\og Minimal strata for {${\rm GL}(m,D)$}\fg}, \emph{J. Reine Angew.
  Math.} \textbf{514} (1999), p.~199--236.

\bibitem{BG}
{\scshape P.~Broussous {\normalfont \smfandname} M.~Grabitz} -- {\og Pure
  elements and intertwining classes of simple strata in local central simple
  algebras\fg}, \emph{Comm. Algebra} \textbf{28} (2000), no.~11, p.~5405--5442.

\bibitem{BH1}
{\scshape C.~J. Bushnell {\normalfont \smfandname} G.~Henniart} -- {\og Local
  tame lifting for {${\rm GL}(N)$}. {I}. {S}imple characters\fg}, \emph{Inst.
  Hautes \'Etudes Sci. Publ. Math.} (1996), no.~83, p.~105--233.

\bibitem{BK}
{\scshape C.~J. Bushnell {\normalfont \smfandname} P.~C. Kutzko} -- \emph{The
  admissible dual of {${\rm GL}({\rm N})$} via compact open subgroups},
  Princeton University Press, Princeton, NJ, 1993.

\bibitem{BK1}
\bysame , {\og Smooth representations of reductive $p$-adic groups: structure
  theory via types\fg}, \emph{Proc. London Math. Soc. (3)} \textbf{77} (1998),
  no.~3, p.~582--634.

\bibitem{BK2}
\bysame , {\og Semisimple types in {${\rm GL}\sb n$}\fg}, \emph{Compositio
  Math.} \textbf{119} (1999), no.~1, p.~53--97.

\bibitem{Grabitz}
{\scshape M.~Grabitz} -- {\og Simple characters for principal orders and their
  matching\fg}, \emph{{\rm Preprint MPIM1999-117, Max Planck Institut, Bonn}}
  (1999).

\bibitem{Grabitz2}
\bysame , {\og Simple characters for principal orders, part {II}\fg},
  \emph{{\rm Preprint MPIM2003-56, Max Planck Institut, Bonn}} (2003).

\bibitem{GSZ}
{\scshape M.~Grabitz, A.~J. Silberger {\normalfont \smfandname} E.-W. Zink} --
  {\og Level zero types and {H}ecke algebras for local central simple
  algebras\fg}, \emph{J. Number Theory} \textbf{91} (2001), no.~1, p.~92--125.

\bibitem{HM}
{\scshape R.~Howe {\normalfont \smfandname} A.~Moy} -- {\og Minimal {$K$}-types
  for {${\rm GL}\sb n$} over a {$p$}-adic field\fg}, \emph{Ast\'erisque}
  (1989), no.~171-172, p.~257--273, Orbites unipotentes et repr\'esentations,
  II.

\bibitem{VS1}
{\scshape V.~S{\'e}cherre} -- {\og Repr\'esentations lisses de {${\rm
  GL}(m,D)$}, {I} : caract\`eres simples\fg}, \emph{Bull. Soc. math. France}
  \textbf{132} (2004), no.~3, p.~327--396.

\bibitem{VS2}
\bysame , {\og Repr\'esentations lisses de {${\rm GL}(m,D)$}, {II} :
  {$\beta$}-extensions\fg}, \emph{Compositio Math.} \textbf{141} (2005),
  p.~1531--1550.

\bibitem{VS3}
\bysame , {\og Repr\'esentations lisses de {${\rm GL}(m,D)$}, {III} : {types
  simples}\fg}, \emph{Ann. Scient. \'Ec. Norm. Sup.} \textbf{38} (2005),
  p.~951--977.

\bibitem{St1}
{\scshape S.~Stevens} -- {\og Double coset decompositions and intertwining\fg},
  \emph{Manu\-scrip\-ta Math.} \textbf{106} (2001), no.~3, p.~349--364.

\bibitem{St3}
\bysame , {\og Semisimple strata for {$p$}-adic classical groups\fg},
  \emph{Ann. Sci. \'Ecole Norm. Sup. (4)} \textbf{35} (2002), no.~3,
  p.~423--435.

\bibitem{St4}
\bysame , {\og Semisimple characters for {$p$}-adic classical groups\fg},
  \emph{Duke Math. J.} \textbf{127} (2005), no.~1, p.~123--173.

\bibitem{Zi5}
{\scshape E.-W. Zink} -- {\og Representation theory of local division
  algebras\fg}, \emph{J. Reine Angew. Math.} \textbf{428} (1992), p.~1--44.

\bibitem{Zi1}
\bysame , {\og More on embeddings of local fields in simple algebras\fg},
  \emph{J. Number Theory} \textbf{77} (1999), no.~1, p.~51--61.

\end{thebibliography}
\end{document}